\documentclass{amsart}
\usepackage[latin1]{inputenc}
\usepackage[english]{babel}
\usepackage{csquotes}

\usepackage{lipsum}  

\usepackage{todonotes}

\usepackage[title]{appendix}

\usepackage{color} 
\usepackage[colorlinks=true,linktoc=all,linkcolor=blue,citecolor=red,filecolor=pink,urlcolor=blue]{hyperref}
\usepackage[hyperref=true,style=alphabetic,citestyle=alphabetic,backend=bibtex,maxnames=100,doi=false,isbn=false,url=false,giveninits=true]{biblatex}
\bibliography{biblio}

\usepackage{wasysym}
\usepackage[hypcap=false]{caption}
\usepackage{amssymb,amsmath,latexsym,graphicx}
\usepackage{enumitem}
\usepackage{subfigure}
\usepackage{tikz}
\usepackage{tikz-cd}
\usepackage{hyperref}
\usetikzlibrary{arrows}

\theoremstyle{plain}                    
\newtheorem{theorem}{Theorem}[subsection]
\newtheorem{lemma}[theorem]{Lemma}
\newtheorem{proposition}[theorem]{Proposition}
\newtheorem{corollary}[theorem]{Corollary}
\newtheorem*{graftconj}{Grafting Conjecture}

\newtheorem{introtheorem}{Theorem}

\newcommand{\theoremnumber}{} % initialize
\newtheorem*{maintheorem}{Theorem \theoremnumber}
\newenvironment{maintheoremc}[1]
  {\renewcommand{\theoremnumber}{#1}%
   \begin{maintheorem}}
  {\end{maintheorem}}

\theoremstyle{definition}

\newtheorem{example}[theorem]{Example}

\theoremstyle{definition}
\newtheorem{remark}[theorem]{Remark}

\numberwithin{equation}{subsection} 

\newcommand{\CC}{\mathbb{C}}
\newcommand{\DD}{\mathbb{D}}
\newcommand{\FF}{\mathbb{F}}
\newcommand{\HH}{\mathbb{H}}
\newcommand{\NN}{\mathbb{N}}
\newcommand{\PP}{\mathbb{P}}
\newcommand{\RR}{\mathbb{R}}
\renewcommand{\SS}{\mathbb{S}}
\newcommand{\ZZ}{\mathbb{Z}}

\newcommand{\C}{\mathcal{C}}

\renewcommand{\P}{\mathcal{P}}
\newcommand{\T}{\mathcal{T}}

\newcommand{\R}{\mathcal{R}}

\newcommand{\Circles}{\mathfrak{C}}

\newcommand{\Tri}{\mathfrak{T}}
\newcommand{\wTri}{\widetilde{\Tri}}

\newcommand{\pslc}{\mathrm{PSL}_2\CC}
\newcommand{\slc}{\mathrm{SL}_2\CC}
\newcommand{\pslr}{\mathrm{PSL}_2\RR}

\newcommand{\psu}{\mathrm{PSU}(2)}

\newcommand{\cp}{\CC\PP^1}
\newcommand{\rp}{\RR\PP^1}
\newcommand{\End}{E}
\newcommand{\Ends}{\mathcal{E}}

\newcommand{\tri}{\mathrm{T}}

\newcommand{\Char}{\R}

\newcommand{\Chart}{\Char^\odot}
\newcommand{\Moduli}{\P}
\newcommand{\Modulic}{\P^\bullet}
\newcommand{\Modulit}{\P^\odot}
\newcommand{\Teich}{\T}

\newcommand{\Sigmac}{\overline{\Sigma}}
\newcommand{\Sigmaw}{\widetilde{\Sigma}}
\newcommand{\Sigmawe}{\widetilde{\Sigma}^{\#}}

\newcommand{\Sw}{\widetilde{S}}
\newcommand{\Swe}{\widetilde{S}^{\#}}

\newcommand{\hol}{\rho}
\newcommand{\sch}{\mathcal S}

\newcommand{\cangle}[2]{\angle_{#1} {#2}}

\renewcommand{\Re}{\text{Re}}
\renewcommand{\Im}{\text{Im}}

\newcommand{\pDisk}{\DD^*}

\newcommand{\maxN}{\hat{N}}

\newcommand{\Quadddp}[1]{\mathcal Q_2({#1})}

\newcommand{\Compl}[2]{\mathrm{M}^{#1}({#2})}
\newcommand{\Binf}[2]{\partial^{#1}_{\infty}({#2})}

\DeclareMathOperator{\Clo}{Cl}
\DeclareMathOperator{\Clot}{Cl_{\#}}
\newcommand{\Cloc}[1]{\Clo_{#1}}
\DeclareMathOperator{\dev}{dev}
\newcommand{\devt}{\dev_{\#}}
\newcommand{\devc}[1]{\dev_{#1}}
\DeclareMathOperator{\Hol}{Hol}
\DeclareMathOperator{\Ind}{I}
\DeclareMathOperator{\Int}{Int}

\DeclareMathOperator{\Rot}{Rot}
\DeclareMathOperator{\tr}{tr}
\DeclareMathOperator{\Gr}{Gr}
\DeclareMathOperator{\Ball}{B}
\DeclareMathOperator{\Ballt}{\Ball_{\#}}
\newcommand{\Ballc}[1]{\Ball_{#1}}

\newcommand{\critrad}[1]{\rho_{#1}}

\def\n@te#1{\textsf{\boldmath \textbf{#1}}\leavevmode}
\def\Note#1{\textcolor{red}{\n@te {}}}

\newcommand\restr[2]{{% we make the whole thing an ordinary symbol
  \left.\kern-\nulldelimiterspace % automatically resize the bar with \right
  #1 % the function
  \vphantom{\big|} % pretend it's a little taller at normal size
  \right|_{#2} % this is the delimiter
  }}

\begin{document}

\title[$\mathbb{CP}^1$--structures on the thrice--punctured sphere]{Tame and relatively elliptic $\mathbb{CP}^1$--structures\\ on the thrice--punctured sphere}

\author{Samuel A. Ballas}
\address{Department of Mathematics, Florida State University, 1017 Academic Way, Tallahassee, FL 32304, USA.}
\email{ballas@math.fsu.edu}

\author{Philip L. Bowers}
\address{Department of Mathematics, Florida State University, 1017 Academic Way, Tallahassee, FL 32304, USA.}
\email{bowers@math.fsu.edu}

\author{Alex Casella}
\address{Department of Mathematics, Florida State University, 1017 Academic Way, Tallahassee, FL 32304, USA.}
\email{alex.casella.usyd@gmail.com}

\author{Lorenzo Ruffoni}
\address{Department of Mathematics, Tufts University, 177 College Ave, Medford, MA 02155, USA.}
\email{lorenzo.ruffoni2@gmail.com}

% \thanks{}
\subjclass[2010]{57M50, 30F30}
\date{\today}

% \dedicatory{}

\begin{abstract}
Suppose a relatively elliptic representation $\rho$ of the fundamental group of the thrice--punctured sphere $S$ is given. We prove that all projective structures on $S$ with holonomy $\rho$ and satisfying a tameness condition at the punctures can be obtained by grafting certain circular triangles. The specific collection of triangles is determined by a natural framing of $\rho$.
In the process, we show that (on a general surface $\Sigma$ of negative Euler characteristics) structures satisfying these conditions can be characterized in terms of their M\"obius completion, and in terms of certain meromorphic quadratic differentials. 
\end{abstract}

\maketitle
\tableofcontents

%%%%%%% prints list of todos	
% \makeatletter
% \providecommand\@dotsep{5}
% \makeatother
% \listoftodos\relax

	%%%%%%%%%%%%%%%%%%%%%%%%%%%%%%%%%
	%%%%%%%%%%%%%%%%%%%%%%%%%%%%%%%%%
	%%%%%%%%%%%%%%%%%%%%%%%%%%%%%%%%%

\section{Introduction}\label{sec:intro}
\addtocontents{toc}{\protect\setcounter{tocdepth}{1}}

This paper deals with the geometry of surfaces which are locally modelled on the geometry of the Riemann sphere $\cp$, and their grafting deformations. 
Throughout the paper, $\Sigma$ denotes an orientable surface with finitely many punctures (and no boundary) and $\Sigmac$ denotes the closed orientable surface where the punctures have been filled in. 
While the main technical core of the paper holds for a general $\Sigma$ with negative Euler-characteristic (see \S \ref{sec:tame_and_rel_elliptic_structures} and \S\ref{sec:complex}), the final chapter \S\ref{sec:thrice_punct_sphere} deals specifically with the case of a thrice--punctured sphere, which we denote by $S$. 

The structures under consideration here are known as complex projective structures or $(\pslc,\cp)$--structures. 
We denote respectively by $\Teich (\Sigma)$ and $\Moduli(\Sigma)$  the deformation spaces of complex and complex projective structures on $\Sigma$. 
We also denote by $\Char(\Sigma)$ the space of representations of $\pi_1(\Sigma)$ into $\pslc$, up to conjugation by $\pslc$. We have natural forgetful maps
$$
\pi:\Moduli (\Sigma)\to \Teich(\Sigma) \quad \text{and} \quad \Hol:\Moduli (\Sigma)\to \Char(\Sigma),
$$
respectively recording the underlying complex structure and holonomy representation.
We refer the reader to \S \ref{sec:tame_and_rel_elliptic_structures} for precise definitions, and to  \cite{DU} for a general survey about $\cp$-structures. For more on the geometry of the deformation space, see \cite{FA20}.

Classic examples of complex projective structures are given by hyperbolic metrics (seen as $(\pslr,\HH^2)$--structures), but a general projective structure is not defined by a Riemannian metric, nor is it completely determined by its holonomy (not even in the Fuchsian case, see for instance \cite{GO87,CDF}). 
However, under some additional conditions $\Hol$ is known to be a local homeomorphism (see \cite{HE75,L93,GM19}), i.e. a structure is at least locally determined by its holonomy.
A major question in the field is the description in geometric terms of all structures having the same holonomy.

\begin{graftconj}[Problems 12.1.1-2 in \cite{GKM00}]
Two complex projective structures have the same holonomy if and only if it is possible to obtain one from the other by some sequence of graftings and degraftings.
\end{graftconj}

Here \textit{grafting} refers to a geometric surgery on $\Sigma$ which consists in cutting $\Sigma$ open along a curve and inserting a domain from $\cp$, and \textit{degrafting} is the inverse operation. 
For the reader familiar with grafting deformations: by grafting we will always mean \textit{projective $2\pi$--grafting}. 
This construction allows one to change a structure without changing its holonomy, and iterating this construction shows that $\Hol$ has infinite fibers. 
The Grafting Conjecture has been verified for closed surfaces: the case of (quasi-)Fuchsian representations is due to Goldman \cite{GO87}, and Baba has addressed the case of generic (i.e. totally loxodromic) representations in a series of papers \cite{BA10,BA12,BA15,BA17}.

Inspired by a specific question about punctured spheres in \cite[Problem 12.2.1]{GKM00}, we propose a study of certain structures on the thrice--punctured sphere, and we prove the Grafting Conjecture in this setting (see part  \ref{intro:other work on grafting} of this introduction for a comparison with related results available in the literature). 
It is worth noticing that the complex projective geometry around a puncture is much more interesting than the underlying complex geometry. As an example, consider the two structures on the thrice-punctured sphere given by the complete hyperbolic metric of finite area and by the inclusion $\cp \setminus \{0,1,\infty\} \subseteq \cp$: they are not isomorphic as complex projective structures, but they have the same underlying complex structure.

The study of holonomy fibers also has an analytic motivation coming from the classical monodromy problems for ODEs, i.e. generalization of Hilbert's XXI problem.
Since the work of Poincar\'e \cite{P08}, projective structures have been known as a geometric counterpart to second--order linear ODEs. In more recent years, some monodromy problems for such ODEs have successfully been approached in terms of holonomy problems for projective structures (see \cite{GKM00,CDHL,G19,GM19b,KA20,CFG22}). 

\vspace{.5cm}
We consider structures satisfying some regularity conditions at the punctures, which can be roughly stated as follows (see \S \ref{subsec:tame_rel_ell} for precise definitions):
\begin{itemize}[leftmargin=0.5cm]
    \item (\textit{tameness}) each local chart has a limit along arcs going off into a puncture; 
    \item (\textit{relative ellipticity}) each peripheral holonomy (i.e. the holonomy around each puncture) is a non-trivial elliptic element in $\pslc$;
    \item (\textit{non-degeneracy}): there is no pair of points $p_\pm \in \cp$ such that the entire holonomy preserves the set $\{p_\pm\}$.
\end{itemize}
Motivating examples of tame structures arise from the study of triangle groups and automorphism groups (as in \cite[Remark 2.13]{FaR19}), and more generally from metrics of constant curvature with cones or cusps.
Tameness is not a generic condition in the space of all complex projective structures, but is a natural case to consider. 
Indeed, it corresponds to the condition that the associated second--order linear ODE has \textit{regular singular points} (see Theorem~\ref{intro:quadratic differentials} below).
It turns out that the peripheral holonomy of a tame structure can only be trivial, parabolic, or elliptic (see Lemma~\ref{lem:devt_is_equivariant}), so the second condition is a generic condition within the space of holonomies of tame structures. 
In particular, it implies that there are no \textit{apparent singularities} (i.e. no puncture has trivial holonomy).

For an arbitrary surface $\Sigma$, we denote by $\Modulit(\Sigma)$ the subspace of $\Moduli(\Sigma)$ consisting of non-degenerate tame and relatively elliptic structures: the white disk in the superscript represents the local invariance under a rotation, and the black dot the possibility to extend the charts to the puncture.
The tameness condition provides a natural choice of a fix point for each peripheral holonomy, i.e. a \textit{framing} for the holonomy representation (see Corollary \ref{cor:framed_hol}). 
We observe that grafting preserves this natural framing, which suggests a more precise formulation of the Grafting Conjecture in the non-compact case.
Our main result in the case of the thrice--punctured sphere $S$ is the following, which confirms the conjecture, in the spirit of Problem 12.2.1 in \cite{GKM00}.

\begin{introtheorem}\label{intro:solution to grafting conjecture}
Two structures in $\Modulit(S)$ have the same framed holonomy if and only if it is possible to obtain one from the other by some combination of graftings and degraftings along ideal arcs.
\end{introtheorem}

Here an arc is \textit{ideal} if it starts and ends at a puncture. 
To the best of our knowledge this is the first result in this direction for the case of non-compact surfaces with non--trivial holonomy around the punctures. 
 
The representations involved here are representations of the free group $\FF_2=\pi_1(S)$ generated by elliptic elements. 
Representations satisfying certain rationality conditions correspond to the classical triangle groups, but the general ones are non-discrete. 
In all cases we construct an explicit list of triangular membranes (i.e. immersions of a triangle in $\cp$) realizing these representations, and identify the ones that are \textit{atomic}: these can be taken as basic building blocks that can be grafted to reconstruct all the projective structures in $\Modulit(S)$. Theorem~\ref{intro:solution to grafting conjecture} is a consequence of the following theorem.

\begin{introtheorem}\label{intro:main theorem}
Every $\sigma \in \Modulit(S)$ is obtained by grafting on an atomic triangular structure with the same framed holonomy.
\end{introtheorem}

Another consequence of Theorem~\ref{intro:main theorem} is a handy description of the moduli space $\Modulit(S)$ with positive real coordinates, which we plan to address in a future work.

When a representation $\rho:\pi_1(S)\to \pslc$ is unitary (i.e. is conjugate into $\psu$), it preserves a spherical metric, and a structure $\sigma \in \Modulit(S)$ is given by a spherical metrics with cone points. This special case of Theorem~\ref{intro:main theorem} is implicit in the proof of  \cite[Theorem 3.8]{MP16}, which constructs such spherical metrics by gluing together spherical triangles and bigons. Grafting a spherical metric results in a spherical metric, with increased angles at the cones. However in general this is not always the case: for example the structure obtained by grafting a hyperbolic structure is not defined by any Riemannian metric.

\vspace{.5cm}
While our results about the Grafting Conjecture are for the case of the thrice--punctured sphere $S$, the main technical core of the paper applies to any non-compact surface $\Sigma$ of negative Euler characteristic, and is of independent interest. It consists of a characterization of structures from $\Modulit(\Sigma)$ in terms of their M\"obius completion (see \S\ref{sec:tame_and_rel_elliptic_structures} and \cite{KP94}) and in terms of meromorphic projective structures (see \S \ref{sec:complex} and \cite{AB18}).
The easy case of structures on a twice-punctured sphere can be worked out concretely, see Remark~\ref{rem:twice}.
In the remaining part of the introduction we present our main results in the general case (see \S \ref{intro:results on general surfaces}), as well as a comparison with other work in the literature about the Grafting Conjecture (see \S \ref{intro:other work on grafting}).

%%%%%%%%%%%%%%%%%%%%%%%%%%%%%%%%%%%%%%%%%%%%%%%%%%%%%%%%%%%%%%

\subsection{Results for general surfaces}\label{intro:results on general surfaces}
The universal cover $\Sigmaw$ of $\Sigma$ is a topological disk. It admits a natural decoration obtained by adding ideal points at infinity ``above'' the punctures. We call these ideal points \textit{ends}. This gives rise to a natural enlargement of $\Sigmaw$ that we call the \textit{end-extension}, and denote by $\Sigmawe$. Part of the paper is concerned with understanding the behavior of the developing map in the limit to an end.

\subsubsection*{M\"obius completion}
Any complex projective structure $\sigma$ on $\Sigma$ can be used to define another natural extension of $\Sigmaw$, known as the \textit{M\"obius completion} $\Compl{\sigma}{\Sigmaw}$, which comes with a (non-canonical) structure of complete metric space (see \cite{KP94}).
For instance, when $\sigma$ is induced by a spherical metric with cone points, $\Compl{\sigma}{\Sigmaw}$ coincides with $\Sigmawe$, while when $\sigma$ is induced by a complete hyperbolic metric of finite area $\Compl{\sigma}{\Sigmaw}$ identifies with the closed disk model for the hyperbolic plane $\mathbb H^2 \cup \rp$ (see Examples \ref{ex:completion_hyperbolic} and \ref{ex:completion_spherical}). 

The topologies on $\Sigmawe$ and on $\Compl{\sigma}{\Sigmaw}$ are not in general compatible.
One of the main technical contributions of this paper is a study of the geometry of the M\"obius completion  $\Compl{\sigma}{\Sigmaw}$ for $\sigma \in \Modulit(\Sigma)$, and of its relation with the end--extension  $\Sigmawe$ (see \S \ref{sec:tame_and_rel_elliptic_structures}). 
Tameness of a structure $\sigma$ implies that its developing map  admits  natural continuous extensions $\devt$ to the end-extension $\Sigmawe$ and $\devc{\sigma}$ to the M\"obius completion $\Compl{\sigma}{\Sigmaw}$. We study the local properties of $\devt$ and $\devc{\sigma}$ around the ends.

\begin{introtheorem}\label{intro:embedding}
Let $\sigma \in \Moduli (\Sigma)$ be non-degenerate and without apparent singularities. Let $j^{\#}:\Sigmaw \to \Sigmawe$ and $j_\sigma:\Sigmaw \to \Compl{\sigma}{\Sigmaw}$ be the  natural embeddings. Then $\sigma \in \Modulit(\Sigma)$ if and only if there exists a continuous open $\pi_1(\Sigma)$--equivariant embedding $j^{\#}_\sigma:\Sigmawe \to \Compl{\sigma}{\Sigmaw}$
that makes the following diagram commute.

\centering
\begin{tikzpicture}[scale=1]
\node (S) at (-2,0) {$\Sigmaw$};
\node (CP1) at (2,0) {$\cp$};
\node (E) at (0,1) {$\Sigmawe$ };
\node (M) at (0,-1) {$\Compl{\sigma}{\Sigmaw}$};
\path[->,font=\scriptsize,>=angle 90]
(S) edge node[above]{$j^{\#}$}    (E)
(S) edge node[below]{$j_\sigma$}    (M)
(E) edge node[above]{$\devt$}    (CP1)
(M) edge node[below]{$\devc{\sigma}$}    (CP1)
(E) edge node[left]{$j^{\#}_\sigma$}    (M);
\end{tikzpicture}
\end{introtheorem}

In this statement, continuity is a consequence of tameness of $\sigma$, and openness is a consequence of relative ellipticity. 

In general, the developing map for a projective structure is a surjection onto $\cp$, in which case it fails to be a global covering map. However, under certain circumstances it is known to be a covering map onto a component of the domain of discontinuity in $\cp$ for its holonomy representation (see for instance \cite[Theorem 1]{K7169}). But in general the holonomy group is not discrete, so it has no domain of discontinuity. The following statement shows that in our context some local covering behavior can be guaranteed around ends.

\begin{introtheorem}\label{intro:max neighborhood}
Let $\sigma \in \Modulit(\Sigma)$, and let $\End$ be an end. Then there is a neighborhood $\maxN_{\End}$ of $\End$ in $\Compl{\sigma}{\Sigmaw}$ onto which the developing map for $\sigma$ restricts to a branched covering map, branching only at $\End$, and with image a round disk in $\cp$.
\end{introtheorem}

These neighborhoods should be regarded as an analogue of the round balls considered in \cite{KP94}, but ``centered'' at ideal points in the M\"obius completion. While Theorem~\ref{intro:max neighborhood} is stated as a local fact, we actually show that such a neighborhood can be chosen to be so large as to have another ideal point on its boundary.
We use the existence of these neighborhoods to define a local geometric invariant, which we call the \textit{index} (see \S \ref{subsec:index_puncture}). This number measures the angle described by the developing map at a puncture, and provides a notion of complexity for an inductive proof of Theorem~\ref{intro:main theorem}.

\subsubsection*{Meromorphic projective structures}
A second major ingredient (once again valid for an arbitrary non-compact surface $\Sigma$) consists of an analytic description of structures in $\Modulit(\Sigma)$ as meromorphic projective structures in the sense of \cite{AB18}. These are projective structures whose developing map is defined by solving certain differential equations with coefficients given by meromorphic quadratic differentials on the closed surface $\Sigmac$ (with poles corresponding to the punctures of $\Sigma$, see \S \ref{subsec:local_theory} for precise definitions).
The local control from Theorem~\ref{intro:max neighborhood} allows us to obtain the following result.

\begin{introtheorem}\label{intro:quadratic differentials}
Let  $\sigma \in \Moduli(\Sigma)$ and let $X \in \Teich(\Sigma)$ be the underlying complex structure. Then $\sigma \in \Modulit(\Sigma)$ if and only if $X$ is a punctured Riemann surface and $\sigma$ is represented by a meromorphic quadratic differential on $X$ with double poles and reduced exponents in $\RR\setminus \ZZ$.
\end{introtheorem}

Here the parametrization of projective structures by quadratic differentials is the classical one in terms of the Schwarz derivative, which here is taken with respect to any compatible holomorphic structure on the closed Riemann surface obtained by filling the punctures (e.g. the constant curvature uniformization).
From this point of view, the index of a structure at a puncture corresponds to the absolute value of the exponents of the quadratic differential, so it can be computed in terms of its residues.

It should also be noted that work of Luo in \cite{L93} guarantees that $\Hol$ is a local homeomorphism for this class of meromorphic projective structures, as there are no apparent singularities. Therefore fibers of $\Hol$ in $\Modulit(\Sigma)$ are discrete, and in particular it makes sense to seek a description of them in terms of a discrete geometric surgery such as the type of grafting that we consider in this paper.

\subsubsection*{Outline of the proof of Theorem~\ref{intro:main theorem}}
% \proof{}
    Let $S$ be the thrice--punctured sphere, and let $\sigma \in \Modulit(S)$, with developing map $\dev$ and holonomy $\rho$.
    By Theorems~\ref{intro:embedding} and \ref{intro:max neighborhood}, $\dev$ extends continuously and equivariantly to the ends, and restricts to a branched covering map on a suitable neighborhood of each end. This allows to define the index of $\sigma$ at each puncture. Then we construct a circular triangle such that the pillowcase obtained by doubling it provides a structure $\sigma_0\in \Modulit(S)$ with holonomy $\rho$. Note that such a triangle is not unique in general. A careful analysis of the framing of $\rho$ defined by $\sigma$  shows that such a triangle can be found with the same framing for $\rho$. On such a triangle, we find a suitable combination of disjoint ideal arcs that are graftable, and we show that if sufficiently many grafting regions are inserted, the resulting structure $\sigma'\in \Modulit(S)$ has the same indices as $\sigma$. By Theorem~\ref{intro:quadratic differentials}, $\sigma$ and $\sigma'$ can be represented by two meromorphic differentials on the Riemann sphere $\cp$ with double poles at $0,1,\infty$. 
    Two such differentials on $\cp$ are completely determined by their residues, and in this case residues can be computed directly from the indices, hence are the same.
    So we conclude that $\sigma = \sigma'$.
      \hfill $\square$    
    
%%%%%%%%%%%%%%%%%%%%%%%%%%%%%%%%%%%%%%%%%%%%%%%%%%%%%%%%%%%%%%

\subsection{Relation to other work about the Grafting Conjecture}\label{intro:other work on grafting}
Following seminal work of Thurston (see \cite{KT92,DU,BA19} and references therein), grafting (in its general version) has been successfully used as a tool to explore the deformation space of $\cp$--structures.
The grafting we consider here preserves the holonomy representation, hence can be used to explore holonomy fibers. The classical case is that of structures on a closed surface with Fuchsian holonomy, which was considered by Goldman (see \cite{GO87}). 
Our work displays some technical differences, that we summarize here for the expert reader.

\subsubsection*{Framing} The main results for closed surfaces in \cite{GO87,BA12,BA15,BA17,CDF2} confirm the Grafting Conjecture, i.e. that two structures with the same holonomy differ by grafting.
In our non-compact case there is a natural framing for the holonomy which needs to be taken into consideration, as it is preserved by grafting (see Lemma~\ref{lem:grafting_properties}). We prove that having the same framed holonomy is not only necessary, but also sufficient, for two structures on the thrice--punctured sphere to differ by grafting.

\subsubsection*{Basepoints for holonomy fibers} When $\rho:\pi_1(\Sigma)\to \pslc$ is Fuchsian, the holonomy fiber $\Hol^{-1}(\rho)$ contains a preferred structure, namely the hyperbolic structure $\mathbb H^2 / \rho (\pi_1(\Sigma))$. This structure serves as a basepoint, i.e. any other structure in $\Hol^{-1}(\rho)$ can be obtained by grafting it (see \cite{GO87}). In this paper, we show that every representation coming from $\Modulit(S)$ is generated by reflections in the sides of a circular triangle in $\cp$. Even when such a representation $\rho$ is non-discrete, the pillowcase obtained by doubling the triangle provides a basepoint in the holonomy fiber $\Hol^{-1}(\rho)$. 
A first guess is that every structure in $\Hol^{-1}(\rho)$ is obtained by grafting this pillowcase. However, this is not the case, because of the aforementioned framing, which is given by the vertices of the triangle.
In \S\ref{subsec:tri_immersions} we identify the list of the structures that can be taken as basepoints in the above sense, which we call atomic. Interestingly, they are not all embedded geodesic triangles for some invariant metric.

\subsubsection*{Type of grafting curves} In the classical Fuchsian case it is enough to perform grafting along simple closed geodesics on the hyperbolic basepoint (see \cite{GO87}).
Here we consider grafting along ideal arcs, i.e. arcs that start and end at punctures. Grafting along open arcs is also known as bubbling in the literature (see \cite{GKM00,CDF,R19a,R19b,FR19}). 
Most structures considered here are not metric, but they still have a well-defined notion of circular arc. We show that in most cases grafting arcs can be chosen to be circular. % (see Figures \ref{fig:multi_grafting} and \ref{fig:edge_graft_degenerate_hyp}). 
    
\subsubsection*{Uniqueness of grafting curves} In the classical Fuchsian case grafting curves are homotopically non-trivial, and are uniquely determined by the structure itself (see \cite{GO87}). Here grafting regions do not carry any topology (they are disks), hence they should not be expected to be canonically associated with the structure. Indeed it is quite common for a structure to arise from different graftings on different atomic structures.

	%%%%%%%%%%%%%%%%%%%%%%%%%%%%%%%%%%%%%%%%%%%%%%%%%%%%%%%%%%%%%%

	\subsection*{Outline of the paper}
	
    Section \ref{sec:cp1} contains background material about the geometry of circles and circular triangles in $\cp$ (see \S\ref{subsec:config_circles} and \S\ref{subsec:elliptic_mobius_transf}). 
    In \S\ref{subsec:tri_immersions} we provide a classification of certain triangular immersions that will serve as the atomic structures for our main grafting results. This classification is referred to in different parts of the paper, and it is summarized in Tables  \ref{tab:hyperbolic}, \ref{tab:spherical} and \ref{tab:Euclidean}. 
    
    Section \ref{sec:tame_and_rel_elliptic_structures} introduces the main geometric definitions, i.e. that of tameness and relative ellipticity. In \S\ref{subsec:mobius_completion} we study the geometry of the M\"obius completion for a general surface and address Theorem~\ref{intro:embedding}. The proof of Theorem~\ref{intro:max neighborhood} is in \S\ref{subsec:canonic_max_nbhd}, where we show that the developing map restricts to a nice branched cover around each end. This is used in \S\ref{subsec:index_puncture} to define the index of a puncture, and in Section \S\ref{sec:complex} to obtain a characterization of tame and relatively elliptic structures in terms of quadratic differentials on a general Riemann surface. In particular we show that the geometric notion of index can be also defined and computed analytically. Theorem~\ref{intro:quadratic differentials} is contained in \S\ref{subsec:meromorphic_structures}. 
    
    Finally, in Section \ref{sec:thrice_punct_sphere} we restrict our attention to the case of the thrice--punctured sphere $S$. 
    In \S\ref{subsec:tri_structures} we define the class of triangular structures on $S$, based on \S \ref{subsec:tri_immersions}, and in \S\ref{subsec:grafting theorem} we prove the main grafting results of Theorems~\ref{intro:solution to grafting conjecture} and~\ref{intro:main theorem}.

	\subsection*{Acknowledgements}
    We thank Gabriele Mondello for some useful conversations, Spandan Ghosh and Subhojoy Gupta for helpful comments on an earlier version of Lemma~\ref{lem:local_theory}, and the anonymous reviewers for their careful reading of the manuscript and their insightful suggestions.

	%%%%%%%%%%%%%%%%%%%%%%%%%%%%%%%%%
	%%%%%%%%%%%%%%%%%%%%%%%%%%%%%%%%%
	%%%%%%%%%%%%%%%%%%%%%%%%%%%%%%%%%

	\section{Basics on complex projective geometry}\label{sec:cp1}
	\addtocontents{toc}{\protect\setcounter{tocdepth}{3}}

	In this chapter we collect some background about the geometry of the Riemann sphere, on which our geometric structures will be modelled, mainly to fix notation and terminology.
	Let $\cp$ denote the the set of complex lines through the origin in $\CC^2$, i.e. the quotient of $\CC^2 \setminus \{0\}$ by scalar complex multiplication. We fix identifications of $\cp$ with the extended complex plane $\CC \cup \{\infty\}$ and the unit sphere $\SS^2$. Through them, $\cp$ inherits a natural complex structure, an orientation, and a spherical metric. A \textit{circle} in $\cp$ is a circle or a line in $\CC \cup \{\infty\}$. Every circle divides $\cp$ into two disks, each of which has a standard identification with the hyperbolic plane which respects the underlying complex structure.
	We denote by $\pslc$ the group of projective classes of $2$--by--$2$ complex matrices of determinant 1. This group acts on $\cp$ by \textit{M\"obius transformations}:
	
	$$
	\pslc \times \cp \to \cp, \begin{bmatrix}
	a & b \\
	c & d
	\end{bmatrix}, z\mapsto \frac{az+b}{cz+d}.
	$$
	
	For elements in $\pslc$, traces and determinants are not well defined. However there is a two-to-one map $\slc \rightarrow \pslc$ such that $\pm A \mapsto [A]$. Therefore, given an element $G \in \pslc$, we can always assume it to be in $\slc$ modulo a sign. It follows that $\det(G)$, $|\tr(G)|$ and $\tr(G)^2$ are well defined quantities. The action of $\pslc$ on $\cp$ is faithful, and simply transitive on triples of pairwise distinct points. In particular, we can always map three distinct points $(p_1,p_2,p_3)$ to $(0,1,\infty)$. M\"obius transformations are conformal, preserve cross ratios and preserve circles. Three distinct points in $\cp$ determine a unique circle through them. Great circles are geodesic circles in the underlying spherical metric. However, elements of $\pslc$ are generally not isometries, and so the set of great circles is not $\pslc$--invariant. \par
	
	A non-trivial element $G \in \pslc$ is classified as follows:
	\begin{itemize}
		\item \emph{Parabolic} if $\tr(G)^2 = 4$.
		\item \emph{Elliptic} if $\tr(G)^2$ is real and $\tr(G)^2 < 4$.
		\item \emph{Loxodromic} otherwise.
	\end{itemize}
	
	%%%%%%%%%%%%%%%%%%%%%%%%%%%%%%%%%
	%%%%%%%%%%%%%%%%%%%%%%%%%%%%%%%%%
	%%%%%%%%%%%%%%%%%%%%%%%%%%%%%%%%%
	
	\subsection{Configurations of circles} \label{subsec:config_circles}
	
	Let $\Circles = (\C_1, \C_2, \C_3)$ be an \emph{(ordered) configuration} of three distinct circles in $\cp$. The configuration $\Circles$ is \emph{non-degenerate} if every pair $\C_i$,$\C_j$ intersects in exactly two points $\{x_{ij},y_{ij}\}$, and the set of pairwise intersection points has at least four elements.  Henceforth, all configurations will be assumed to be non-degenerate. Also notice that by definition $\Circles$ is an ordered triple.
	
\begin{figure}
		\centering
		\includegraphics[width=\textwidth]{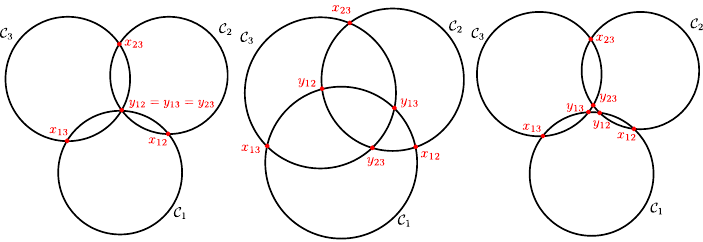}
		\caption{A Euclidean, spherical and hyperbolic configuration of circles, from left to right.}
		\label{fig:type_of_config}
	\end{figure}

	A configuration of circles is \emph{Euclidean} if the circles have a common intersection point. In this case  there are exactly four intersection points. If the configuration is not Euclidean, since every circle divides $\cp$ into two disjoint regions, then $\C_1$ separates $\{x_{23},y_{23}\}$ if and only if $\C_2$ separates $\{x_{13},y_{13}\}$ if and only if $\C_3$ separates $\{x_{12},y_{12}\}$. In that case, we say that the configuration $\Circles$ is \emph{spherical}. Otherwise, it is $\emph{hyperbolic}$ (cf. Figure~\ref{fig:type_of_config}).
	
	\begin{remark}\label{rem:number_of_configurations}
	    A configuration of circles induces a CW-structure on $\cp$, in which the $2$--cells are either bigons, triangles or quadrilaterals; in the spherical case the structure is simplicial and isomorphic to an octahedron.
	    Given two configurations of circles $\Circles^i = (\C_1^i, \C_2^i, \C_3^i)$ of the same kind (Euclidean, spherical or hyperbolic), there is always (at least) one CW--isomorphism of $\cp$ mapping $\C_k^1$ to $\C_k^2$. For spherical and hyperbolic configurations, it is enough to consider orientation preserving CW--isomorphisms. On the other hand, if $\Circles = (\C_1, \C_2, \C_3)$ is a Euclidean configuration of circles, 
	    there is no orientation preserving CW--isomorphism mapping $(\C_1, \C_2, \C_3)$ to $(\C_1, \C_3, \C_2)$: the obstruction being the cyclic order of the circles at the common intersection point.
	    
	\end{remark}
	
\begin{figure}
		\centering
		\includegraphics[width=\textwidth]{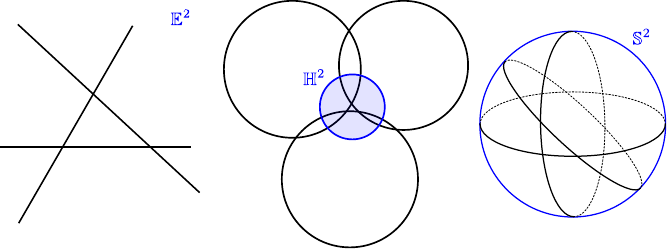}
		\caption{Euclidean, hyperbolic and spherical configurations are related to the corresponding geometries.}
		\label{fig:geom_of_config}
	\end{figure}

	The connection between a configuration of circles and the corresponding geometries is well known. We recall it in the next result (cf. Figure~\ref{fig:geom_of_config}).
	
	\begin{lemma}\label{lem:geometric_config_circles}
		Let $\Circles$ be a configuration of three circles. \begin{itemize}
			\item If $\Circles$ is Euclidean, let $y$ be the common intersection point. Then $\cp\setminus \{y\}$ admits a Euclidean metric for which the circles in $\Circles$ are geodesics.
			\item If $\Circles$ is spherical, then there is a M\"obius transformation $G \in \pslc$ such that $G \cdot \Circles$ are great circles for the underlying spherical metric.
			\item If $\Circles$ is hyperbolic, then there is a unique circle $\C_\HH$ orthogonal to every circle in $\Circles$. In particular, each connected component $D_\HH$ of $\cp \setminus \C_\HH$ admits a hyperbolic metric for which the intersections of $\Circles$ with $D_\HH$ are geodesics.
		\end{itemize} 
	\end{lemma}
	
	\begin{figure}
		\centering
		\includegraphics[width=4cm]{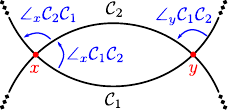}
		\caption{The anticlockwise angle between two circles at a point of intersection.}
		\label{fig:angles_config}
	\end{figure}
	
	Any two distinct circles $\C_1,\C_2$ in a configuration $\Circles$ intersect in two points. If $x$ is a point of intersection, then we can use the orientation of $\cp$ to determine the \emph{anticlockwise angle} $\cangle{x}{\C_1\C_2}$ from $\C_1$ to $\C_2$ at $x$ (cf. Figure~\ref{fig:angles_config}). We have that
	$$
	\cangle{x}{\C_2\C_1} = \pi - \cangle{x}{\C_1\C_2} = \cangle{y}{\C_1\C_2},
	$$
	where $y$ is the other point of intersection of $\C_1$ and $\C_2$. It is a simple exercise in complex projective geometry to show that a configuration of circles is uniquely determined (up to M\"obius transformations) by the ordered triple of angles at three points.
	
	\begin{lemma}\label{lem:points_angles_circles}
		For $i \in \{1,2\}$, let $\Circles^i = (\C_1^i, \C_2^i, \C_3^i)$ be a configuration of circles. For every pair of circles in $\Circles^i$ let $x_{jk}^i \in \C_j^i \cap \C_k^i$ be an intersection point such that
		$$
		\cangle{x_{12}^1}{\C_1^1\C_2^1} = \cangle{x_{12}^2}{\C_1^2\C_2^2}, \quad \cangle{x_{23}^1}{\C_2^1\C_3^1} = \cangle{x_{23}^2}{\C_2^2\C_3^2}, \quad \text{and} \quad
		\cangle{x_{13}^1}{\C_1^1\C_3^1} = \cangle{x_{13}^2}{\C_1^2\C_3^2}.
		$$
		Then there is a M\"obius transformation $M \in \pslc$ such that $M\cdot \Circles^1 = \Circles^2$ with $M \cdot x_{jk}^1 = x_{jk}^2$.
	\end{lemma}

	\subsection{Elliptic M\"obius transformations}\label{subsec:elliptic_mobius_transf}
	In this section we prove a correspondence between configurations of circles and certain triples of elliptic M\"obius transformations (see Corollary~\ref{cor:circles_reps}).
	
	As defined above, a non-trivial M\"obius transformation $G\in \pslc$ is said to be elliptic if $\tr(G)^2$ is real and $\tr(G)^2 < 4$. An elliptic transformation fixes exactly two points of $\cp$. Let $G\in \pslc$ be elliptic. The \emph{rotation angle} $\Rot(G,x) \in (0,2\pi)$ \emph{of $G$ at a fixed point $x$} is the angle of anticlockwise rotation of $G$ at $x$ (more precisely of $dG_x$ on $T_x\cp$). If $x,y$ are the fixed points of $G$, a M\"obius transformation mapping $x,y$ to $0,\infty$ conjugates $G$ to the element of $\pslc$
	\begin{equation} \label{eq:elliptic_normalisation}
	\begin{bmatrix}
	e^{i\frac{\Rot(G,x)}{2}}     & 0 \\
	0   & e^{-i\frac{\Rot(G,x)}{2}}
	\end{bmatrix}.
	\end{equation}
	
	The definition of rotation angle implies the following result.
	
	\begin{lemma}\label{lem:rots}
		Let $G\in \pslc$ be elliptic with fixed points $\{x,y\}$. Then
		$$
		\Rot(G,y)=2\pi-\Rot(G,x)=\Rot(G^{-1},x).
		$$
	\end{lemma}
	
	The \emph{rotation invariant} of an elliptic transformation $G$ is the unordered pair $\Rot(G):=\{\Rot(G,x),\Rot(G,y)\}$.
	
	\begin{lemma}\label{lem:rottrace}
		Let $G\in \pslc$ be elliptic, and let $\theta \in (0,2\pi)$. Then $\theta \in \Rot(G)$ if and only if $4\cos^2(\frac{\theta}{2})=\tr^2(G)$.
	\end{lemma}
	\proof
	Both the rotation angle and the trace operator are invariant under conjugation, thus we may assume that $G$ is normalised as in~\eqref{eq:elliptic_normalisation}. The equation $4\cos^2(\frac{\theta}{2})=\tr^2(G)$ has precisely two solutions in $(0,2\pi)$, of the form 
	$$\theta_1 = 2\arccos \left(\frac{|\tr(G)|}{2} \right) \qquad \text{and} \qquad \theta_2 = 2\pi - 2\arccos \left(\frac{|\tr(G)|}{2} \right),$$
	where we fix a determination of $\arccos$ in $[0,\pi]$. A direct computation shows that $\Rot(G) = \{\theta_1,\theta_2\}$ concluding the proof.
	\endproof
	
	Given the fixed points of $G$, the rotation invariant is enough to determine $G$ up to inversion, while the rotation angle is a complete invariant.
	
	\begin{lemma}\label{lem:rot_conj}
		Let $G, H \in \pslc$ be two elliptic transformations. Then 
		\begin{enumerate}
			\item \label{item:conj} $\Rot(G)=\Rot(H) \iff \tr^2(G)=\tr^2(H) \iff G,H$ are conjugate.
			\item If $G,H$ have the same fixed points $\{x,y\}$, then
			$$\Rot(G)=\Rot(H) \iff G=H^{\pm 1},$$
			and in particular
			$$\Rot(G,x)=\Rot(H,x) \iff G=H.$$
		\end{enumerate}
	\end{lemma}
	\proof $\ $
	\begin{enumerate}
		\item Two elliptics with the same rotation invariants must have the same trace squared by the previous Lemma~\ref{lem:rottrace}. But this is a complete invariant of conjugacy classes for semisimple elements of $\pslc$.
		\item Since $G$ and $H$ share the same fixed points, we can simultaneously normalise them as in~\eqref{eq:elliptic_normalisation}. Both statements follow from comparing the two normal forms.
	\end{enumerate}
	\endproof
	
	Next we analyse the connection between elliptic transformations, whose product is elliptic, and configurations of circles in $\cp$. First, we recall the following result from~\cite[lemma 3.4.1]{GKM00}. 
	
	\begin{lemma}\label{lem:from_elliptics_to_circles}
		Let $G,H \in \pslc$ be elliptic transformations with at most one common fixed point, and such that the product $GH$ is elliptic. Then the fixed points of $G$ and $H$ are contained in a unique circle $\C_{G,H}$.
	\end{lemma}
	
	We recall that given any two distinct circles $\C_1,\C_2$ intersecting at a point $x$, the (anticlockwise) angle from $\C_1$ to $\C_2$ at $x$ is denoted by $\cangle{x}{\C_1\C_2}$ (cf. \S\ref{subsec:config_circles}).
	
	\begin{lemma}\label{lem:reflection_circles}
		Let $\C_1,\C_2$ be distinct circles in $\cp$ meeting exactly at two points $x,y$. Let $J_i$ denote the reflection in $\C_i$. Then the product $G=J_2J_1$ is an elliptic transformation fixing $x,y$ with
		$$
		\Rot(G,x)=2\cangle{x}{\C_1\C_2}, \qquad \text{ and } \qquad \Rot(G,y)=2\cangle{y}{\C_1\C_2}.
		$$
	\end{lemma}
	\proof
	Since M\"obius transformations are conformal, we can normalise so that $x=0$ and $y=\infty$. Under the standard identification $\cp = \CC \cup \{\infty\}$, we can further normalise so that $\C_1=\RR \cup \{\infty\}$. Then $\C_2$ is a Euclidean line through $0$ and $\infty$. In this setting
	$$
	J_1(z)=\bar z, \qquad \text{and} \qquad J_2(z)=e^{i2(\cangle{x}{\C_1\C_2})}\bar z,
	$$
	and the statement follows from a direct computation.
	\endproof
	
	Henceforth we fix the following notation. Given $G,H$ distinct elliptic transformations whose product $GH$ is elliptic, we denote by $\{p_G,q_G\}$ (resp. $\{p_H,q_H\}$) the fixed points of $G$ (resp. $H$), by $\C_{G,H}$ the unique circle through $\{p_G,q_G,p_H,q_H\}$ (cf.~Lemma~\ref{lem:from_elliptics_to_circles}), and by $J_{G,H}$ the reflection about $\C_{G,H}$.
	
	\begin{lemma}\label{lem:circles_to_transformations}
		Let $(A,B,C)$ be an ordered triple of elliptic transformations with at most one common fixed point, and such that $ABC=1$. Then
		\begin{enumerate}
			\item \label{item:circles_to_transformations_1} $\C_{A,C}\cap \C_{A,B}=\{p_A,q_A\}$.
			\item \label{item:circles_to_transformations_2}
			$2\cangle{p_A}{\C_{A,B}\C_{A,C}} = \Rot(A,p_A) \ $ and $\ 2\cangle{q_A}{\C_{A,B}\C_{A,C}} = \Rot(A,q_A)$.
			\item \label{item:circles_to_transformations_3} $A = J_{A,C}J_{A,B}$.
		\end{enumerate}
	\end{lemma}
	\proof
	We begin by noticing that two of the three elliptic transformations share a common fixed point $p$ if and only if $p$ is fixed by all three of them. Hence there are either four or six distinct fixed points. Then the first statement \eqref{item:circles_to_transformations_1} follows from Lemma~\ref{lem:from_elliptics_to_circles}.
	
	Next, we recall that M\"obius transformations are conformal, thus without loss of generality we can simultaneously normalize $(A,B,C)$ so that $(p_A,q_A,p_B) = (0,\infty,1)$. It follows that $\C_{A,B} = \RR \cup \{\infty\}$. If we let $\theta := \frac{\Rot(A,0)}{2}$, then the three elliptic transformations take the following forms
	
	$$A=
	\begin{bmatrix}
	e^{i\theta}  & 0  \\ 0  & e^{-i\theta}
	\end{bmatrix}, \quad B=
	\begin{bmatrix}
	a  & b  \\ c  & 
	\overline{a}
	\end{bmatrix},\quad \text{ and }  \quad C^{-1} = AB =
	\begin{bmatrix}
	ae^{i\theta} & be^{i\theta} \\ ce^{-i\theta} & \overline{a} e^{-i\theta}
	\end{bmatrix},
	$$
	where $|\tr B| = 2|\Re(a)|<2$ (the relation between the diagonal elements of $B$ is implied by the fact that $C$ is elliptic). We remind the reader that we are always taking representatives in $\slc$ modulo a sign. Using that $\det(B) =1$ and that $B$ fixes $1$, it follows that $b,c$ are purely imaginary. In particular, there are choices of signs for which
	$$
	b=-i \Im(a)\pm \sqrt{\Re^2(a)-1}, \qquad \text{and} \qquad c=i \Im(a)\pm \sqrt{\Re^2(a)-1}.
	$$
	
	We claim that $C^{-1}$ has 
	fixed points of the form $te^{i\theta}$ for $t \in \RR \setminus \{0\}$. Since $C$ and $C^{-1}$ have the same fixed points, this will imply that $\cangle{0}{\C_{A,B}\C_{A,C}} = \theta = \frac{\Rot(A,0)}{2}$. To this end, we look for real solutions of the equation 
	$$
	te^{i\theta}=AB \cdot  te^{i\theta}=\frac{ate^{2i\theta}+be^{i\theta}}{ct+\overline{a}e^{-i\theta}} \iff ct^2-2i \Im(ae^{i\theta})t-b=0.
	$$
	Since $b,c$ are purely imaginary, this polynomial has real roots if and only if its discriminant 
	$-4\Im(ae^{i\theta})^2+4bc$
	is negative. But that follows from
	$$
	1=\det(AB)=\|a\|^2-bc=\|ae^{i\theta}\|^2-bc=\Re(ae^{i\alpha})^2+\Im(ae^{i\theta})^2-bc,
	$$
	and
	$$
	2>|\tr(AB)|=|2\Re(ae^{i\theta})|.
	$$
	This concludes the proof of the first part of \eqref{item:circles_to_transformations_2}, while the rest follows from the definition of the anticlockwise angle between two circles and Lemma~\ref{lem:rots}.
	
	For the last statement of the lemma, recall that $G := J_{A,C}J_{A,B}$ is an elliptic M\"obius transformation with fixed points $\{p_A,q_A\}$ (cf. Lemma~\ref{lem:reflection_circles}). Then $G$ has the same fixed points and rotation angles as $A$, thus $G=A$ by Lemma~\ref{lem:rot_conj}.
	\endproof
	
	Lemmas~\ref{lem:from_elliptics_to_circles},~\ref{lem:reflection_circles} and ~\ref{lem:circles_to_transformations} have the following straightforward consequence.
	\begin{corollary}\label{cor:circles_reps}
		There is a bijection
		$$
		\left\{
		\begin{array}{c}
		\text{Configurations of } \\
		\text{three circles.}
		\end{array}
		\right\}
		\longleftrightarrow
		\left\{
		\begin{array}{c}
		\text{Ordered triples of elliptic} \\
		\text{transformations with at most one} \\
		\text{common fixed point and product $1$.}
		\end{array}
		\right\},
		$$
		where $(\C_1,\C_2,\C_3) \mapsto (J_3J_1,J_1J_2,J_2J_3)$ and $(A,B,C) \mapsto (\C_{A,B},\C_{B,C},\C_{A,C})$.
	\end{corollary}

	%%%%%%%%%%%%%%%%%%%%%%%%%%%%%%%%%
	%%%%%%%%%%%%%%%%%%%%%%%%%%%%%%%%%
	%%%%%%%%%%%%%%%%%%%%%%%%%%%%%%%%%
	
	\subsection{Triangular immersions} \label{subsec:tri_immersions}
	In this section we define certain immersions of the standard $2$--simplex in $\cp$. Lemmas~\ref{lem:tiny_angles}, \ref{lem:triangular_reg_for_small_angles} and \ref{lem:additional_euclidean} prove the existence of immersions with certain requirements on the angles at the vertices. These are the ones we call atomic, and are listed in Table~\ref{tab:hyperbolic},~\ref{tab:spherical} and~\ref{tab:Euclidean}. Then we study some invariants of such immersions, and conclude in Corollary~\ref{cor:atomic_angles_relations} that they are essentially determined by the image of the vertices, up to a minor ambiguity.
	
	Let $\triangle := \{(x_1,x_2,x_3) \in \mathbb{R}^3_{\geq 0} \ | \ x_1 + x_2 + x_3 = 1  \}$ be the \emph{standard $2$--dimensional simplex}. Let $\{V_1,V_2,V_3\} \subset \triangle$ be its set of vertices so that 
	$V_1=(1,0,0)$, $V_2=(0,1,0)$, $V_3=(0,0,1)$, and let $e_{ij} \subset \triangle$ be the edge between $V_i$ and $V_j$. We endow $\triangle$ with the orientation induced form the ordering $(V_1,V_2,V_3)$ of its vertices.\par
	
	A \emph{triangular immersion} is an orientation preserving immersion $\tau : \triangle \rightarrow \cp$ such that each $\tau(e_{ij})$ is contained in a circle. In particular, we require $\tau$ to be locally injective everywhere except at the vertices. When every $\tau(e_{ij})$ is contained in a great circle, the triangle $\triangle$ inherits a spherical metric with  geodesic boundary and cone angles at the vertices. This is usually referred to as a \emph{spherical triangular membrane} in the literature~\cite{ER04,MP16}. Triangular immersions are relevant to this paper as they produce natural examples of $\cp$--structures on the thrice--punctured sphere (cf. \S\ref{subsec:tri_structures} for details.)
	
	Henceforth, we will often make the abuse of notation of referring by $\tau$ both the triangular immersion and its image in $\cp$, when it is not necessary to make a distinction. The image of the vertices (resp. edges) of $\triangle$ are the \emph{vertices} (resp. \emph{edges}) of $\tau$. Since edges of $\tau$ are arcs of circles, $\tau$ has well defined \emph{angles} at the vertices. When $\tau$ is not locally injective at a vertex, the angle is larger than $2\pi$, and $\tau$ should be thought as ``spreading over'' $\cp$. The orientation of $\triangle$ and the ordering of its vertices induce an orientation on $\tau$, and an ordering of its vertices and of its angles (which agree with the orientation and ordering induced by the orientation of $\cp$).\\
	
	Configurations of circles and triangular immersions are related to one another. If $\tau$ is a triangular immersion, each one of its edges extends to a unique circle giving a (possibly degenerate) configuration $\Circles_\tau$ of three circles. In this case we say that $\Circles_\tau$ \emph{supports} $\tau$. When $\Circles_\tau$ is non-degenerate, we say that $\tau$ is \emph{non-degenerate}. When the interior of the image of $\tau$ is disjoint from $\Circles_\tau$, we say that $\tau$ is \emph{enclosed in} $\Circles_\tau$. These are exactly those triangular immersions whose (interior of the) images are the connected components of $\cp \setminus \Circles_\tau$. Necessary and sufficient conditions on the angles of $\tau$ for it to be enclosed in $\Circles_\tau$ are well known, but we provide a short proof as we could not find a direct reference.
	
	\begin{lemma}\label{lem:tiny_angles}
		Let $(a,b,c)$ be an ordered triple of angles in $(0,\pi)^3$.
		\begin{enumerate}
			\item (Euclidean Triangles.) There is a  Euclidean configuration of circles $\Circles$ and a triangular immersion $\tau$ enclosed in $\Circles$ with angles $(a,b,c)$ if and only if one of the following conditions is satisfied:
			\begin{equation}\label{eq:euclidean_tri}
			a+ b+  c = \pi, \quad  -a+ b+  c = \pi, \quad  a- b+  c = \pi, \quad  a+ b-  c = \pi.
			\end{equation}
			
			\item (Hyperbolic Triangles.) There is a  hyperbolic configuration of circles $\Circles$ and a triangular immersion $\tau$ enclosed in $\Circles$ with angles $(a,b,c)$ if and only if:
			\begin{equation} \label{eq:hyperbolic_tri}
			a+b+c < \pi.
			\end{equation}
			
			\item (Spherical Triangles.) There is a spherical configuration of circles $\Circles$ and a triangular immersion $\tri$ enclosed in $\Circles$ with angles $(a,b,c)$ if and only if $(a,b,c)$ satisfies:
			\begin{equation}\label{eq:spherical_tri}
			a+b+c > \pi, \qquad \text{and} \qquad 
			\begin{cases}
			a + \pi > b + c,\\
			b + \pi > c + a,\\
			c + \pi > a + b.
			\end{cases}
			\end{equation}
		\end{enumerate}
	\end{lemma}
	
	\proof $\ $
	
	\begin{figure}[t]
		\centering
		\includegraphics[width=0.9\textwidth]{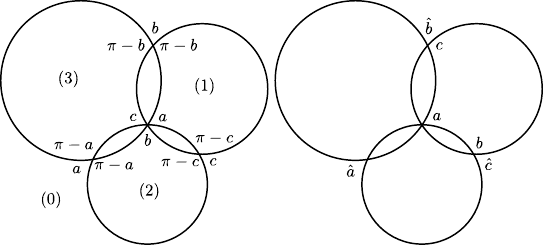}
		\caption{Two Euclidean configurations. Both support an enclosed triangular immersion with angles $(a,b,c) \in (0,\pi)^3$, such that either $a+b+c = \pi$ (left) or $-a+b+c = \pi$ (right).}
		\label{fig:tiny_angles_e}
	\end{figure}
	
	\begin{figure}[t]
		\centering
		\includegraphics[width=0.40\textwidth]{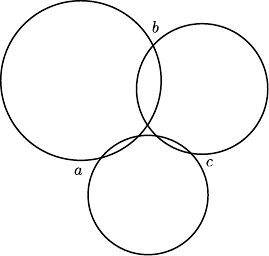}
		%\hspace{0.1\textwidth}
	    \includegraphics[width=0.40\textwidth]{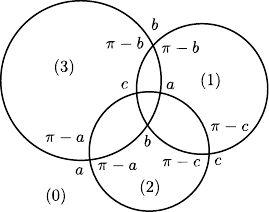}
		\caption{A hyperbolic configuration and a spherical configuration. They both support an enclosed triangular immersion with angles $(a,b,c) \in (0,\pi)^3$, such that either $a+b+c = \pi$ (left) or \eqref{eq:spherical_tri} is satisfied (right).}
		\label{fig:tiny_angles_hs}
	\end{figure}
	
	\begin{enumerate} 
		\item Let $\tau$ be a triangular immersion enclosed in a Euclidean configuration of circles $\Circles_\tau$. Then there is a common intersection point $y$, and $\cp \setminus \{y\}$ admits a Euclidean metric for which the circles in $\Circles_\tau$ are geodesics (cf. Lemma~\ref{lem:geometric_config_circles}). In this setting, it is easy to check that each one of the four triangular immersions that are enclosed in $\Circles_\tau$ have angles
		$$
		(0) \ (a,b,c), \quad 
		(1) \ (a, \pi-c, \pi-b), \quad (2) \ (\pi-c,b, \pi-a), \quad \text{and} \quad (3) \ (\pi-b, \pi-a, c),
		$$
		each one satisfying exactly one of the equalities in~\eqref{eq:euclidean_tri} (cf. Figure~\ref{fig:tiny_angles_e}).
		
		The converse implication is well known for $a+b+c=\pi$. If $-a+b+c=\pi$, we consider the angles $\hat{a} = a$, $\hat{b} = \pi-c$ and $\hat{c} = \pi - b$. Clearly $(\hat{a},\hat{b},\hat{c}) \in (0,\pi)^3$ and $\hat{a}+\hat{b}+\hat{c}=\pi$, therefore there is a Euclidean triangle with angles $(\hat{a},\hat{b},\hat{c})$ supported by some configuration of circles. One of the other enclosed triangular immersions has angles $(a,b,c)$ (cf. Figure~\ref{fig:tiny_angles_e}). The same strategy applies to the other cases.
		
		\item Let $\tau$ be a triangular immersion enclosed in a hyperbolic configuration of circles $\Circles_\tau$. Let $\C_\HH$ be the circle that is orthogonal to the family $\Circles_\tau$ (cf. Lemma~\ref{lem:geometric_config_circles}). In this case there are precisely two triangular immersions that are enclosed in $\Circles_\tau$, and they are both disjoint from $\C_\HH$. It follows that $\tau$ is a hyperbolic triangle in one of the two connected components of $\cp \setminus \C_\HH$, thus the inequality~\eqref{eq:hyperbolic_tri} is a consequence of the formula for hyperbolic area of triangles (cf. Figure~\ref{fig:tiny_angles_hs}). The converse implication is~\cite[Theorem 3.5.9]{RA06}.
		
		\item Finally, let $\tau$ be a triangular immersion enclosed in a spherical configuration of circles $\Circles_\tau$. By Lemma~\ref{lem:geometric_config_circles}, we can realize this configuration of circles by great circles. So every triangular region $\tau$ enclosed in $\Circles_\tau$ is a geodesic triangle for the standard spherical metric. By the area formula for spherical triangles, we have that 
		$$
		a+b+c=\pi+\text{Area}(\tau)>\pi.
		$$
		The other inequalities~\eqref{eq:spherical_tri} are obtained by applying Gauss-Bonnet to the enclosed triangular regions adjacent to $\tau$ (cf. Figure~\ref{fig:tiny_angles_hs}), whose angles are
		$$
		(1) \ (a, \pi-c, \pi-b), \quad (2) \ (\pi-c,b, \pi-a), \quad \text{and} \quad (3) \ (\pi-b, \pi-a, c).
		$$
		
		The converse implication is a simple adaptation of~\cite[Theorem 3.5.9]{RA06} using the law of cosines in spherical geometry (cf.~\cite[Exercise 2.5.8]{RA06})
	\end{enumerate}
	\endproof
	
	\begin{remark}\label{rem:uniqueness_tri_imm}
	    For convenience, Lemma~\ref{lem:tiny_angles} is stated just in terms of the existence of a triangular immersion $\tau$. Despite we will not need it, we remark that it is a simple consequence of Lemma~\ref{lem:points_angles_circles} that $\tau$ is also unique up to M\"obius transformations. The same is true for the following results.
	\end{remark}
	
\begin{figure}[t]
		\centering
		\includegraphics[width=\textwidth]{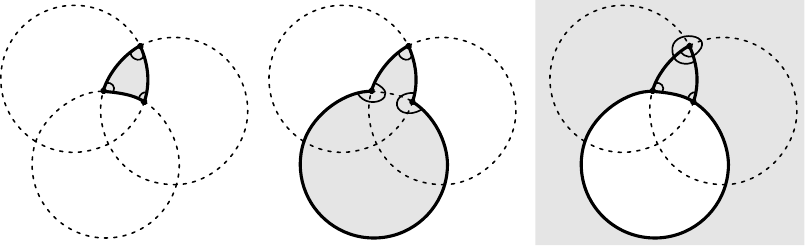}
		\caption{A triangular immersion can be manipulated to new triangular immersions by adding an entire disk, or by taking a full turn around a vertex.}
		\label{fig:tri_manipulations}
	\end{figure}

	Given an enclosed triangular immersion $\tau$, there are two simple operations that one can perform to construct new triangular immersions supported by the same configuration of circles. The first one consists in extending $\tau$ by a full disk, by ``pushing'' an edge of $\tau$ to its complement in its supporting circle (Figure~\ref{fig:tri_manipulations}).
	This operation increases the two angles adjacent to the pushed edge by $\pi$. The second manipulation involves making a full turn around a vertex, by extending the opposite edge to cover its entire supporting circle (Figure~\ref{fig:tri_manipulations}). This operation increases the angle at the highlighted vertex by $2\pi$. It will be remarked later on how these operations are related to grafting the associated triangular structure (cf. Example~\ref{ex:edge_graft}).
	
	On the other hand, there are triangular immersions that do not arise from these operations, whose existence we prove now.

	\begin{lemma}\label{lem:triangular_reg_for_small_angles}
		Let $(a,b,c)$ be an ordered triple of angles such that
		$$
		a \in (0,\pi) \cup (\pi,2\pi), \qquad \text{ and } \qquad b,c \in (0,\pi).
		$$
		Then there is a configuration of circles $\Circles$ and a triangular immersion $\tau$ supported by $\Circles$ with angles $(a,b,c)$.
	\end{lemma}
	
\begin{figure}[t]
		\centering
		\includegraphics[width=0.40\textwidth]{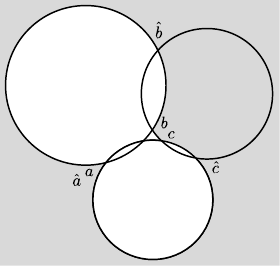}
		%\hspace{0.1\textwidth}
	    \includegraphics[width=0.40\textwidth]{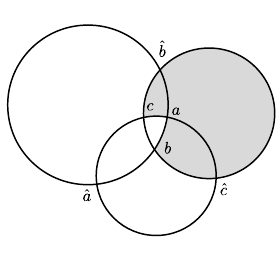}
		\caption{On the left, a triangular immersion on a hyperbolic configuration with angles $(a,b,c) \in (0,\pi)^3$ such that $a+b+c > \pi$ but $a + \pi < b + c$. On the right, a triangular immersion supported by a spherical configuration.}
		\label{fig:small_angles_hs}
	\end{figure}

	\proof
	First suppose $a \in (0,\pi)$. Those cases where $(a,b,c)$ satisfies one of the conditions \eqref{eq:euclidean_tri}, \eqref{eq:hyperbolic_tri} or \eqref{eq:spherical_tri} from Lemma~\ref{lem:tiny_angles} are covered by that lemma. Hence suppose $a+b+c > \pi$, but at least one of the other inequalities in~\eqref{eq:spherical_tri} is not satisfied. Up to permuting $a,b,c$ we may assume that $a + \pi < b + c$. Let
	$$
	\hat{a} = a, \quad \text{and} \quad \hat{b} = \pi - b, \quad \text{and} \quad \hat{c} = \pi -c.
	$$
	Then $\hat{a} + \hat{b} + \hat{c} = a + (\pi - b) + (\pi - c) < \pi$ by assumption, therefore by Lemma~\ref{lem:tiny_angles} there is a hyperbolic configuration of circles $\Circles$ and a triangular immersion $\hat \tau$ enclosed in $\Circles$ with angles $(\hat{a},\hat{b},\hat{c})$. Figure~\ref{fig:small_angles_hs} (on the left) shows that the same configuration of circles supports a triangular immersion with angles $(a,b,c)$.\par
	
\begin{figure}[t]
		\centering
		\includegraphics[width=0.9\textwidth]{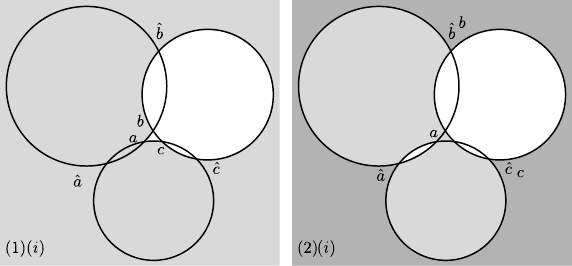}
		\linebreak
		\includegraphics[width=0.9\textwidth]{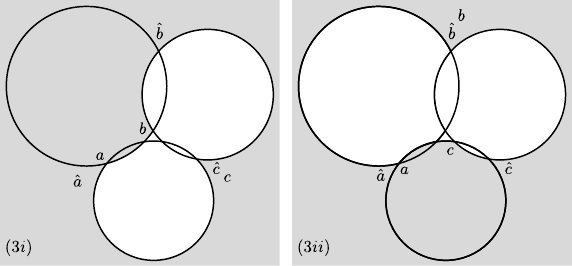}
		\caption{Different triangular immersions with angles $(a,b,c)$, supported by hyperbolic configurations. We remark that the one depicted in (2) covers the darker triangle twice. }
		\label{fig:small_angles_h2}
	\end{figure}

\begin{figure}[t]
		\centering
		\includegraphics[width=0.9\textwidth]{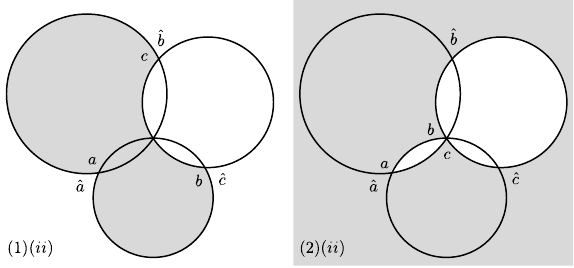}
		\linebreak
		\includegraphics[width=0.9\textwidth]{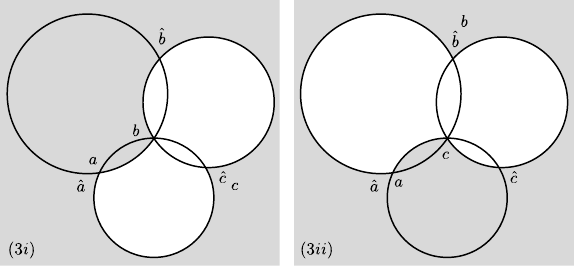}
		\caption{Different triangular immersions with angles $(a,b,c)$, supported by Euclidean configurations.}
		\label{fig:small_angles_e2}
	\end{figure}

	Now suppose $a \in (\pi,2\pi)$. Consider the following relations:
	$$
	(1)
	\begin{array}{cc}
	(i) \quad a + b + c > 3\pi,\\
	(ii) \quad a + b + c = 3\pi,\\
	\end{array}, \qquad (2)
	\begin{array}{cc}
	(i) \quad a - b - c > \pi,\\
	(ii) \quad a - b - c = \pi,\\
	\end{array} \qquad
	\begin{array}{cc}
	(3i) \quad a - b + c \leq \pi,\\
	(3ii) \quad a + b - c \leq \pi.
	\end{array}
	$$
	We observe that these three groups of inequalities are mutually exclusive, as any two of them imply the following contradictions:
	$$
	\begin{array}{c}
	(1) + (2) \Longrightarrow a \geq 2\pi, \\
	(3i) + (3ii) \Longrightarrow a \leq \pi,
	\end{array}
	 \quad
	\begin{array}{c}
	(1) + (3i) \Longrightarrow b \geq \pi,\\
	(1) + (3ii) \Longrightarrow c \geq \pi,
	\end{array}
	 \quad
	\begin{array}{c}
	(2) + (3i) \Longrightarrow c\leq 0,\\
	(2) + (3ii) \Longrightarrow b\leq 0.
	\end{array}
	$$
	If one of those inequalities is satisfied, we define
	$$
	\begin{cases}
	\hat{a} = 2\pi - a, \quad \hat{b} = \pi - b, \quad \text{and} \quad \hat{c} = \pi - c, & \text{ if } (1)(i) \text{ is satisfied,}\\
	\hat{a} = 2\pi - a, \quad \hat{b} = \pi - c, \quad \text{and} \quad \hat{c} = \pi - b, & \text{ if } (1)(ii) \text{ is satisfied,}\\
	\hat{a} = 2\pi - a, \quad \hat{b} = b, \quad \text{and} \quad \hat{c} = c, & \text{ if } (2)(i) \text{ is satisfied,}\\
	\hat{a} = 2\pi - a, \quad \hat{b} = c, \quad \text{and} \quad \hat{c} = b, & \text{ if } (2)(ii) \text{ is satisfied,}\\
	\hat{a} = a - \pi, \quad \hat{b} = \pi - b, \quad \text{and} \quad \hat{c} = c, & \text{ if } (3i) \text{ is satisfied,}\\
	\hat{a} = a - \pi, \quad \hat{b} = b, \quad \text{and} \quad \hat{c} = \pi - c, & \text{ if } (3ii) \text{ is satisfied.}
	\end{cases}
	$$
	In each case, the assumption implies that $\hat{a} + \hat{b} + \hat{c} \leq \pi$, therefore Lemma~\ref{lem:tiny_angles} applies to give a Euclidean or hyperbolic configuration of circles $\Circles$ and a triangular immersion $\hat \tau$ enclosed in $\Circles$ with angles $(\hat{a},\hat{b},\hat{c})$. Figures~\ref{fig:small_angles_h2} and \ref{fig:small_angles_e2} show that the same configuration of circles supports a triangular immersion with angles $(a,b,c)$.\par
	
	Finally, let $(\neg 1),(\neg 2),(\neg 3i),(\neg 3ii)$ be the opposite of the inequalities $(1),(2),(3i),(3ii)$, and suppose $(a,b,c)$ satisfies all of $(\neg 1),(\neg 2),(\neg 3i),(\neg 3ii)$. We define
	$$
	\hat{a} = 2\pi - a, \quad \text{and} \quad  \hat{b} = \pi - b, \quad \text{and} \quad \hat{c} = \pi-c.
	$$
	This time $\hat{a} + \hat{b} + \hat{c} = 4\pi -a - b - c > \pi$ because of $(\neg 1)$. Moreover,
	\begin{alignat*}{2}
		\hat{a} + \pi = 3\pi - a &> 2\pi - b - c = \hat{b} + \hat{c}, \qquad &&\text{ by } (\neg 2) ,\\
		\hat{b} + \pi = 2\pi - b &> 3\pi - a - c = \hat{a} + \hat{c}, \qquad &&\text{ by } (\neg 3i),\\
		\hat{c} + \pi = 2\pi - c &> 3\pi - a - b = \hat{a} + \hat{b}, \qquad &&\text{ by } (\neg 3ii).
	\end{alignat*}
	By Lemma~\ref{lem:tiny_angles}, there is a spherical configuration of circles $\Circles$ and a triangular immersion $\hat \tau$ enclosed in $\Circles$ with angles $(\hat{a},\hat{b},\hat{c})$. See Figure~\ref{fig:small_angles_hs} for a triangular immersion with angles $(a,b,c)$ supported by the same configuration $\Circles$.\par
	\endproof
	
	Due to the degenerate nature of Euclidean configurations, there is one additional case that needs to be considered,  which we address in the next lemma.
	
	\begin{lemma}\label{lem:additional_euclidean}
		Let $(a,b,c)$ be an ordered triple of angles such that
		$$
		a \in (2\pi,3\pi), \quad b,c \in (0,\pi) \qquad \text{ and } \qquad a-b-c = \pi.
		$$
		Then there is a configuration of circles $\Circles$ and a triangular immersion $\tau$ supported by $\Circles$ with angles $(a,b,c)$.
	\end{lemma}
	
\begin{figure}[t]
    		\centering
    		\includegraphics[width=0.45\textwidth]{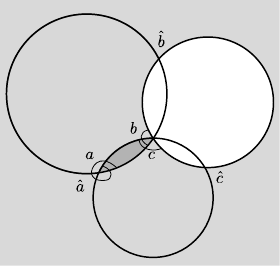}
    		\caption{The additional triangular immersion mentioned in Lemma~\ref{lem:additional_euclidean}. Notice that $a > 2\pi$, hence the darker bigon is covered twice.}
    		\label{fig:additional_euclidean}
	    \end{figure}

	\proof
	    Let
    	$$
    	\hat{a} = a-2\pi, \quad \text{and} \quad \hat{b} = \pi-b, \quad \text{and} \quad \hat{c} = \pi-c.
    	$$
    	Then $\hat{a},\hat{b},\hat{c} \in (0,\pi)$ and $\hat{a} + \hat{b} + \hat{c} = a-2\pi + \pi - b + \pi - c = \pi$, therefore by Lemma~\ref{lem:tiny_angles} there is a Euclidean configuration of circles $\Circles$ and a triangular immersion $\hat \tau$ enclosed in $\Circles$ with angles $(\hat{a},\hat{b},\hat{c})$. Figure~\ref{fig:additional_euclidean} shows that the same configuration of circles supports a triangular immersion with angles $(a,b,c)$.
    	
	\endproof
	
	The triangular immersions constructed in the proofs of Lemmas~\ref{lem:tiny_angles},~\ref{lem:triangular_reg_for_small_angles} and~\ref{lem:additional_euclidean}, which are depicted in Figures~\ref{fig:small_angles_hs},~\ref{fig:small_angles_h2},~\ref{fig:small_angles_e2}, and~\ref{fig:additional_euclidean}, are the starting point to construct all complex projective structures of interest in this paper. 	For this reason, we will refer to them as the \emph{atomic} triangular immersions. They are \emph{Euclidean/hyperbolic/spherical} depending on the type of the underlying configuration of circles.
	In Lemmas~\ref{lem:triangular_reg_for_small_angles} and~\ref{lem:additional_euclidean} exactly one angle is allowed to be larger than $\pi$, and we have assumed that to be the first one for simplicity. This normalization is inessential, and the same statements and proofs hold if one chooses a different angle to be the large one. This should be regarded as a change of marking (i.e. a permutation of the vertices of the simplex on which the triangular immersions are defined), and we call \textit{atomic triangular immersion} any triangular immersion obtained in this way. Theorem \hyperref[thm:mainthm]{B} and Corollary \ref{cor:atomic=ungrafted} will show that, in a precise sense, this is indeed the minimal collection of triangular immersions to be considered.
	
	We remark that the proofs of these lemmas are explicit, and construct a concrete collection of triangular immersions. Notice that for every triple of real numbers $(a,b,c)$, two of which are in $(0,\pi)$ and one is $(0,\pi) \cup (\pi,2\pi) \cup (2\pi,3\pi)$, there is a unique atomic triangular immersion with those angles. This allows us to organize the atomic triangular immersions in Tables~\ref{tab:hyperbolic},~\ref{tab:spherical} and ~\ref{tab:Euclidean}. We now define the other features listed in those tables. 
	
	Let $\tau:\triangle\to \cp$ be an atomic triangular immersion, and let $\Circles_\tau$ be the configuration of circles that supports it. The configuration $\Circles_\tau$ is either of spherical, Euclidean or hyperbolic \textit{type}. The \textit{target angles} of $\tau$ are the numbers $(\hat a,\hat b,\hat c)$ defined as follows.
	    \begin{itemize}
	        \item If $(a,b,c)$ satisfies the hypothesis of Lemma~\ref{lem:tiny_angles} then $(\hat a,\hat b,\hat c)=(a,b,c)$.
	        \item If $(a,b,c)$ does not satisfy the hypothesis of Lemma~\ref{lem:tiny_angles} then $(\hat a,\hat b,\hat c)$ is defined as in the proofs of Lemmas~\ref{lem:triangular_reg_for_small_angles} and~\ref{lem:additional_euclidean}, depending on what conditions are satisfied, and up to permuting the angles as appropriate. 
	   \end{itemize}
	The target angles of $\tau$ satisfy the hypothesis of Lemma~\ref{lem:tiny_angles}. Therefore there is a triangular immersion $\hat \tau$ with angles $(\hat a,\hat b,\hat c)$, which we call the \textit{target triangular immersion}. If $\Circles_{\tau}:=(\C_{12},\C_{23},\C_{13})$, it follows from the construction that $\hat \tau$ is  supported either by $\Circles_\tau$ or by $\Circles_\tau^* := (\C_{12},\C_{13},\C_{23})$, but the latter only happens in the Euclidean cases of Figure~\ref{fig:tiny_angles_e} (right) and Figures~\ref{fig:small_angles_e2} (1) and~(2). In addition, $\hat \tau$ is always enclosed (while $\tau$ may not be). 
	All the above pictures representing the atomic triangular immersions have been normalized so that $\hat \tau(\triangle)$ contains the point at infinity in its interior.
	    
	For pairwise distinct $i,j,k \in \{1,2,3\}$, consider the circle $\C_{ij}\in \Circles_\tau$ supporting $\hat \tau (e_{ij})$; the intersection $\C_{ij}\cap \C_{jk}$ consists of two points: one is $\hat \tau (V_j)$, and we define $\hat \tau (V_j)'$ to be the other one. The collection $\{\hat \tau (V_j),\hat \tau (V_j)' \ | \ j=1,2,3 \}$ accounts for all the points of intersection of the circles in $\Circles_\tau$, which are the possible vertices for $\tau$. Note that by construction we always have $ \{\tau(V_1),\hat\tau(V_1) \} \subseteq \C_{12}\cap \C_{13}$.
	We say a vertex $\tau(V_j)$ of $\tau$ is \textit{positive} if there exists $k$ such that $\tau(V_j)=\hat \tau(V_k)$, i.e. if it coincides with a vertex of $\hat \tau$, and we say it is \textit{negative} otherwise. 
	This defines a triple of \textit{signs} $(s_1(\tau),s_2(\tau),s_3(\tau))\in \{ \pm\}^3$ associated to $\tau$.
	In the Euclidean case, we additionally decorate this triple: we define it to be $(s_1(\tau),s_2(\tau),s_3(\tau))$ when $\hat \tau$ is supported by $\Circles_{\tau}$, and to be $(s_1(\tau),s_2(\tau),s_3(\tau))^*$ when $\hat \tau$ is supported by $\Circles_{\tau}^*$.

	\begin{remark}
	The Euclidean case (see Table~\ref{tab:Euclidean}) displays all possible cases for the triple of signs, including the extra $*$ decoration, with the only exception of the cases in which all vertices are negative. This cannot happen as it would mean that $\tau$ maps all vertices to the common intersection point of the configuration of circles, but this never happens for an atomic triangular immersion.
	The extra $*$ decoration is not needed for the hyperbolic and spherical cases as they are less degenerate than the Euclidean ones, in the sense that circles in $\Circles$ have six distinct intersection points, which allows for more flexibility in the definition of the atomic immersions. See Tables \ref{tab:hyperbolic} and \ref{tab:spherical}.
	In the hyperbolic case we find all possible cases for the signs. In the spherical case we only see the triples $(\pm,\pm,\pm)$. This is because a spherical configuration of circles has only triangular complementary regions (while the complement of a hyperbolic configuration has different shapes, with only two triangles). As a result it is much easier for a spherical atomic triangular immersion to be enclosed, and equal to its own target triangular immersion.
	\end{remark}

\begin{lemma}\label{lem:target_find}
	 Let $\tau$ be an atomic triangular immersions supported by a configuration of circles $\Circles$. 
	 Then $\hat \tau$ is uniquely determined by $\Circles$ and the vertices of $\tau$.
\end{lemma}
\proof
Let $\Circles=(\C_{12},\C_{23},\C_{13})$ and recall that  we have $\tau(e_{jk})\subseteq \C_{jk}$ for $j,k=1,2,3$, by definition of what it means for a triangular immersion to be supported by a configuration of circles. Moreover by construction $ \{\tau(V_1),\hat\tau(V_1) \} \subseteq \C_{12}\cap \C_{13}$.
	 
	 Suppose that $\Circles$ is Euclidean. Then $\hat \tau$ is the unique enclosed triangular immersion mapping to the Euclidean triangle such that $\hat\tau(V_1)=(\C_{12}\cap \C_{13}) \setminus \{\infty\}$. 
	 
	 Next, if $\Circles$ is hyperbolic, then let $\C=\C_\HH$ be the dual circle from Lemma~\ref{lem:geometric_config_circles}. If $\Circles$ is spherical, then let $\C$ be a circle which separates the vertices of $\tau$ from the other intersection points of circles in $\Circles$. In either case $\hat \tau$ is the unique enclosed triangular immersion which has image disjoint from $\C$, is supported by $\Circles$, and such that $\{\tau(V_1),\hat\tau(V_1) \} \subseteq \C_{12}\cap \C_{13}$. We additionally remark that $\hat \tau$ is always on the left of $\C$ with respect to the orientation induced by $\Circles$.
\endproof

	 \begin{corollary}\label{cor:atomic_angles_relations}
	    Let $\Circles$ be a configuration of circles.
	    Let $\tau_1,\tau_2$ be two atomic triangular immersions supported by $\Circles$, such that $\tau_1(V_j) = \tau_2(V_j)$ for all $j \in \{1,2,3\}$. Then $\hat \tau_1 =\hat \tau_2$.
	    Moreover, if $(a_i,b_i,c_i)$ are the angles of $\tau_i$, then exactly one of the following happens:
	 \begin{enumerate}
	     \item \label{item:same_angles_0} $(a_1,b_1,c_1)=(a_2,b_2,c_2)$ and $\tau_1 = \tau_2$;
	     \item $(a_1-a_2,b_1-b_2,c_1-c_2)=(\pi,-\pi,0)$ up to permutation.
	 \end{enumerate}
	 \end{corollary}
	 \proof
	    The first assertion follows directly from Lemma~\ref{lem:target_find}. As a direct consequence, $\tau_1,\tau_2$ have the same target angles and the same triple of signs. A direct inspection of Tables \ref{tab:hyperbolic}, \ref{tab:spherical}, \ref{tab:Euclidean} proves the desired relations between the angles, just by imposing equalities of the respective target angles. In particular, recall that atomic triangular immersions are uniquely determined by their angles, hence $(a_1,b_1,c_1)=(a_2,b_2,c_2)$ implies $\tau_1 = \tau_2$.
	 \endproof
	 
\begin{figure}[t]
    	\centering
    	\includegraphics[width=0.9\textwidth]{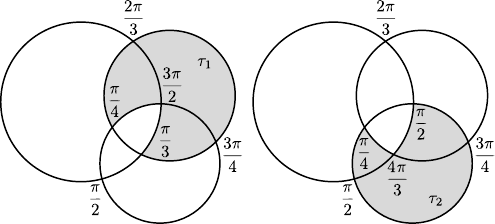}
    	\caption{Two atomic triangular immersions supported by the same (spherical) configuration of circles, and with the same signs $(-,-,-)$.}
    	\label{fig:mickey_mouse}
	\end{figure}

	 \begin{example}
	     Let $\tau_1,\tau_2$ be two atomic triangular immersions with angles
	     $$
	     (a_1,b_1,c_1) = \left( \frac{3\pi}{2},\frac{\pi}{3},\frac{\pi}{4} \right) \quad \text{and} \quad (a_2,b_2,c_2) = \left( \frac{\pi}{2},\frac{4\pi}{3},\frac{\pi}{4} \right).
	     $$
	     These immersions correspond to the second and third row of Table~\ref{tab:spherical}, respectively. They are supported by the same spherical configuration of circles $\Circles$, with target angles $(\hat a,\hat b,\hat c) = \left( \frac{\pi}{2},\frac{2\pi}{3},\frac{3\pi}{4} \right)$, and share the same signs $(-,-,-)$. In particular, $\hat \tau_1 =\hat \tau_2$. Furthermore, $\tau_1$ can be transformed into $\tau_2$ by first adding a disk and then removing another disk (cf. Figure~\ref{fig:mickey_mouse}).
	 \end{example}
	 
	 \begin{remark}\label{rem:rigid_signs}
	 Some of the sign invariants in each of the Tables \ref{tab:hyperbolic}, \ref{tab:spherical}, \ref{tab:Euclidean} occur exactly once. 
	 If two triangular immersions have same such signs, then they are equal by Corollary~\ref{cor:atomic_angles_relations} case \eqref{item:same_angles_0}. This applies for instance to atomic triangular immersions arising from Lemma~\ref{lem:additional_euclidean}, depicted in Figure~\ref{fig:additional_euclidean}.
	 \end{remark}
    
	%%%%%%%%%%%%%%%%%%%%%%%%%%%%%%%%%
	%%%%%%%%%%%%%%%%%%%%%%%%%%%%%%%%%
	%%%%%%%%%%%%%%%%%%%%%%%%%%%%%%%%%
	
	\section{Tame and relatively elliptic \texorpdfstring{$\cp$}--structures}\label{sec:tame_and_rel_elliptic_structures}
	
	In this chapter we define the geometric structures of interest in this paper, and study the geometry they induce on the universal cover. The reader can find the proofs of Theorems~\hyperref[thm:embedding]{C} and \hyperref[thm:maxneighborhood]{D} in \S\ref{subsec:mobius_completion} and \S\ref{subsec:canonic_max_nbhd} respectively.

	Let $\Sigmac$ be a closed oriented surface and let $\{x_1,\dots,x_n\} \subset \Sigmac$ be $n$ distinct points such that the \emph{punctured surface} $\Sigma := \Sigmac \setminus \{x_1,\dots,x_n\}$ has negative Euler characteristic. If $g$ is the genus of $\Sigmac$, this is equivalent to $2g + n > 2$, and it implies that $\Sigma$ admits a complete hyperbolic metric of finite area. The points $\{x_1,\dots,x_n\}$ are the \emph{punctures} of $\Sigma$.
	
	A \emph{complex projective structure} ($\cp$--structure in short) on $\Sigma$ is a maximal atlas of charts into $\cp$ with transition maps in $\pslc$ (see \cite{G67,DU}). A $\cp$--structure can be described by a \emph{developing pair} $(\dev,\hol)$  consisting of  a \emph{developing map} and a \emph{holonomy representation}
	$$
	\dev : \widetilde{\Sigma} \rightarrow \cp, \qquad \hol : \pi_1(\Sigma) \rightarrow \pslc,
	$$
	satisfying the equivariance condition
	$$
	\dev( \gamma \cdot x) = \hol(x) \cdot \dev(x), \qquad \text{for all } \ x \in \widetilde{\Sigma}, \quad \gamma \in \pi_1(\Sigma).
	$$
	
	There is a natural equivalence relation on the set of complex projective structures on a surfaces for which two pairs $(\dev,\hol)$ and $(\dev',\hol')$ are equivalent if there is $A\in \pslc$ so that   $\dev'=A\circ \dev$ and $\rho'=A\rho A^{-1}$ (up to isotopy  of $\Sigma$). The \emph{deformation space of marked $\cp$--structures on $\Sigma$} is the space of equivalence classes of complex projective structures and it is denoted by $\Moduli(\Sigma)$. 
	We denote by $\Char(\Sigma)$ the \textit{space of conjugacy classes of representations} of $\pi_1(\Sigma)$ into $\pslc$. We prefer not to use the GIT quotient because some of the representations of interest in this paper are reducible.
	The \emph{holonomy map} is the forgetful map
	$$
	\Hol : \Moduli(\Sigma) \rightarrow \Char(\Sigma), \qquad [(\dev,\hol)] \mapsto [\hol].
	$$

	Every $\cp$--structure has a natural underlying complex structure (or equivalently a conformal structure). We define  $\Modulic(\Sigma) \subset \Moduli(\Sigma)$ to be the subset of $\cp$--structures on $\Sigma$ whose underlying conformal structure around every puncture is the complex punctured disk $\pDisk := \{ z \in \CC \ | \ 0 < |z| < 1 \}$. 
	
	The space of interest in this paper is the subspace $\Modulit(\Sigma)$ of $\Modulic(\Sigma)$ of those structures whose developing map is \emph{tame} and whose holonomy representation is \emph{relatively elliptic}. We will define these terms in \S\ref{subsec:tame_rel_ell}.

	%%%%%%%%%%%%%%%%%%%%%%%%%%%%%%%%%
	%%%%%%%%%%%%%%%%%%%%%%%%%%%%%%%%%
	%%%%%%%%%%%%%%%%%%%%%%%%%%%%%%%%%

	\subsection{Ends, framing, and grafting}\label{subsec:tame_rel_ell}
	Let $\Sigmaw$ be the topological universal cover of $\Sigma$, and choose an identification $\Sigmaw \cong \HH^2$ coming from a uniformization of $\Sigma$ as a complete hyperbolic surface of finite area. An \emph{end} $\End$ of $\Sigmaw$ is defined to be the fixed point of a parabolic deck transformation in the boundary of $\HH^2$ in the closed disk model. For every puncture $x$ of $\Sigma$, we denote by $\Ends_x(\Sigmaw)$ the \emph{set of ends covering} $x$ (see Remark~\ref{rem:end_cover_puncture} for more details), and by
	$$
	\Ends(\Sigmaw) :=\cup_x \Ends_x(\Sigmaw),
	$$
	the \emph{set of all ends}. The \emph{end-extension} of $\Sigmaw$ is the topological space $\Sigmawe=\Sigmaw\cup \Ends(\Sigmaw)$, equipped with the topology generated by all open sets of $\Sigmaw$ together with the \textit{horocyclic neighborhoods} of the ends, i.e. sets of the form $N=N_0\cup \{\End\}$ where $N_0$ is an open disk in the closed disk model for $\HH^2$ which is tangent to the boundary at $\End$. The action of $\pi_1(\Sigma)$ on $\Sigmaw$ naturally extends to a continuous (neither free nor proper) action on $\Sigmawe$. The quotient of $\Ends(\Sigmaw)$ by this action is precisely the set of punctures of $\Sigma$.
	
	\begin{remark}\label{rem:end_cover_puncture}
		Ends cover the punctures of $\Sigma$ in the sense that the universal cover projection $\Sigmaw\to \Sigma$ admits a continuous extension to a map $\Sigmawe\to \Sigmac$. 
		In particular, a sequence of points $x_n \in \Sigmaw$ converges to an end $\End\in \Ends(\Sigmaw)$ if and only if its projection to $\Sigma$ is a sequence of points converging (in $\Sigmac$) to the puncture covered by $\End$. This happens if and only if $x_n$ eventually enters every horocyclic neighborhood of $\End$.
	\end{remark}
	
	\begin{remark}\label{rmk:two_topologies}
		Notice that $\Sigmaw$ is open and dense in $\Sigmawe$, but this is not the same topology as the one induced from the closed disk model for $\mathbb{H}^2$. Indeed the topology of $\Sigmawe$ is strictly finer; the natural inclusion of $\Sigmawe$ into the closed disk is continuous but not open. Furthermore, the topology induced on the collection of ends is discrete, so $\Sigmawe$ is not compact. Actually not even locally compact, as ends do not have compact neighborhoods.
	\end{remark}
	
	Recall that a \emph{peripheral element} $\delta_x \in \pi_1(\Sigma)$ is the homotopy class of a \emph{peripheral loop} (also denoted by $\delta_x$) around the puncture $x$. If $\End_x$ is an end covering $x$, then $\delta_x$ is a generator of the stabilizer of $\End_x$ in $\pi_1(\Sigma)$. We make the convention that $\delta_x$ is the \emph{positive} peripheral element if the corresponding peripheral loop is \emph{positively} oriented, namely it turns anticlockwise around $x$ (with respect to the orientation of $\Sigma$). This convention is chosen to match the convention that the angle between two circles is also taken in the anticlockwise direction.

	Let $\sigma \in \Moduli(\Sigma)$ be represented by a developing pair $(\dev,\hol)$. We say that $\sigma$ is:
	\begin{itemize}
		\item \emph{Tame at a puncture $x$}: if $\dev$ admits a continuous extension
		$$
		(\devt)_x :\Sigmaw \cup \Ends_x(\Sigmaw) \rightarrow \cp.
		$$ 
		\item \emph{Tame}: if $\dev$ admits a continuous extension
		$$
		\devt : \Sigmaw \cup \Ends(\Sigmaw) \rightarrow \cp.
		$$
		Note that this is equivalent to $\sigma$ being tame at each puncture.
		\item \emph{Relatively elliptic}: if the holonomy representation is relatively elliptic, i.e. the holonomy of every peripheral element is an elliptic M\"obius transformation.
		\item \emph{Degenerate}: if the holonomy representation is degenerate in the sense of \cite[Definition 2.4]{G19}, i.e. if either one of the following happens:
		\begin{itemize}
			\item there are two points $p_\pm \in \cp$ such that the entire holonomy preserves the set $\{p_\pm\}$ and the holonomy of every peripheral element fixes $p_\pm$ individually;
			\item there exists a point $p\in \cp$ such that the entire holonomy fixes $p$ and the holonomy of every peripheral element is parabolic or identity.
		\end{itemize}

	\end{itemize}
	The property of being degenerate is related (but not equivalent) to the more classical notions of \textit{reducible} or \textit{elementary} representations. In the case of punctured spheres, a degenerate representation is always reducible; on the other hand a representation generated by rotations of the Euclidean plane around different points is reducible but non-degenerate (see \cite[\S2.4]{G19}  for a discussion).
	
	The above notions are invariant under conjugation of representations in $\pslc$ and post-composition of developing maps by M\"obius transformations, thus they do not depend on the choice of representative pair $(\dev,\hol)$. The \emph{deformation space of $\cp$--structure on $\Sigma$ which are tame, relatively elliptic and non-degenerate} is $\Modulit(\Sigma)$. The image of $\Modulit(\Sigma)$ under the holonomy map is $\Chart(\Sigma) :=\Hol(\Modulit(\Sigma))$.

	\begin{lemma}\label{lem:devt_is_equivariant}
		Let $\sigma \in \Moduli(\Sigma)$ and let $(\dev,\hol)$ be a developing pair. Let $x$ be a puncture and suppose that $\sigma$ is tame at $x$. Let $\End_x$ be an end covering $x$ and let $\delta_x \in \pi_1(S)$ be a peripheral element fixing it. Then:
		\begin{enumerate}
			\item the map $(\devt)_x$ is $\hol$--equivariant. In particular, the transformation $\hol(\delta_x)$ fixes $\devt(\End_x)$;
			
			\item the transformation $\hol(\delta_x)$ is either trivial, parabolic or elliptic.
			
		\end{enumerate}
	\end{lemma}
	\proof $\ $	
	\begin{enumerate}
		\item Follows by equivariance of $\dev$ and continuity of the extension $\devt$.
		
		\item Let $p:=(\devt)_x(\End_x)$ be one of the fixed points of $\hol(\delta_x)$, and assume by contradiction that $\hol(\delta_x)$ is hyperbolic or loxodromic. Then it  has another fixed point $q$ and there is a $\hol(\delta_x)$--invariant simple arc $\ell$ joining them. Let $\eta$ be an initial segment of $\ell$ starting at $p$ and ending at some other point $y$ on $\ell$, and lift it to an arc $\widetilde{\eta}$ starting at $\End_x$. Consider the family of arcs $\widetilde{\eta}_n := \delta_x^n \cdot \widetilde{\eta}$, for $n\in \ZZ$. Up to replacing $\delta_x$ with its inverse, the sequence   $\{(\devt)_x(\widetilde{\eta}_n)\}$ converges to the whole curve $\ell$ as $n\to + \infty$, and  shrinks to $p$ as $n\to - \infty$. Hence for all $n \in \ZZ^+$ there is a point $x_n\in \widetilde{\eta}_n$ developing to $y$. Then we have $x_n\to \End_x$ in the topology of $\Sigmawe$, but also $(\devt)_x(x_n)=y \not= p$, which contradicts the continuity of $(\devt)_x$ at $\End_x$. 
	\end{enumerate}
	\endproof
	
	We will see in \S\ref{sec:complex}  that if $\sigma\in \Modulit(\Sigma)$ then $\sigma \in \Modulic(\Sigma)$, i.e. the underlying complex structure is that of a punctured Riemann surface. More precisely, $\sigma$ can be defined by a suitable meromorphic quadratic differential with double poles
	(cf. Theorem~\hyperref[thm:tamerelell=merodoublenonint]{E}). However $\Modulit(\Sigma)$ is strictly contained in $\Modulic (\Sigma)$, as the following examples show.
	
	\begin{example}\label{ex:tame_relell}
		
	We now collect examples of structures in $\Modulic(\Sigma)$ which are or are not in $\Modulit(\Sigma)$. These examples show that being tame and having relatively elliptic holonomy are independent concepts.
	\begin{itemize}[leftmargin=0.5cm]
		
	\item All structures induced by Euclidean or hyperbolic metrics with cone points of angles $2\pi \theta$ are in $\Modulit(\Sigma)$, when $\theta \not \in \NN$. For spherical metrics one has to additionally require that they do not have coaxial holonomy (see \cite{MP16}).
	
	\item The structure induced by a complete hyperbolic metric of finite area is tame, but its holonomy is not relatively elliptic because peripherals have parabolic holonomy. Hence it is in $\Modulic(\Sigma)$ but not in $\Modulit(\Sigma)$.
	
	\item Let $\sigma_0$ be the structure induced by a constant curvature metric with cone points of angles $2\pi \theta$, for $\theta \not \in \NN$.
	Remove disks centered at the cones, turn them into crowns and perform infinitely many graftings along arcs joining the crown tips. The resulting structure is in $\Modulic(\Sigma)$ and has relative elliptic holonomy, but it is not tame, hence it is not in $\Modulit(\Sigma)$. This construction is described in \cite{GM19}, where it is shown that these structures arise from meromorphic quadratic differentials with poles of order at least $3$ on punctured Riemann surfaces. Compare Example~\ref{ex:non tame crowns}.
	
	\item Let $\sigma_0\in \Moduli(\Sigmac) $ be the complex projective structure induced by a hyperbolic metric on the closed surface $\Sigmac$. Pick a simple closed geodesic and let $\sigma_n$ be the structure obtained by grafting along it $n$ times. For $n\to \infty$ we obtain a punctured surface $\Sigma$ with two punctures (possibly disconnected if the geodesic is separating) which is endowed with a complex projective structure in $\Modulic(\Sigma)$ (see \cite{HE11}). However it is not tame, and peripherals have hyperbolic holonomy, so it is not in  $\Modulit(\Sigma)$. Compare Example~\ref{ex:non tame pinch geodesics}.
	
	\end{itemize}
	\end{example}

    We conclude this section by observing that structures in $\Modulit(\Sigma)$ carry some additional piece of information which can be regarded as a decoration of the holonomy representation.
    A \textit{framing} for a representation $\rho:\pi_1(\Sigma)\to \pslc$ consists of a choice of a fixed point in $\cp$ for the holonomy about each puncture (compare \cite{AB18,G19}). When considering representations up to conjugacy (as we do), a framing can equivalently be defined as a $\rho$--equivariant map $\mathcal F:\Ends(\Sigmaw) \to \cp$ from the space of ends to $\cp$.
    A framing is said to be \textit{degenerate} if one of the following occurs (compare \cite[\S2.5]{G19}):
    \begin{itemize}
        \item $\mathcal F(\Ends(\Sigmaw))$ consists of two points, preserved as a set by every element, and fixed individually by the holonomy at every puncture;
        \item $\mathcal F(\Ends(\Sigmaw))$ consists of one point, fixed by every element, and the holonomy at every puncture is either parabolic or the identity.
    \end{itemize}
    Every framing of a every non-degenerate representation is non-degenerate (cf. \cite[Prop.3.1]{G19}). In general, a $\cp$--structure can be framed in different ways, by arbitrarily picking the fixed point for each peripheral curve. However, tame structures can be canonically framed.
    
    \begin{corollary}\label{cor:framed_hol}
        Let $\sigma \in \Modulit(\Sigma)$. Then the extension of a developing map provides a non-degenerate canonical framing for the holonomy.
    \end{corollary}
    \proof
        Let $(\dev,\hol)$ be a developing pair defining $\sigma$. By Lemma \ref{lem:devt_is_equivariant} we know $\dev$ extends naturally to  a map $\devt$ on the space of ends. The restriction $\mathcal F =\restr{\devt}{\Ends(\Sigmaw)}$ provides the desired framing. 
        The framing is non-degenerated because $\hol$ itself is a non-degenerate representation.
    \endproof
    
    In the following, whenever dealing with a structure $\sigma \in \Modulit(\Sigma)$, we assume that this natural framing $\mathcal F$ has been chosen for its holonomy representation, and refer to the pair $(\hol,\mathcal F)$ as its \textit{framed holonomy}. 

    \vspace{.5cm}
    In this paper we are mostly interested in a surgery that can be used to deform $\cp$--structures and explore their moduli space. It was introduced by  Maskit (see \cite{MAS69}) and later developed in unpublished work of Thurston (see \cite{KT92,DU,BA19} and references therein for some accounts). The specific version we are interested in is designed to create new structures from old ones without changing their holonomy. For convenience we define it just in the setting of $\cp$--structures in $\Modulit(\Sigma)$.
	
	Let $\sigma \in \Modulit(\Sigma)$, and let $\eta:I\to \Sigma$ be an ideal arc (i.e. with endpoints in the set of punctures). We say $\eta$ is \textit{graftable} if it is simple and injectively developed, i.e  $\devt$ is injective on some (every) lift of $\eta$ to $\Sigmawe$, all the way to the ends. In particular, the two endpoints develop to two distinct points.
	When $\eta$ is graftable, the developed image of any of its lifts $\devt(\widetilde{\eta})$ is a simple arc in $\cp$, hence $\cp\setminus \devt(\widetilde{\eta})$ is a topological disk, endowed with a natural $\cp$--structure, which we call a \textit{grafting region}.
	
\begin{figure}[t]
    	\centering
    	\includegraphics[width=0.9\textwidth]{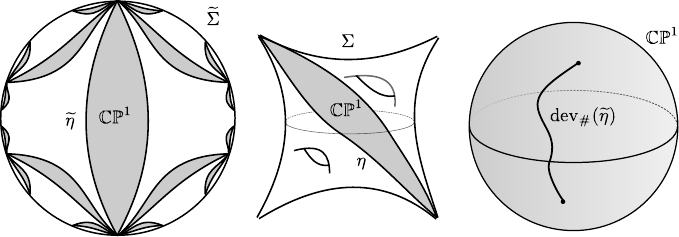}
    	\caption{The structure on $\Sigmaw$ induced by grafting a structure along a curve $\eta$.}
    	\label{fig:grafting}
	\end{figure}

	Let $\sigma \in \Modulit(\Sigma)$ and let $\eta:I\to \Sigma$ be a graftable arc. The \textit{grafting} of $\sigma$ along $\eta$ is the $\cp$--structure $\Gr(\sigma,\eta)$ obtained by the following procedure: for each lift $\widetilde{\eta}$ of $\eta$ to the universal cover, cut $\Sigmaw$ along $\widetilde{\eta}$ and glue in a copy of the disk $\cp\setminus \devt(\widetilde{\eta})$ using $\devt$ as a gluing map. The obvious inverse operation is called \textit{degrafting}.
	The structure on $\Sigmaw$ induced by $\Gr(\sigma,\eta)$ looks like the union of the one induced by $\sigma$ together with an equivariant collection of grafting regions, glued along all the possible lifts of $\eta$ (cf. Figure~\ref{fig:grafting}).
	
	\begin{remark}\label{rem:multiple_grafting}
	If two graftable arcs $\eta,\eta'$ have the same endpoints and are isotopic  through graftable curves, then $\Gr(\sigma,\eta)=\Gr(\sigma,\eta')$, and $\Gr(\sigma,\eta)$ is graftable again along $\eta$ (see \cite[Lemma 2.8]{CDF} or \cite[\S 2]{R19a} for details).
	On the other hand, if $\eta,\eta'$ are disjoint, then $\Gr(\sigma,\eta)$ (resp. $\Gr(\sigma,\eta')$) is graftable along $\eta'$ (resp. $\eta$), and $\Gr(\Gr(\sigma,\eta),\eta') = \Gr(\Gr(\sigma,\eta'),\eta)$. 
	\end{remark}
	
	More generally, a grafting surgery can be defined along any graftable measured lamination on a $\cp$--structure, and the reader familiar with grafting deformations will identify the type of grafting introduced here as a type of projective $2\pi$--grafting  (see \cite{KT92,DU,BA19} for details).  
    We record   the following statement for future reference.
    
	\begin{lemma}\label{lem:grafting_properties}
		Let $\sigma \in \Modulit(\Sigma)$ and $\eta:I\to \Sigma$ be a graftable arc. Then
		\begin{enumerate}
		    \item $\Hol(\Gr(\sigma,\eta))=\Hol(\sigma)$ (i.e. grafting preserves the holonomy),
		    \item $\Gr(\sigma,\eta)\in \Modulit(\Sigma)$,
		    \item grafting does not change the developed images of the punctures (i.e. grafting preserves the framed holonomy).
		\end{enumerate}
	\end{lemma}
    \proof
    The first statement is well-known in the literature for this type of grafting (see for instance \cite{BA19} and references therein). The statements about tameness and framing follow by pasting together the developing map for $\sigma$ and the natural embedding of the grafting regions in $\cp$.
    \endproof

	%%%%%%%%%%%%%%%%%%%%%%%%%%%%%%%%%
	%%%%%%%%%%%%%%%%%%%%%%%%%%%%%%%%%
	%%%%%%%%%%%%%%%%%%%%%%%%%%%%%%%%%
	
	\subsection{The M\"obius completion}\label{subsec:mobius_completion}
	
	In this section we prove Theorem~\hyperref[thm:embedding]{C}. Henceforth we fix a complex projective structure $\sigma \in \Moduli(\Sigma)$ with developing pair $(\dev,\rho)$. First of all we recall the definition of a natural projective completion of $\Sigmaw$ defined in terms of $\sigma$ (see \cite{KP94} for details).
	Let $g_0$ be a conformal Riemannian metric on $\cp$ (e.g. the standard spherical metric). Let $g:=\dev^*(g_0)$ 
	be the metric on $\Sigmaw$ obtained by pullback, and let $d$ be the associated distance function, i.e.
	$$d (x,y):= \inf\{ \ell_{g}(\eta) \ | \ \eta:[0,1]\to \Sigmaw \textrm{ is a rectifiable arc from $x$ to $y$} \} $$
	where $\ell_{g}(\eta)$ denotes the length of $\gamma$ with respect to the metric $g$. Notice that $g$ is generally not invariant under deck transformations.
	By construction $(\Sigmaw,d)$ is a path-connected length space.	It is locally path-connected, but not necessarily geodesic. Moreover it is locally compact, but in general not proper, nor complete.
	
	The \textit{M\"obius completion} $\Compl{\sigma}{\Sigmaw}$ of $\Sigmaw$ with respect to $\sigma$ is defined to be the metric completion of $(\Sigmaw,d)$. The subspace $\Binf{\sigma}{\Sigmaw}:=\Compl{\sigma}{\Sigmaw}\setminus \Sigmaw$ is called the \textit{ideal boundary} of $\Sigmaw$ with respect to $\sigma$.
	We collect the following facts from \cite[\S 2]{KP94}:
	\begin{enumerate}
	    \item  Different choices of the metric $g_0$ on $\cp$ or of the developing map for $\sigma$ result in   metrics on $\Sigmaw$ having the same underlying uniform structure. So $\Compl{\sigma}{\Sigmaw}$ does not depend (up to homeomorphism) on these choices.
	    
	    \item $\dev:\Sigmaw\to\cp$ extends continuously to a map $\devc{\sigma}: \Compl{\sigma}{\Sigmaw}\to \cp$.
	    
	    \item The action of  $\pi_1(\Sigma)$ by deck transformations extends to an action   by homeomorphisms on the M\"obius completion.
	\end{enumerate}
	
	\begin{lemma}\label{lem:equivariance_mobius_compl}
	The map $\devc{\sigma}$ is $\rho$--equivariant.
	\end{lemma}
	\proof
	Let $\xi \in \Binf{\sigma}{\Sigmaw}$ and let $x_n \in \Sigmaw$ a Cauchy sequence converging to $\xi$. Then by continuity of $\devc{\sigma}$
	\begin{align*}
	\devc{\sigma}(\gamma \cdot \xi) &= \devc{\sigma}(\gamma \cdot \lim_{n \to \infty} x_n) = \lim_{n \to \infty} \devc{\sigma}(\gamma \cdot  x_n)\\
	&= \rho(\gamma) \cdot \lim_{n \to \infty} \devc{\sigma}( x_n) = \rho(\gamma) \cdot\devc{\sigma}(\xi).
	\end{align*}
	\endproof
	
	\begin{lemma}\label{lem:completion_is_length}
	$\Compl{\sigma}{\Sigmaw}$ is a complete, path-connected and locally path-connected length space. 
	\end{lemma}
	\proof
	Completeness is trivial by construction. The completion of a length space is a length space (see for instance \cite[I.3.6(3)]{BH99}).
	Since $\Sigmaw$ is path-connected and $\Compl{\sigma}{\Sigmaw}$ is a length space, it follows that $\Compl{\sigma}{\Sigmaw}$ is path-connected. Analogously one can obtain that $\Compl{\sigma}{\Sigmaw}$ is locally path-connected.
	\endproof

	The following examples describe more explicitly the M\"obius completion for projective structures defined by certain constant curvature metrics. Notice they are both examples of hyperbolic M\"obius structures with respect to the terminology introduced in \cite[\S 2]{KP94}.
	
	\begin{example}\label{ex:completion_hyperbolic}
	    Let $\sigma=(\dev,\rho)$ be defined by a complete hyperbolic metric of finite area on $\Sigma$. 
	    In this case $\Compl{\sigma}{\Sigmaw}$ is homeomorphic to a closed disk, and $\Binf{\sigma}{\Sigmaw}$ to a circle. Ideal points are either ends, or limit points of complete lifts of closed geodesics.
	    Indeed, $\rho:\pi_1(\Sigma)\to \pslr$ is an isomorphism onto Fuchsian group, and $\dev:\Sigmaw \to \cp$ is a $\rho$--equivariant diffeomorphism with an open hemisphere.  
	\end{example}
	
	\begin{example}\label{ex:completion_spherical}
	    Let $\sigma=(\dev,\rho)$ be defined by a spherical metric on $\Sigma$, with cone singularities at the punctures. 
	    In this case $\Compl{\sigma}{\Sigmaw}=\Sigmawe$ is homeomorphic to the end-extension, and $\Binf{\sigma}{\Sigmaw}=\Ends (\Sigmaw)$.
	    Indeed, the action of $\pi_1(\Sigma)$ on $\Sigmaw$ preserves a spherical metric and admits a fundamental domain $D$ given by a geodesic spherical polygon having finite area $A$ and all the vertices in the set of ends. Notice that each pair of non-intersecting edges of this polygon has positive finite distance, and let $L>0$ be the minimum of such distances.
	    Pick $\xi \in \Binf{\sigma}{\Sigmaw}$, and a rectifiable curve of finite length $\gamma:[0,1)\to \Compl{\sigma}{\Sigmaw}$ tending to $\xi$. 
	    If $\gamma$ intersects finitely many fundamental domains, then it is eventually contained in a single one, hence $\xi$ must be an end. 
	    If $\gamma$ intersects infinitely many domains $D_n$, then the length of the arcs $\gamma\cap D_n$ converges to zero, so is eventually less than $L$. In particular, eventually all the domains $D_n$ share a common vertex. By construction this vertex is an end  and $\gamma$ converges to it, which forces $\xi$ to be an end.
	\end{example}

	\begin{lemma}\label{lem:ideal_points_are_accessible}
	For all $x\in \Sigmaw$, $\xi \in \Binf{\sigma}{\Sigmaw}$ and $c >0$ there is a continuous curve $\eta_c:[0,1)\to \Sigmaw$ such that $\eta_c(0)=x$, $\lim_{t\to 1} \eta_c(t)=\xi$ and $d(x,\xi)\leq \ell(\eta_c)\leq d(x,\xi)+c$.
	\end{lemma}
	
	\proof
    By definition $d(x,\xi) = \lim_{n \to \infty} d(x,y_n)$ for any Cauchy sequence $\{y_n\}$ converging to $\xi$. Let $\{y_n\}$ be a Cauchy sequence in $\Sigmaw$ converging to $\xi$ such that
    $$
    d(x,y_i) \leq  d(x,\xi) + \frac{1}{i}, \qquad \text{and} \qquad d(y_i,y_{i+1}) \leq \frac{1}{i^2}.
    $$
    Such sequence can be easily constructed from any Cauchy sequence by taking an appropriate subsequence. 
    Since $\Sigmaw$ is a length space, for all $k$ there is a continuous curve $\gamma_k: [0,1] \rightarrow \Sigmaw$ such that $\gamma_k(0) = y_k$, $\gamma_k(1) = y_{k+1}$ and
    $$
    \ell(\gamma_k) \leq  d(y_k,y_{k+1}) + \frac{1}{k^2} = \frac{2}{k^2}.
    $$
    By concatenating these curves, for every $i$, we obtain a continuous curve $\gamma_i : [0,1) \rightarrow \Sigmaw$ such that $\gamma_i(0) = y_i$, $\lim_{t\to 1} \gamma_i(t) = \xi$, and
    $$
    \ell(\gamma_i) \leq \sum_{k=i}^\infty \frac{2}{k^2} =: T_i.
    $$
    In particular, $\lim_{i \to \infty} \ell(\gamma_i) = \lim_{i \to \infty} T_i = 0$.
    Finally, let $\eta_i : [0,1) \rightarrow \Sigmaw$ be a continuous curve such that $\eta_i(0) = x$, $\eta_i(1) = y_i$, and
    $$
    \ell(\eta_i) \leq d(x,y_i) + \frac{1}{i^2}.
    $$
    Let $f_i : [0,1) \rightarrow \Sigmaw$ be the continuous curve obtained by concatenating $\eta_i$ with $\gamma_i$. Then $f_i$ is a continuous curve such that $f_i(0) = x$, $\lim_{t\to 1} f_i(t)=\xi$ and
    $$
    d(x,\xi)\leq \ell(f_i) = \ell(\eta_i) + \ell(\gamma_i) \leq d(x,y_i) + \frac{1}{i^2} + T_i \leq d(x,\xi) + \frac{1}{i} + \frac{1}{i^2} + T_i.
    $$
    Now let $i$ such that $\frac{1}{i} + \frac{1}{i^2} + T_i < c$ and take $f_c := f_i$.
	\endproof

	\begin{lemma}\label{lem:ideal_points_locally_pathconnected}
	Let $\xi \in \Binf{\sigma}{\Sigmaw}$ and $\varepsilon >0$. Then $\Ballc{\sigma}(\xi,\varepsilon)\cap \Sigmaw$ is path-connected.
	\end{lemma}
	\proof
	First of all let us show that each path-component $N$ of $\Ballc{\sigma}(\xi,\varepsilon)\cap \Sigmaw$ contains points arbitrarily close to $\xi$. Pick a base point $x\in N$, and let $R=d(x,\xi)$; notice $R<\varepsilon$. By Lemma~\ref{lem:ideal_points_are_accessible} for all $c>0 $ we can pick a  continuous curve $\eta_c:[0,1)\to \Sigmaw$ such that $\eta_c(0)=x$, $\lim_{t\to 1} \eta_c(t)=\xi$ and $R\leq \ell(\eta_c)\leq R+c$. For each $t\in[0,1)$ we have
	$$
	d(\eta_c(t),\xi)\leq \ell(\eta_c([t,1))) \leq \ell(\eta_c([0,1))) \leq R+c.
	$$
	In particular, for $c<\frac{1}{2}(\varepsilon-R)$ we get that $d(\eta_c(t),\xi)<\varepsilon$, i.e. $\eta_c$ is entirely contained in $\Ballc{\sigma}(\xi,\varepsilon)\cap \Sigmaw$. Since it is a curve starting at $x$, it is then entirely contained in $N$; since it converges to $\xi$ we get $\lim_{t\to 1}d(\eta_c(t),\xi)=0$.\par
	
	Suppose by contradiction that $\Ballc{\sigma}(\xi,\varepsilon)\cap \Sigmaw$ admits at least two different path-components $N_1,N_2$. 
% 	Since $\Compl{\sigma}{\Sigmaw} $ is locally path-connected, we can shrink $\varepsilon>0$ to ensure that $\Ballc{\sigma}(\xi,\varepsilon)$ is path-connected.
	Let $x_k\in N_k$ be two points such that $d(x_k,\xi)<\frac{\varepsilon}{4}$.
	In particular, $d(x_1,x_2)<\frac{\varepsilon}{2}$. Since $(\Sigmaw,d)$ is a length space, for every $\delta>0$ we can find a continuous curve $\gamma_\delta:[0,1]\to \Sigmaw$ joining $x_1$ to $x_2$ of length at most $\frac{\varepsilon}{2}+\delta$.
	Let now $z\in \gamma_\delta$. Without loss of generality let us assume that $d(z,x_1)\leq d(z,x_2)$, so that by triangle inequality we get
	$$
	d(z,\xi)\leq   d(z,x_1)  + d(x_1,\xi) \leq  \frac{1}{2}\left(\frac{\varepsilon}{2}+\delta\right) + \frac{\varepsilon}{4} = \frac{\varepsilon}{2}+\frac{\delta}{4}.
	$$
	In particular, for each $\delta<\varepsilon$ we get that the curve $\eta_\delta$ is at distance at most $\varepsilon$ from $\xi$. In particular, it is entirely contained in $\Ballc{\sigma}(\xi,\varepsilon)\cap \Sigmaw$, which   contradicts the fact that $x_1,x_2$ are in distinct path-components.
	
	\endproof
	
	Our next goal is to define a  cyclic order on $\Binf{\sigma}{\Sigmaw}$, which will induce a total order on $\Binf{\sigma}{\Sigmaw} \setminus \{\xi\}$, for any $\xi \in \Ends(\Sigmaw)$.
	
	\begin{lemma}\label{lem:binf_decomposition}
	For any pair of distinct points $(\xi_0,\xi_1)\in \Binf{\sigma}{\Sigmaw}$ there exists a simple continuous curve $\gamma:(0,1)\to \Sigmaw$ such that $\lim_{t\to 0}\gamma(t)=\xi_0$, $\lim_{t\to 1}\gamma(t)=\xi_1$. 

	Moreover, for any such curve $\gamma$, the space $\Compl{\sigma}{\Sigmaw} \setminus \Clo(\gamma)$ has exactly two path-components, which we call the left and right components $C_L(\gamma),C_R(\gamma)$ with respect to the orientation of $\gamma$. The induced partition of $\Binf{\sigma}{\Sigmaw}$ as
	$$
	\{\xi_0,\xi_1\}\cup(\Binf{\sigma}{\Sigmaw}\cap C_L(\gamma))\cup (\Binf{\sigma}{\Sigmaw}\cap C_R(\gamma))
	$$
	only depends on the ordered pair $(\xi_0,\xi_1)$ and not on $\gamma$.
	\end{lemma}
	\proof
	Existence of $\gamma$ is clear, for instance by Lemma~\ref{lem:ideal_points_are_accessible}. Let us show that its complement consists of exactly two path-components. $\Sigmaw \setminus \gamma$ clearly has exactly two path components, so $\Compl{\sigma}{\Sigmaw} \setminus \Clo(\gamma)$ has at most two components (again by Lemma~\ref{lem:ideal_points_are_accessible}). We need to show that no ideal point can be joined by an arc to both components. This follows from Lemma~\ref{lem:ideal_points_locally_pathconnected}.\par
	
	To show that the induced decomposition of $\Binf{\sigma}{\Sigmaw}$ does not depend on the choice of $\gamma$, just notice that any two such curves are isotopic relatively to their endpoints in $\Sigmaw$.
	\endproof

	Hence we denote by $C_L(\xi_0,\xi_1):=\Binf{\sigma}{\Sigmaw}\cap C_L(\gamma)$ and $C_R(\xi_0,\xi_1)=\Binf{\sigma}{\Sigmaw}\cap C_R(\gamma)$ for any curve $\gamma$ as in Lemma~\ref{lem:binf_decomposition}.
	We define the following ternary relation on $\Binf{\sigma}{\Sigmaw}$. If $\xi_0,\xi_1,\zeta\in \Binf{\sigma}{\Sigmaw}$ then we say they are in relation (denoted $[\xi_0,\zeta,\xi_1]$) if  $\zeta\in C_R(\xi_0,\xi_1)$, i.e. $\zeta$ is on the right of $\gamma$. 
	
	\begin{remark}\label{rem:invariant_cyclic_order}
	    This relation defines a $\pi_1(\Sigma)$--invariant cyclic order on $\Binf{\sigma}{\Sigmaw}$.
	\end{remark}

    The goal of  the rest of this section is to explore the features of the M\"obius completion and the ideal boundary in the case of structures from $\Modulit(\Sigma)$.

	\begin{proposition}\label{prop:tame_iff_extension}
	A structure $\sigma$ is tame if and only if the natural embedding $j_\sigma:\Sigmaw\hookrightarrow \Compl{\sigma}{\Sigmaw}$ extends to a $\pi_1(\Sigma)$--equivariant continuous embedding $j^{\#}_\sigma:\Sigmawe \hookrightarrow \Compl{\sigma}{\Sigmaw}$. 
	Moreover in this case $\devt=\devc{\sigma} \circ j_\sigma^{\#}$.
	
	\end{proposition}
	\proof
	First assume the existence of a $\pi_1(\Sigma)$--equivariant continuous embedding $j_\sigma^{\#}:\Sigmawe \hookrightarrow \Compl{\sigma}{\Sigmaw}$. As remarked above there exists a continuous extension $\devc{\sigma}$ of $\dev$ to $\Compl{\sigma}{\Sigmaw}$. Then $\devc{\sigma} \circ j^\sigma$ provides a continuous extension of $\dev$ to $\Sigmawe$, i.e. $\sigma$ is tame.
	
	Conversely let $\sigma$ be tame, let $\End \in \Ends(\Sigmaw)$ and $p_\End=\devt(\End)$. Since $\dev$ extends continuously to $\End$, for all $\varepsilon>0$ the set $N_\varepsilon=(\devt)^{-1}(B(p_\End,\varepsilon))$ is an open neighborhood of $\End$ in $\Sigmawe$, containing points at distance at most $\varepsilon$ from $\End$. Therefore we can construct a Cauchy sequence $x_n$ in $\Sigmaw$ converging to $\End$ (in $\Sigmawe$). We can associate to $E$ the limit of $x_n$ in the completion $\Compl{\sigma}{\Sigmaw}$. Suppose $y_n$ is another Cauchy sequence in $\Sigmaw$ converging to $\End$ (in $\Sigmawe$). By definition of the topology on $\Sigmaw$, continuity of $\devt$ at $\End$ implies that $\devt(x_n)$ and $\devt(y_n)$ both converge to $p_\End$. Hence $y_n$ eventually enters each neighborhood $N_\varepsilon$. As a result we get $d(x_n,y_n)\leq 2\varepsilon$, which implies that the two sequences give rise to the same point in the completion. This defines the desired extension, which is (sequentially) continuous.
	Injectivity follows from the fact that any two ends are at a positive distance from each other.
	Moreover $\devt=\devc{\sigma} \circ j_\sigma^{\#}$ because they agree on the dense subset $\Sigmaw$ and $\cp$ is Hausdorff.
	\endproof
	
	In particular, tame structures have infinitely many ideal points, hence they are of hyperbolic type with respect to the classification in \cite{KP94}. Moreover it should be noticed that ends do not have compact neighborhoods, so that the completion fails to be locally compact or proper.  %On the other hand $\Sigmaw$ is locally compact but not complete.

	\begin{example}\label{ex:non tame crowns}
	Gupta and Mj in  \cite{GM19} consider structures obtained by grafting crowned hyperbolic surfaces, and show that the local structure at the crown can be modeled by a meromorphic differential with a pole of sufficiently high order.
	For such a structure, every sequence going off to a puncture gives rise to an ideal point in the M\"obius completion, but sequences converging in different Stokes sectors  develop to sequences converging to different limit points in $\cp$, hence give rise to different ideal points in the M\"obius completion. 
	They are not tame structures (as observed in Example \ref{ex:tame_relell}), and the space of ends does not embed continuously in their ideal boundary.
	Notice that Lemma~\ref{lem:ideal_points_locally_pathconnected} applies to each individual ideal point, while  the intersection of $\Sigmaw$ with the neighborhood of an end can fail to be connected.
    \end{example}	
	
	\begin{example}\label{ex:non tame pinch geodesics}
	For a more extreme behavior, take a closed hyperbolic surface, and graft it along a geodesic pants decomposition infinitely many times. The underlying complex structure is being pinched along each pants curve, and in the limit the structure decomposes into a collection of thrice--punctured spheres (see \cite[\S 6]{HE11}). There, punctures do not give rise to well-defined ideal points; indeed, the structure has hyperbolic peripheral holonomy, hence it is not tame (by Lemma~\ref{lem:devt_is_equivariant}).
    \end{example}	
	
	\begin{remark}\label{rmk:parabolic_not_open}
	In general the embedding $j_\sigma^{\#}$ in Proposition~\ref{prop:tame_iff_extension} is not open. For instance consider the tame relatively parabolic structure induced by a complete finite area hyperbolic metric. In this case the completion is the closed disk, and we have already observed in Remark~\ref{rmk:two_topologies} that inclusion of the space of ends in it is not open. We will show below in Proposition~\ref{prop:open_horocycles} that having relatively parabolic holonomy is actually the only obstruction to the openness of $j_\sigma^{\#}$.
	\end{remark}

	%%%
	
	For a point $p\in \Compl{\sigma}{\Sigmaw} $ we define the balls 
	\begin{align*}
	    \Ball (p,r) &:= \{ z\in \Sigmaw \ | \  d(p,z)<r \},\\
	    \Ballt (p,r) &:= \{ z\in \Sigmawe \ | \  d(p,z)<r \},\\
	    \Ballc{\sigma} (p,r) &:= \{ z\in \Compl{\sigma}{\Sigmaw} \ | \  d(p,z)<r \}.
	\end{align*}
    By Proposition~\ref{prop:tame_iff_extension}, $\Ball (p,r) \subseteq \Ballt (p,r) \subseteq \Ballc{\sigma} (p,r)$ for any $p$ and $r$, and these balls are open. For small values of $r$ they also enjoy extra properties.

    By Proposition~\ref{prop:tame_iff_extension} we know we can embed the space of ends in the ideal boundary $\Binf{\sigma}{\Sigmaw}$ of the M\"obius completion $\Compl{\sigma}{\Sigmaw}$. So it makes sense for a given subset $Z$ of $\Sigmaw$ to consider its closure $\Clot(Z)$ in $\Sigmawe$ or $\Cloc{\sigma}(Z)$ in $\Compl{\sigma}{\Sigmaw}$; by completeness of  $\Compl{\sigma}{\Sigmaw}$, the latter is the same as the metric completion of $Z$ with respect to some choice of metric as in the previous sections. In either case, the \emph{(topological) boundary} of a subset $Z$ is the difference between its \emph{closure} and its \emph{interior} $\partial Z:=\Clo(Z) \setminus \Int(Z)$.
    
\begin{lemma}\label{lem:small_balls}
For each $p\in \Sigmaw$ let $R=d(p,\Binf{\sigma}{\Sigmaw})$. Then for all $r<R$   the following hold.
\begin{enumerate}
	    \item $\Ball (p,r)= \Ballt(p,r) =  \Ballc{\sigma} (p,r) $,
	    
	    \item $\Clo(\Ball (p,r)) = \Clot (\Ballt(p,r)) =  \Cloc{\sigma} (\Ballc{\sigma}(p,r)) $  is complete.
	    
\end{enumerate}
\end{lemma}
\proof
Since the metric structure on $\Sigmaw$ is induced by the  Riemannian metric $g_\End$, for sufficiently small radius, the metric balls are just balls for the Riemannian metric $g_\End$. In particular, they are all disjoint from the ideal boundary, hence they coincide and their closure is complete and contained in $\Sigmaw$. 
\endproof

\begin{lemma}\label{lem:isometric_balls}
For each $p\in \Sigmaw$ let $R=d(p,\Binf{\sigma}{\Sigmaw})$; then for all $r\leq R$ the developing map induces an isometry between $\Cloc{\sigma}(\Ballc{\sigma}(p,r))$ and $\Clo (\Ball(\dev(p),r))$.
\end{lemma}
\proof
Let $I$ be the set of $r \in [0,R]$ such that the developing map induces an isometry between $\Cloc{\sigma}(\Ballc{\sigma}(p,r))$ and $\Clo (\Ball(\dev(p),r))$. We are going to show that $I$ is not empty, open on the right and closed on the right to conclude that $I = [0,R]$.
\begin{itemize}
    \item $[0,\epsilon) \subset I$ for $\epsilon >0$ small enough. This is because $\dev$ is a local isometry at $p$.
    
    \item If $[0,r) \subset I$ then $[0,r] \subset I$. Notice that the developing map induces an isometry between $\Cloc{\sigma}(\Ballc{\sigma}(p,r - \frac{1}{n}))$ and $\Clo (\Ball(\dev(p),r - \frac{1}{n}))$ for all $n>0$. This is enough to deduce that the developing map induces an isometry between $\Ballc{\sigma}(p,r)$ and $\Ball(\dev(p),r)$. Since $\Cloc{\sigma}(\Ballc{\sigma}(p,r))$ is complete, and the metric completion is unique, the developing map induces an isometry between $\Cloc{\sigma}(\Ballc{\sigma}(p,r))$ and $\Clo (\Ball(\dev(p),r))$.
    
    \item If $[0,r] \subset I,r<R$ then $[0,r+\epsilon] \subset I$ for $\epsilon >0$ small enough. Given that $r \in I$, then the developing map induces an isometry between $\Cloc{\sigma}(\Ballc{\sigma}(p,r))$ and $\Clo (\Ball(\dev(p),r))$. In particular, $ \partial \Ballc{\sigma}(p,r)$ is compact. Since $r < d(p,\Binf{\sigma}{\Sigmaw})$, there is an $\epsilon$--neighborhood of $\partial \Ballc{\sigma}(p,r)$ on which $\dev$ is an isometry and $r+\epsilon \in I$.
\end{itemize}
\endproof

We call $\Cloc{\sigma}(\Ballc{\sigma}(p,R))$ the \textit{maximal ball} centered at $p$.
It is a maximal round ball containing $p$, in the sense of \cite{KP94}. Our goal in \S \ref{subsec:canonic_max_nbhd} is to construct analogous ``round neighborhoods'' of all the ends, in the case of elliptic holonomy.  We will need the following preliminary results.

\begin{proposition}\label{prop:open_horocycles}
    Let $\End\in \Ends(\Sigmaw)$, let $\sigma$ be tame at $\End$, and let $N$ be an open horocyclic neighborhood of $\End$. Then $j_\sigma^{\#}(N)$ is open if and only if $\End$ has non-parabolic holonomy.
    \end{proposition}

	\proof
	Let $\delta_\End$ be the peripheral element fixing $\End$, let $R_\End := \hol(\delta_\End)$, and let $p_\End=\devt(\End)$. 
	By Lemma~\ref{lem:small_balls}, every point in $j_\sigma^{\#}(N) \cap \Sigmaw$ is in the interior of $j_\sigma^{\#}(N)$, so we only need to check whether $\End$ is in the interior of $j_\sigma^{\#}(N)$.\par
	
	First, consider the case $R_\End$ is parabolic. 
	Pick a point $x\in \partial j_\sigma^{\#}(N)$ that does not develop to $p_\End$, e.g. on the image via $j_\sigma^{\#}(N)$ of a horocycle bounding $N$.
	Then $d(\End,\delta_\End^n(x))\to 0$, i.e. $\delta_\End^n(x)\to \End$ in $\Compl{\sigma}{\Sigmaw}$. So the sequence $\delta_\End^n(x)$ must eventually enter in every open neighborhood of $\End$ in $\Compl{\sigma}{\Sigmaw}$. However it clearly does not enter in $j_\sigma^{\#}(N)$ by construction, which shows $j_\sigma^{\#}(N)$ is not open.
	
	So let us now assume $R_\End$ is non--parabolic; by Lemma~\ref{lem:devt_is_equivariant} we know that since $\sigma$ is tame at $\End$, $R_\End$ is either the identity or elliptic.
	Since $\delta_\End$ acts cocompactly on the boundary $\partial N$ of $N$ and $\devt$ is a local diffeomorphism along $\partial N$, we have that $\devt^{-1}(p_\End)\cap \partial N$ is finite in any $\delta_\End$--fundamental domain. In particular, we can equivariantly modify $N$ to a $\delta_\End$--invariant neighborhood $W \subset N$ of $\End$, such that $\partial W$ stays at finite distance from $\partial N$. By construction $\End$ is the only end in the closure of $W$.
	
	When $R_\End$ is trivial or elliptic, the set $\devt(\partial W)$ has compact closure in $\cp \setminus \{p_\End\}$. In particular, it sits in the annulus $\{ z \in \cp \ | \ R_1 \leq d_0(p_\End,z) \leq R_2 \}$, for some suitable radii $0<R_1 \leq R_2$. 
    For $r<R_1$ consider the open $R_\End$--invariant ball $D_r\subseteq \cp$ of radius $r$ around $p_\End$, as well as the open ball $\Ballc{\sigma}(\End,r)$.
	Observe that $\devc{\sigma} (\Ballc{\sigma}(\End,r))$ is contained in $D_r$, and so is disjoint from $\devt(\partial W)$. 
    We claim $\Ballc{\sigma}(\End,r)\subseteq j_\sigma^{\#}(W) \subset j_\sigma^{\#}(N)$. By contradiction let $x\in \Ballc{\sigma}(\End,r)\setminus j_\sigma^{\#}(W)$. Then connect $x$ to $\End$ by a continuous arc $\gamma$ contained in $\Ballc{\sigma}(\End,r)$ (which is possible since we are in a length space). Then $\gamma$ has to cross $\partial j_\sigma^{\#}(W)$, since $\partial W$ separates $\End$ from the complement of $W$ in $\Sigmawe$. Then $\devc{\sigma}(\gamma)$ meets $\devt(\partial W)=\devc{\sigma}(\partial j_\sigma^{\#}(W))$, which leads to the desired contradiction.
	\endproof
	
	\begin{figure}[t]
		\centering
		\includegraphics[width=.6\textwidth]{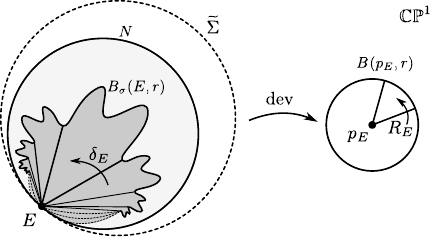}
		\caption{A horocycle containing a ball, in the elliptic case.}
		\label{fig:canonica_max_nbhd}
	\end{figure}	
	
	We summarize the results of this section in the following statement.
	
	\begin{maintheoremc}{C}\label{thm:embedding}
    Let $\sigma \in \Moduli (\Sigma)$ be non-degenerate and without apparent singularities. Let $j^{\#}:\Sigmaw \to \Sigmawe$ and $j_\sigma:\Sigmaw \to \Compl{\sigma}{\Sigmaw}$ be the  natural embeddings. Then $\sigma \in \Modulit(\Sigma)$ if and only if there exists a continuous open $\pi_1(\Sigma)$--equivariant embedding $j^{\#}_\sigma:\Sigmawe \to \Compl{\sigma}{\Sigmaw}$
    that makes the following diagram commute
    
    \centering
    \begin{tikzpicture}[scale=0.8]
    \node (S) at (-2,0) {$\Sigmaw$};
    \node (CP1) at (2,0) {$\cp$};
    \node (E) at (0,1) {$\Sigmawe$ };
    \node (M) at (0,-1) {$\Compl{\sigma}{\Sigmaw}$};
    \path[->,font=\scriptsize,>=angle 90]
    (S) edge node[above]{$j^{\#}$}    (E)
    (S) edge node[below]{$j_\sigma$}    (M)
    (E) edge node[above]{$\devt$}    (CP1)
    (M) edge node[below]{$\devc{\sigma}$}    (CP1)
    (E) edge node[left]{$j^{\#}_\sigma$}    (M);
    \end{tikzpicture}
    \end{maintheoremc}
    \proof
    First assume $\sigma\in \Modulit(\Sigma)$. Since $\sigma$ is tame, by Proposition~\ref{prop:tame_iff_extension} we know that $j_\sigma:\Sigmaw\hookrightarrow \Compl{\sigma}{\Sigmaw}$ extends to a $\pi_1(\Sigma)$--equivariant continuous embedding $j^{\#}_\sigma:\Sigmawe \hookrightarrow \Compl{\sigma}{\Sigmaw}$, and that $\devt=\devc{\sigma} \circ j_\sigma^{\#}$. 
    To check that $j^{\#}_\sigma$ is open we argue as follows. Observe that the restriction of $j_\sigma^{\#}$ to $\Sigmaw$ is just the natural embedding of $\Sigmaw$ in its completion, which is open. So we only need to check the ends. Let $\End$ be an end; without loss of generality we can assume that an open neighborhood of $\End$ in $\Sigmawe$ is an open horocycle $N$. Since $\sigma$ is relatively elliptic, Proposition~\ref{prop:open_horocycles} implies that $j^{\#}_\sigma(N)$ is an open neighborhood of $j^{\#}_\sigma(\End)$ in $\Compl{\sigma}{\Sigmaw}$.

    Conversely, assume the existence of the extension $j_\sigma^{\#}$ as in the statement. Its continuity implies tameness of $\sigma$ by Proposition~\ref{prop:tame_iff_extension}. 
    Let $\End$ be an end.
    By Lemma~\ref{lem:devt_is_equivariant} we know that the holonomy of $\sigma$ at $\End$ is either trivial, parabolic or elliptic. The first case is excluded by the hypothesis that $\sigma$ has no apparent singularities, and the second case by the hypothesis that $j_\sigma^{\#}$ is open, together with Proposition~\ref{prop:open_horocycles}. Therefore $\sigma$ is relatively elliptic. It is also assumed to be non-degenerate, hence we can conclude that $\sigma\in \Modulit(\Sigma)$.
    \endproof

    \begin{corollary}\label{cor:ends are discrete}
       If $\sigma \in \Moduli(\Sigma)$, then $\Ends(\Sigmaw)$ is a discrete subspace of $\Binf{\sigma}{\Sigmaw}$.
    \end{corollary}
    \proof
    Let $\End$ be an end.
    By Theorem~\hyperref[thm:embedding]{C} any horocyclic neighborhood $N$ of $\End$ is open for the topology of $\Compl{\sigma}{\Sigmaw}$, under the natural embedding $j_\sigma^{\#}$. By definition, $N$ does not contain any other point of $\Binf{\sigma}{\Sigmaw}$, hence $\End$ is an open point.
    \endproof    
	
	\begin{corollary}\label{cor:proper_free_action}
	Let $\sigma$ be tame and relatively elliptic.
	For every end $\End \in \Ends(\Sigmaw)$, the action of the peripheral subgroup $\langle \delta_\End \rangle$ on $\Compl{\sigma}{\Sigmaw} \setminus \{ \End \}$ is proper and free.
	\end{corollary}
	\proof
    The action on the M\"obius completion extends the action by deck transformations, so the statement is trivial for points in $\Sigmaw$. By Proposition~\ref{prop:open_horocycles}, both metric balls and horocyclic neighborhoods provide fundamental systems of neighborhoods of the ends in the completion. So one can see that the action of $\delta_\End$ on the subspace $\Ends(\Sigmaw) \setminus \End$ is proper and free. The case of a general ideal point follows from this fact together with the existence of an $\delta_\End$--invariant cyclic order on the ideal boundary (cf. Remark~\ref{rem:invariant_cyclic_order}).
	\endproof

	%%%%%%%%%%%%%%%%%%%%%%%%%%%%%%%%%
	%%%%%%%%%%%%%%%%%%%%%%%%%%%%%%%%%
	%%%%%%%%%%%%%%%%%%%%%%%%%%%%%%%%%
		
	\subsection{Local properties of the developing map at an end}
	\label{subsec:canonic_max_nbhd} 
    The main goal of this section is to prove Theorem~\hyperref[thm:maxneighborhood]{D}, about the behavior of developing maps around $\End$ for a structure $\sigma \in \Modulit(\Sigma)$. 
    If $\sigma$ has developing pair $(\dev,\rho)$, and if $\End\in \Ends(\Sigmaw)$, then let $p_{\End}:=\devt(\End)\in \cp$ and let $\delta_{\End}\in\pi_1(\Sigma)$ be a peripheral element fixing $\End$. 
    Then $R_{\End}:=\hol(\delta_{\End})$ is an elliptic M\"obius transformation fixing $p_{\End}$ (cf.~Lemma~\ref{lem:devt_is_equivariant}); let $q_\End$ denote the other fixed point of $R_\End$. 
	We will construct a family of $\delta_\End$--invariant neighborhoods of $\End$ which develop to $R_\End$--invariant round disks in $\cp$, and on which $\devt$ restricts to a branched covering (branching only at $\End$).
	
	While the results of the previous sections relied (but did not depend), on the choice of the background metric $g_0$ on $\cp$, we now want to exploit the fact that the peripheral holonomy is elliptic to pick a convenient metric. The topological structure of the M\"obius completion is not affected by this (e.g. ideal points, etc), but finer metric statements (e.g. the shape and properties of individual metric balls) are.
	Let $g_0$ be the unique $R_\End$--invariant spherical round metric on $\cp$ for which the fixed points $p_\End,q_\End$ of $R_\End$ are antipodal points at distance $1$. Let us denote by $g_\End=\dev^*(g_0)$ the Riemannian metric and by $d_\End$ the  distance function induced on $\Sigmaw$. By construction, the M\"obius completion is the metric completion of $(\Sigmaw,d_\End)$.

\begin{lemma}\label{lem:quotient_metric}
    Let $U\subseteq \Compl{\sigma}{\Sigmaw}$ be a $\delta_\End$--invariant  neighborhood of $\End$. Then the distance between $\delta_\End$--orbits defines a metric on $U/\langle \delta_\End \rangle$ with respect to which the quotient map
    $$
    \pi_\End:U \setminus \{\End\}\to (U\setminus \{\End\})/\langle\delta_\End \rangle
    $$
    is a locally isometric covering map.
\end{lemma}
\proof
    Let $\pi_\End(u),\pi_\End(v)\in U/ \langle \delta_\End \rangle$. Then their distance is defined to be
    $$
    d(\pi_\End(u),\pi_\End(v)):= \inf \{ d(\delta^n_\End(u),\delta^m_\End(v)) \ | \ m,n \in \ZZ \}.
    $$
    Since  the action on $U\setminus \{\End\}$ is isometric, free and proper (cf. Corollary~\ref{cor:proper_free_action}),  by  \cite[Proposition I.8.5]{BH99} we get our statement in the complement of the end. To include the end it is enough to notice that it is an isolated fix point and that no orbit accumulates to it, since the holonomy is elliptic. 
\endproof
	
\begin{lemma}\label{lem:cocompact_action}
    Let $U\subseteq \Compl{\sigma}{\Sigmaw}$  be a $\delta_\End$--invariant neighborhood of $\End$ on which  $\delta_\End$ acts cocompactly. Then the following holds.
    \begin{enumerate}
        \item $U$ is complete.
        \item If $V\subseteq U$ is closed and $\delta_\End$--invariant, then $V$ is complete and $\delta_\End$ acts on $V$ cocompactly.
    \end{enumerate}
\end{lemma}
\proof
$\ $
\begin{enumerate}
    \item   Let $x_n\in U $ be a Cauchy sequence.  Let us denote by $F_n$ a (coarse) compact fundamental domain for the action $\langle \delta_\End \rangle \curvearrowright U$ containing $x_n$. If the sequence of $F_n$ eventually stabilizes to some $F$, then eventually the sequence $x_n$ lies entirely in $F$, hence converges in it by compactness. So let us assume that the sequence $F_n$ does not stabilize. We claim that since $x_n$ is a Cauchy sequence this forces $d_\End (x_n,\End)$ to decrease to zero, i.e. $x_n$ converges to $\End$. Indeed, since the holonomy is elliptic and the metric invariant, if $|n-m| $ is large enough then the shortest curve between a point in $F_n$ and a point in $F_m$ goes through $\End$.
	
	\item If $V$ is closed then it is complete by completeness of $U$. Let $F$ be a (coarse) compact fundamental domain for the action $\langle \delta_\End \rangle \curvearrowright U$. Since $V$ is invariant we get $V/\langle \delta_\End \rangle = (V\cap F)/\langle \delta_\End \rangle$, and this is compact because $V\cap F$ is. 
\end{enumerate}
\endproof

We have seen in Proposition~\ref{prop:open_horocycles} that, when the holonomy is elliptic, horocycles contain metric balls (cf. Figure~\ref{fig:canonica_max_nbhd}). We now describe a sufficient condition on a metric ball to be fully contained in a horocycle. Notice that the following statement fails in the case of parabolic holonomy (see Remark~\ref{rmk:parabolic_not_open}).

\begin{lemma}\label{lem:max_radius}
    For each $\End\in \Ends(\Sigmaw)$ let $\critrad{\End} := d_\End(\End,\Binf{\sigma}{\Sigmaw} \setminus \{ \End\})$. Then $\critrad{\End} > 0$ and for all $0 < r < \critrad{\End}$, there is a proper horocyclic neighborhood of $\End$ containing $\Ballc{\sigma}(\End,r)$. 
    %Moreover, $\devc{\sigma}(\Ballc{\sigma}(\End,\critrad{\End})) \subsetneq \Ball(p_\End,1)$.
\end{lemma}

Henceforth we call $\critrad{\End} := d_\End(\End,\Binf{\sigma}{\Sigmaw} \setminus \{ \End\})$ the \emph{critical radius} of $\End$.

\proof
    For the first part of the lemma, let $V$ be a proper (i.e. $\Clot(V) \subsetneq \Sigma \cup \{\End\}$) horocycle based at $\End$. By Proposition~\ref{prop:open_horocycles}  $V$ is open, so there is $r>0$ such that $\Ballc{\sigma}(\End,r) \subseteq V$. We claim that $\Ballc{\sigma}(\End,r) \subset \Sigma \cup \{\End\}$, from which it follows that $\critrad{\End} \geq r > 0$. Recall that $\delta_\End$ acts cocompactly on $\Clot(V)$, therefore $\Clot(V)$ is complete by Lemma~\ref{lem:cocompact_action}. It follows that
    $$
    \Cloc{\sigma}(\Ballc{\sigma}(\End,r)) \subseteq \Cloc{\sigma}(V) = \Clot(V) \subsetneq \Sigma \cup \{\End\}.
    $$
    
    Next, let $r<\critrad{\End}$. Suppose by contradiction that, for every proper horocyclic neighborhood $N$ of $\End$, there was a point $x \in \Ballc{\sigma}(\End,r) \setminus N$. Fix $\{N_k\}$ a sequence of proper horocyclic neighborhoods of $\End$ such that $N_k \subset N_{k+1}$ and $\cup N_k = \Sigmaw \cup \{\End\}$. Let $x_k \in \Ballc{\sigma}(\End,r) \setminus N_k$. For every $k$, let $r_k := d_\End(\End, \partial N_k)$. As $\delta_\End$ acts cocompactly and by isometries on $\partial N_k$, there is some point on $\partial N_k$ at distance $r_k$ from $E$. The fact that $\End \notin \partial N_k$ and the sequence $\{N_k\}$ is nested further implies that:
    $$
    r_k >0, \qquad r_k \leq r_{k+1} \qquad \text{and} \qquad \lim_{k \to \infty} r_k = \critrad{\End}.
    $$
    Notice that the second inequality is due to the fact that $\partial N_k$ separates $\End$ from $\partial N_{k+1}$. Similarly, $\partial N_k$ separates $\End$ from $x_k$, therefore $r_k \leq d_\End(\End,x_k) < r$, hence in the limit we get $\critrad{\End}=\lim_{k \to \infty} r_k \leq r$, in contradiction with the choice of $r$.

\endproof

\begin{corollary}\label{cor:complete_balls}
	For each $\End\in \Ends(\Sigmaw)$ and $0 < r < \critrad{\End}$, $\Ballt(\End,r)= \Ballc{\sigma} (\End,r)$ and $\Clot (\Ballt(\End,r))$ is complete. Moreover, $\Clot (\Ballt(\End,r)) =  \Cloc{\sigma} (\Ballc{\sigma}(\End,r))$.
\end{corollary}
\proof
    By Lemma~\ref{lem:max_radius}, the ball $\Ballt(\End,r)$ is contained in a proper horocyclic neighborhood $V$ of $\End$. It follows that $\Ballc{\sigma}(\End,r) \subset \Sigmaw \cup \{\End\}$ and so $\Ballt(\End,r) = \Ballc{\sigma} (\End,r)$.
    
    By Lemma~\ref{lem:cocompact_action}, the closed ball $\Clot (\Ballt(\End,r))$ is complete.\par
    
    Finally, since $\Clot (\Ballt(\End,r))$ contains $\Ballt(\End,r)$ and is complete, it must contain the completion of $\Ballt(\End,r)$. 
	Since $\Compl{\sigma}{\Sigmaw}$ is complete we have that  $ \Cloc{\sigma} (\Ballc{\sigma}(\End,r)) $ coincides with the completion of $\Ballc{\sigma}(\End,r)$. But we also know that $\Ballc{\sigma}(\End,r)=\Ballt(\End,r)$. So $ \Cloc{\sigma} (\Ballc{\sigma}(\End,r)) $ coincides with the completion of $\Ballt(\End,r)$, and it is therefore contained in $\Clot (\Ballt(\End,r))$.
\endproof  

Recall that a metric space $Z$ is \emph{star-shaped at a point} $x \in Z$ if, for every $y \in Z$ there is a geodesic in $Z$ connecting $x$ to $y$.

\begin{lemma}\label{lem:star_shaped}
    For each $\End\in \Ends(\Sigmaw)$ and $0 < r < \critrad{\End}$, the open ball $\Ballc{\sigma}(\End,r)$ is star-shaped at $\End$.
\end{lemma}

\proof
    Let $x \in \Ballc{\sigma}(\End,r)$ and let $r' := d_\End(x,\End) < r$. By Lemma~\ref{lem:ideal_points_are_accessible}, for all $c>0 $ we can pick a  continuous curve $\eta_c:[0,1)\to \Sigmaw$ such that $\eta_c(0)=x$, $\lim_{t\to 1} \eta_c(t)=\End$ and $r'\leq \ell(\eta_c)\leq r'+c$. For each $t\in[0,1)$ we have
	$$
	d(\eta_c(t),\End)\leq \ell(\eta_c([t,1))) \leq \ell(\eta_c([0,1))) \leq r'+c.
	$$
	In particular, for $c<r-r'$ we get that $d_\End(\eta_c(t),\End)<r$, i.e. $\eta_c$ is entirely contained in $\Ballc{\sigma}(\End,r)\cap \Sigmaw$. Let $\gamma_n:[0,1)\to  \Ballc{\sigma}(\End,r)$ be the curve obtained for $c=\frac{1}{n}$.
	
	Consider the quotient $\pi_\End:\Cloc{\sigma}(\Ballc{\sigma}(\End,r))\to \Cloc{\sigma}(\Ballc{\sigma}(\End,r))/\langle \delta_\End \rangle=:Y$. It follows from Lemma~\ref{lem:quotient_metric} that $\pi_\End$ is a branched covering map onto a  metric space, branching only at $\End$; let us denote by $d_Y$ the distance in $Y$. Moreover by Lemma~\ref{lem:max_radius} the ball $\Ballc{\sigma}(\End,r)$ is properly contained in a horocycle. Since $\delta_\End$ acts cocompactly on horocycles, it follows that $Y$ is compact by Lemma~\ref{lem:cocompact_action}.
	Notice that since $\delta_\End$ acts by isometries and $\End$ is the only fixed point, we also have that $r'=d_\End(x,\End)=d_Y(\pi_\End(x),\pi_\End(\End))$.\par
	
	Projecting the curves $\gamma_n$ to the quotient we obtain curves $\pi_\End \circ \gamma_n:[0,1)\to Y$ such that $\pi_\End \circ \gamma_n(0)=\pi_\End(x)$, $\lim_{t\to 1} \pi_\End \circ \gamma_n(t)=\pi_\End(\End)$ and
	$$
	d_Y(\pi_\End(\End),\pi_\End(x))=r'\leq \ell(\pi_\End \circ \gamma_n)\leq r'+\frac{1}{n}.
	$$
	In particular, by Arzel\`a-Ascoli we can extract a uniform limit $\overline{\gamma}:[0,1]\to Y$. By the above length inequality we obtain 
    $$d_Y(\pi_\End(\End),\pi_\End(x))=r'= \ell(\overline{\gamma}) = \lim_{n\to\infty} \ell(\pi_\End \circ \gamma_n)  $$

	i.e. $\overline{\gamma}$ is a geodesic from $\pi_\End(x)$ to $\pi_\End(\End)$. Notice that it goes through $\pi_\End(\End)$ only at one endpoint; so we can lift it to a curve $\gamma:[0,1)\to \Cloc{\sigma}(\Ballc{\sigma}(\End,r))$ starting at $x$ and limiting to $\End$, of the same length $r'$. By the same argument as the beginning, $\gamma$ is completely contained in the open ball $\Ballc{\sigma}(\End,r)$, so this is the desired geodesic.
\endproof

We now consider the restriction of $\devc{\sigma}$ to a ball around an end $\End \in \Ends(\Sigmaw)$, i.e.
$$
\devc{\sigma} : \Ballc{\sigma}(\End,r) \to \Ball(p_\End,r),
$$
and we find the values of $r$ for which it is a covering map, branching only at $\End$. 
The proof is reminiscent of (and based on) the classical fact that a local isometry from a complete Riemannian manifold to a connected one is a covering map. Notice that in our setting $\devc{\sigma}$ is not locally isometric at $\End$ (not even locally injective), and on the other hand $\Cloc{\sigma}( \Ballc{\sigma}(\End,r)) \setminus \{\End\}$ is not complete. The proof shows how to deal with this, and also provides quantitative control on the critical radius.

\begin{proposition}\label{prop:branched_covering}
    For each $\End\in \Ends(\Sigmaw)$ we have that $\critrad{\End} \leq 1$. Moreover:
    \begin{enumerate}
        \item \label{item:branched_covering_1} for each $0 < r \leq \critrad{\End}$, $\devc{\sigma}$ maps $\partial \Ballc{\sigma}(\End,r)$ to $\partial \Ball(p_\End,r)$;
        \item \label{item:branched_covering_2} for each $0 < r \leq \critrad{\End}$, $\devc{\sigma} : \Ballc{\sigma}(\End,r) \to \Ball(p_\End,r)$ is a branched covering map, branching only at $\End$.
    \end{enumerate}
\end{proposition}

\proof
We are first going to prove statements~\eqref{item:branched_covering_1} and~\eqref{item:branched_covering_2} for $r \leq \min\{\critrad{\End}, 1\}$, and then we will show that $\critrad{\End} \leq 1$.

We begin with the following observation. Suppose $r < \critrad{\End}$ and let $x \in \partial \Ballc{\sigma}(\End,r)$. Then $d_\End(\End,x)=r$ and $d_0(p_\End,\devc{\sigma}(x))\leq r$. Let $r'>0$ such that $r < r' < \critrad{\End}$. Then $x \in \Ballc{\sigma}(\End,r')$ and by Lemma~\ref{lem:star_shaped} there exists a geodesic $\gamma_r$ from $x$ to $\End$ contained in $\Ballc{\sigma}(\End,r')$. Observe that $r =\ell(\gamma_r) = \ell(\devc{\sigma}(\gamma_r))$. Notice $\dev$ is a local isometry on $\Sigmaw$, so $\gamma_r$ maps to a geodesic in $\cp$. \par

Next, additionally assume that $r < 1$, the diameter of $\cp$. Then the curve $\gamma_{r}$ maps to a simple geodesic arc, starting from $p_\End$ and avoiding $q_\End$, of length $r<1$. Since the choice of $x$ above was arbitrary, it follows that $\devc{\sigma}(\partial \Ballc{\sigma}(\End,r)) \subseteq \partial \Ball(p_\End,r)$. In particular, it avoids $q_\End$. This concludes the proof of \eqref{item:branched_covering_1} in the case where $r < \min\{\critrad{\End}, 1\}$. The limiting case $r = \min\{\critrad{\End}, 1\}$ follows by continuity of the developing map.

We now start the proof of (\ref{item:branched_covering_2}). To begin with, we claim that when $r< \min\{\critrad{\End}, 1\}$, each component of  $\partial \Ballc{\sigma}(\End,r)$ is isometric to a complete line.
Since $r < 1$, $\partial \Ball(p_\End,r)$ is a circle in $\cp$. Since $r < \critrad{\End}$, we have that $\partial \Ballc{\sigma}(\End,r) \subset \Sigmaw$ and $\devc{\sigma}$ is a local homeomorphism on it. In particular, $\partial \Ballc{\sigma}(\End,r)$ is a $1$--dimensional submanifold of $\Sigmaw$; moreover it is closed in $\Cloc{\sigma}(\Ballc{\sigma}(\End,r))$, hence complete by Corollary~\ref{cor:complete_balls}. Then $\devc{\sigma}$ induces a local isometry from the complete manifold $\partial \Ballc{\sigma}(\End,r)$ to the connected manifold $\partial \Ball(p_\End,r)$; it follows that it is a Riemannian covering map. Notice that $\langle \delta_\End \rangle$ is an infinite cyclic group acting on $\partial \Ballc{\sigma}(\End,r)$ properly and freely by Corollary~\ref{cor:proper_free_action}, hence each component of  $\partial \Ballc{\sigma}(\End,r)$ must be isometric to a complete line. 

Now we claim that, for all $0 < r \leq \min\{\critrad{\End}, 1\}$,  $\devc{\sigma} : \Ballc{\sigma}(\End,r) \to \Ball(p_\End,r)$ is a branched covering map, branching only at $\End$. First notice that 
\begin{equation}\label{eq:no_preimages_of_end}
    \Ballc{\sigma}(\End,r) \setminus \{\End\} =  \Ballc{\sigma}(\End,r) \setminus \{\devc{\sigma}^{-1}(p_\End)\}.
\end{equation}
Indeed suppose $z\in \Ballc{\sigma}(\End,r)$ is another point developing to $p_\End$, then there is $r' < r \leq \critrad{\End}$ such that $z\in \Ballc{\sigma}(\End,r')$, and a geodesic $\gamma$ from $z$ to $\End$ contained in $\Ballc{\sigma}(\End,r')$. Since $\devc{\sigma}(\End) = \devc{\sigma}(z) = p_\End$, this geodesic $\gamma$ has to cover at least a great circle through $p_\End$ in $\cp$, hence $d_\End(\End,z)\geq 2$. But $r' < r \leq1$ forbids this. In particular, we get a well defined local homeomorphism
	$$
	\varphi := \restr{\devc{\sigma}}{\Ballc{\sigma}(\End,r) \setminus \{\End\}}: \Ballc{\sigma}(\End,r) \setminus \{\End\} \to \Ball(p_\End,r) \setminus \{p_\End\}.
	$$
It is enough to show that this is a covering map. We are going to show that every point in $\Ball(p_\End,r) \setminus \{p_\End\}$ is evenly covered. Let $y \in \Ball(p_\End,r) \setminus \{p_\End\}$ and let $r_y := d_0(p_\End,y)$. Notice that $0 < r_y < r$. Since $\devc{\sigma}$ is a covering map between $\partial \Ballc{\sigma}(\End,r_y)$ and $\partial \Ball(p_\End,r_y)$, there is $\epsilon_y >0$, such that $\Ball(y,\epsilon_y) \cap \partial \Ball(p_\End,r_y)$ is evenly covered. Let
$$
\delta_y := \min\{\epsilon_y, r- r_y, r_y\}.
$$
Notice that the ball $\Ball(y,\delta_y)$ is entirely contained in $\Ball(p_\End,r)\setminus \{p_\End\}$. Then we claim that $\Ball(y,\delta_y)$ is evenly covered. Let $z \in \devc{\sigma}^{-1}(y) \cap \Ballc{\sigma}(\End,r)$. By definition of $\delta_y$, $\Ballc{\sigma}(z, \delta_y)$ is entirely contained in $\Ballc{\sigma}(\End,r) \setminus \{\End\}$. In particular, it is smaller than the maximal ball centered at $z$, so it is isometrically mapped to $\Ball(y,\delta_y)$ by $\devc{\sigma}$ (cf. Lemma~\ref{lem:isometric_balls}). This implies that if $z' \in \devc{\sigma}^{-1}(y)\cap \Ballc{\sigma}(\End,r)$ is different from $z$, then $\Ballc{\sigma}(z, \delta_y) \cap \Ballc{\sigma}(z', \delta_y) = \varnothing$. This concludes the proof of \eqref{item:branched_covering_2} in the case $r\leq \min\{\critrad{\End}, 1\}$.

Now suppose by contradiction that $\critrad{\End} > 1$. Then there is $r$ such that $1 < r < \critrad{\End}$, and the open ball $\Ballc{\sigma}(\End,r)$ is star-shaped at $\End$ (cf. Lemma~\ref{lem:star_shaped}). Moreover, the developing map maps $\partial \Ballc{\sigma}(\End,1)$ to $q_\End$, the only point at distance $1$ from $p_\End$ in $\cp$.  Since the open ball $\Ballc{\sigma}(\End,r)$ is entirely contained in $\Sigmaw$, and contains $\partial \Ballc{\sigma}(\End,1)$, the developing map is a local homeomorphism on $\partial \Ballc{\sigma}(\End,1)$. In particular $\partial \Ballc{\sigma}(\End,1)$ is discrete. On the other hand, for every $r' < 1 = \min\{\critrad{\End}, 1\}$ we can apply the first part of the proof where we proved that $\devc{\sigma}$ maps $\partial \Ballc{\sigma}(\End,r')$ to $\partial \Ball(p_\End,r')$, and
$$
\lim_{r' \to 1^- } \partial \Ball(p_\End,r') = \{q_\End\}.
$$
This implies that, for radii $r'<1=\min\{ \critrad{\End}, 1\}$ sufficiently close to $1$, $\partial \Ballc{\sigma}(\End,r')$ is a disjoint union of circles, contradicting that each connected component is isometric to a complete line.
\endproof 

\begin{maintheoremc}{D} \label{thm:maxneighborhood}
Let $\sigma \in \Modulit(\Sigma)$, and let $\End$ be an end. Then there is a neighborhood $\maxN_{\End}$ of $\End$ in $\Compl{\sigma}{\Sigmaw}$ onto which the developing map for $\sigma$ restricts to a branched covering map, branching only at $\End$, and with image a round disk in $\cp$.
\end{maintheoremc}

\proof
We can just take  $\maxN_{\End}$ to be any ball $\Ballc{\sigma}(\End,r)$ satisfying the conditions of Proposition~\ref{prop:branched_covering}.
\endproof

Let $\End\in \Ends ( \Sigmaw)$, and let $\critrad{\End}$ be its critical radius. The open metric ball $\maxN_{\End}=\Ballc{\sigma}(\End,\critrad{\End})$ plays the role of a \textit{canonical maximal neighborhood} of $\End$, similar to the maximal round balls in \cite{KP94}.
Indeed, it develops to a round ball in $\cp$, and by definition of $\critrad{\End}$, the boundary of $\maxN_{\End}$ contains an ideal point. 
However, note that we have normalized things ``locally'' at $\End$, by fixing the $R_\End$--invariant round metric on $\cp$ for which the fixed points $p_\End,q_\End$ of the holonomy at $\End$ are antipodal points of distance 1 (here $R_\End=\hol(\delta_\End)$ denotes the peripheral holonomy at $\End$).
Then $\maxN_{\End}$ is defined as a metric ball for the induced metric $g_\End$ on $\Compl{\sigma}{\Sigmaw}$.
If $\End '$ is a different end, then the metric ball around $\End '$ (with respect to $g_\End$) does not necessarily agree with $\maxN_{\End '}$, which would be defined as a metric ball for the metric $g_{\End '}$.

Moreover one can observe that when $r<\critrad{\End}$ the ball  $\partial \Ballc{\sigma}(\End,r)$  contains a horocycle and  is contained in a horocycle. Therefore each component of its boundary is contained in the lune between two horocycles.
Since $\Ballc{\sigma}(\End,r)$ is star-shaped at the end, and $\partial \Ballc{\sigma}(\End,r)$ is invariant under the action of the peripheral $\delta_\End$, we can see that $\partial \Ballc{\sigma}(\End,r)$ is actually connected and isometric to a complete line. In particular, $\partial \Ballc{\sigma}(\End,r)$ is the universal cover of $\partial \Ball (p_\End,r)$, and $\Ballc{\sigma}(\End,r)\setminus \{\End\}$ is isometric to the universal cover of $\Ball (p_\End,r)\setminus \{p_\End\}$.

\begin{remark}
If $\sigma$ is the tame and relatively parabolic structure induced by complete hyperbolic metric of finite area, then $\dev$ is a global diffeomorphism, and horocycles develop to round disks. In particular, Theorem~\hyperref[thm:maxneighborhood]{D} holds for such a structure.
However, there is no analogue of Theorem~\hyperref[thm:maxneighborhood]{D} in the general parabolic case. For example, consider the structure obtained by grafting $\sigma$ along an ideal arc, and let $\End$ be an end covering one of the endpoints of the grafting arc. If $U$ is any $\delta_\End$--invariant neighborhood of $\End$, then $V=\devc{\sigma}(U)$ is invariant under a parabolic transformation and contains its fixed point $p_\End=\devt(\End)$ in its interior. This forces $V=\cp$.
In particular, we see that the local homeomorphism (analogous to the one considered in the proof of Proposition~\ref{prop:branched_covering})
$$
	\varphi := \restr{\devc{\sigma}}{U\setminus \devc{\sigma}^{-1}(p_\End)}: U\setminus \devc{\sigma}^{-1}(p_\End) \to V\setminus \{p_\End\}=\CC
	$$
	cannot be a covering map, because it is not injective and the image is simply connected.
\end{remark}

Throughout this section we have worked under the normalization in which the fixed points $p_\End,q_\End \in \cp$ of the rotation $R_\End$ are antipodal points at distance $1$.
As established in Proposition~\ref{prop:branched_covering}, it follows that the critical radius of an end $E$ satisfies $\critrad{\End}\leq 1$.
We conclude this chapter by discussing what happens when a tame structure $\sigma$ has an end $\End$  with elliptic holonomy and $\critrad{\End}=1$.

 \begin{remark}[Structures on a twice-punctured sphere]\label{rem:twice}
 Suppose $\sigma$ is tame and has an end $\End$  with elliptic holonomy and $\critrad{\End}=1$. 
 By \eqref{item:branched_covering_1} in Proposition~\ref{prop:branched_covering} all the points on the boundary of $\maxN_{\End}$ must develop to $q_\End$. It follows from the  proof of Proposition~\ref{prop:branched_covering} that in this case the boundary of $\maxN_{\End}$ cannot contain any isolated points in $\Sigmaw$. As a result, $\maxN_{\End}=\Sigmaw$.
 By tameness, this forces all the ends different from $\End$ to develop to $q_\End$. We claim that in this case $\Sigma$ must be a twice-punctured sphere, and $\sigma$ is the structure associated to a power map $z\mapsto z^\alpha$ for some $\alpha \in \RR\setminus \ZZ$.
 To see this, assume by contradiction that there is a peripheral element $\gamma \in \pi_1(\Sigma)$  distinct from any power of the peripheral element $\delta_\End$ which fixes $\End$. Then $\gamma$ moves $\End$ to another end $\gamma \End \neq \End$.
 By equivariance and tameness of the developing map (see Lemma~\ref{lem:devt_is_equivariant}) we see that
 $$ q_\End = \devt (\gamma \End) = \hol (\gamma) \devt(\End) = \hol (\gamma) p_\End $$
 On the other hand, $\gamma$ fixes an end $\End ' \neq \End$. 
 It follows that $\devt(\End ')=q_\End=\hol(\gamma)q_\End$.
 We get $\hol (\gamma) p_\End = \hol(\gamma)q_\End$, which is absurd. Therefore all peripheral elements are powers of a fixed one. 
 But the only orientable surface in which this happens is a sphere with two punctures.
 Notice that this surface has zero Euler characteristic, $\Sigmaw$ identifies with $\mathbb C$, and we can normalize things so that $\dev(z)=e^{az}$, $p_\End=0$ and $q_\End=\infty$, for some $a\in \CC^*$. 
 Deck transformations are generated by $z\mapsto z+2\pi i$, and the holonomy by $w\mapsto e^{2\pi i a}w$; ellipticity of the holonomy means  $a \in \RR\setminus \ZZ$. The M\"obius completion is obtained by adding just two ideal points, for $\Re(z)\to \pm \infty$, mapping to $q_\End=\infty$ and $p_\End=0$ respectively. Structures of this type can be defined by a spherical metric with two cone points and coaxial holonomy.
 \end{remark}
 
 \begin{remark}
 Spherical metrics with cone points and coaxial holonomy exist also on a surface of negative Euler characteristic (see \cite{ER04,MP16}), and provide  examples  of  structures  with degenerate holonomy. However such a structure must have some apparent singularities (i.e. punctures with trivial holonomy, see \cite{G19}), whose presence forces the critical radius to be strictly less than 1 at every end with elliptic holonomy.  Indeed, if $\End$ is an end with elliptic holonomy, then the family of neighborhoods $\Ballc{\sigma}(\End,r)$ must hit another end (possibly one covering an apparent singularity) before $r=1$. As an illustrative example, consider the structure obtained by puncturing an additional point on a sphere endowed with a spherical metric with two cone points.
\end{remark}

	%%%%%%%%%%%%%%%%%%%%%%%%%%%%%%%%%
	%%%%%%%%%%%%%%%%%%%%%%%%%%%%%%%%%
	%%%%%%%%%%%%%%%%%%%%%%%%%%%%%%%%%

	\subsection{The index of a puncture}\label{subsec:index_puncture} 
	
	Using the neighborhoods constructed in the previous section (namely Theorem~\hyperref[thm:maxneighborhood]{D}), we can define a numerical invariant of the complex projective structure for each puncture, which we call the \emph{index}. This is essentially the angle that the developed image of a peripheral curve makes around the image of a corresponding end. 
	
	Let $\sigma \in \Modulit(\Sigma)$ be a structure represented by a pair $(\dev,\hol)$.
	Let $x$ be a puncture of $\Sigma$, and let $\eta$ be a positive peripheral curve in $\Sigma$ around $x$. 
	This can be chosen so that for any end $\End$ covering $x$ (in the sense of Remark~\ref{rem:end_cover_puncture}), the lift of $\eta$ which is asymptotic to $\End$ is entirely contained in a neighborhood $V_\End$ of $\End$  on which $\dev$ is a branched covering map, branching only at $\End$ (see Theorem~\hyperref[thm:maxneighborhood]{D}).

    Let us fix an end $\End$, and let $\delta_\End \in \pi_1(\Sigma)$ be the positive peripheral deck transformation fixing $\End$. 
	We recall that $p_\End:=\devt(\End)$ is one of the two fixed points for the elliptic transformation $\hol(\delta_\End)$ (see Lemma~\ref{lem:devt_is_equivariant}). Let us normalize so that $\hol(\delta_\End)$ fixes $0$ and $\infty$.
	Let $\widetilde \eta \subset V_\End$ be the lift of $\eta$ in $V_\End$, and choose $\widetilde \eta_0 \subset \widetilde \eta$ to be a fundamental domain for the action $\langle \delta_\End \rangle \curvearrowright \widetilde \eta$. 	Let $\zeta:=\dev ( \widetilde \eta_0) \subset \dev (V_\End) \setminus \{0\}$.  Notice that freely homotoping $\eta$ deeper into the puncture results in a homotopy of $\zeta$ in the complement of $p_\End = 0$, because there are no other preimages of $p_\End$ in $V_\End$  (cf. \eqref{eq:no_preimages_of_end} in Proposition~\ref{prop:branched_covering}).
	
% 	Let $\omega$ be the unique  meromorphic $1$--form on $\cp$ with simple poles at $p_\End$ and $q_\End$, and residue $\Res_{p_\End}(\omega)=1$. 

	The \emph{index of the structure $\sigma$ at the puncture $x$} is defined to be the number
	$$
% 	\Ind_\sigma(x) :=\left\| \int_{\eta_0} \omega \right\|.
% 	\Ind_\sigma(x) :=\Im \left( \int_{\eta_0} \omega \right).
	\Ind_\sigma(x) := \Im \left( \int_{\zeta} \frac{dz}{z} \right). 
	$$
	When clear from the context, we will usually drop the $\sigma$ and write $\Ind(x) = \Ind_\sigma(x)$.
	
	We remark explicitly that this definition does not depend on any of the choices involved. Indeed, let us choose a parametrization $\zeta:[0,1]\to \CC\setminus \{0\}, \ \zeta(s)=r(s)e^{i\theta(s)}$, where $\theta:[0,1]\to \RR$ is a determination of the argument function on $\CC \setminus \{0\}$, and $r : [0,1]\to \RR$. 
	A direct computation in local coordinates shows that
	$$
	\int_{\zeta} \frac{dz}{z}  = \log \left( \frac{r(1)}{r(0)}\right)+i(\theta(1)-\theta(0)).
	$$ 
	Notice that since $\eta$ is chosen to be a peripheral curve, its holonomy is elliptic. Therefore we get
	$r(1)e^{i\theta(1)} = e^{i\varphi}  r(0)e^{i\theta(0)}$,
    where $\varphi$ is such that $\hol(\delta_\End)z=e^{i\varphi}z$. It follows that 
	$$
	\Ind_\sigma(x) = 2\pi k+ \varphi,
	$$ 
	where $k\in \ZZ$  counts the number of times $\zeta$ turns around $0$ anticlockwise. Notice that the index is always positive, since $\delta_\End$ was chosen to be a positive peripheral.

    \begin{remark}\label{rem:grafting_indices}
    Let $\sigma \in \Modulit(\Sigma)$, and let $x$ and $y$  be   punctures.
    If $\eta$ is a graftable arc joining $x$ to $y$, then the following holds:
    \begin{itemize}
        \item if $x\neq y$ then $\Ind_{ \Gr (\sigma,\eta)}(x)=\Ind_\sigma(x)+2\pi$
    and $\Ind_{ \Gr (\sigma,\eta)}(y)=\Ind_\sigma(y)+2\pi$,
    
    \item if $x= y$ then $\Ind_{ \Gr (\sigma,\eta)}(x)=\Ind_\sigma(x)+4\pi$.
    \end{itemize}
    \end{remark}

	%%%%%%%%%%%%%%%%%%%%%%%%%%%%%%%%%
	%%%%%%%%%%%%%%%%%%%%%%%%%%%%%%%%%
	%%%%%%%%%%%%%%%%%%%%%%%%%%%%%%%%%

	\section{The complex analytic point of view}\label{sec:complex}

	The theory of $\cp$--structures enjoys fundamental interactions with the study of second--order linear ODEs on complex domains, namely through the use of the Schwarzian derivative.
    The purpose of this chapter is to describe the complex analytic counterpart to the structures in $\Modulit(\Sigma)$ (see Theorem~\hyperref[thm:tamerelell=merodoublenonint]{E}). These are described by meromorphic quadratic differentials satisfying certain conditions on their Laurent expansion around poles.

	\subsection{Local theory at regular singularities}\label{subsec:local_theory}
    We start by reviewing the classical theory for the convenience of the reader, with a particular focus to the behavior around singularities of the coefficients  (see \cite{I44,H69}). This will provide the local model for our structures around the punctures.
	
	Let us consider a holomorphic function $q:\pDisk\to \CC$ on the punctured unit disk $\pDisk=\{z\in \CC \ | \ 0<|z|<1\}$ with a double pole at the origin with \textit{leading coefficient} $a$, i.e. a function of the form $q(z)=\frac{a}{z^2}+O\left(\frac{1}{z}\right)$. We will consider the second--order linear ODE
	\begin{equation}\label{eq:IIODE}
	u''+\frac{1}{2}qu=0, \quad \textrm{ for } \quad u:\pDisk\to \CC,
	\end{equation}
	as well as the \textit{Schwarz equation}
	\begin{equation}\label{eq:schwarz}
	\sch f=q, \quad \textrm{ for } \quad f:\pDisk\to \cp,
	\end{equation}
	where the operator
	$$ \sch f = \left(\frac{f''}{f'} \right)'-\frac{1}{2}\left( \frac{f''}{f'}\right)^2
	$$ 
	is the \textit{Schwarzian derivative}. The main properties of $\sch$ are the following:
	\begin{enumerate}
	    \item (\textit{invariance})  $\sch f=0$ if and only if $f$ is the restriction of some M\"obius transformation.
	    \item (\textit{cocycle})  if $f,g$ are locally injective  holomorphic functions for which the composition is defined, then 
	    $\sch (f\circ g)=g^*(\sch f) + \sch g$.
	\end{enumerate}
	
	The relationship between the two equations above is well known (see \cite[Appendix D]{H69}), and can be summarized as follows: if $u_1,u_2$ are linearly independent solutions for (\ref{eq:IIODE}), then $f=\frac{u_1}{u_2}$ is a solution for (\ref{eq:schwarz}); conversely, any solution for (\ref{eq:schwarz}) is obtained in this way. In both cases, since the domain of the equation is not simply-connected, these equations can have non trivial monodromy, i.e. solutions are to be considered as multi-valued functions, or as single-valued functions on a suitable covering domain.\par
	The classical theory of linear ODEs (see \cite[\S 15.3]{I44}, or \cite[\S 5]{AB18} for a more recent treatment) provides an explicit description of the local solutions of (\ref{eq:IIODE}). Firstly, the \textit{indicial equation} of (\ref{eq:IIODE}) is given by
	$$r(r-1)+\frac{a}{2}=0.$$
	Let $r_1,r_2\in \CC$ be its solutions; then one has two cases:
	\begin{enumerate}
		\item if $r_1-r_2\not \in \ZZ$ then (\ref{eq:IIODE}) has two linearly independent solutions of the form $u_k(z)=z^{r_k}h_k(z)$ for $k=1,2$, where $h_k$ is holomorphic on $\DD$, $h_k(0)\neq 0$;
		\item if $r_1-r_2 \in \ZZ$ then (\ref{eq:IIODE}) has two linearly independent solutions of the form $u_1(z)=z^{r_1}h_1(z)$  and $u_2(z)=z^{r_2}h_2(z)+Cu_1(z)\log(z) $   where $C\in \CC$, and $h_k$ is holomorphic on $\DD$, $h_k(0)\neq 0$ for $k=1,2$.
	\end{enumerate}

	An analogous dichotomy for solutions of (\ref{eq:schwarz}) is easier to state if we write the leading coefficient in the form $a=\frac{1-\theta^2}{2}$, where $ \theta =\pm \sqrt{1-2a}$ will be called the \textit{reduced exponent} of $q$ at $z=0$. 
	With respect to the terminology used in \cite{AB18}, the \textit{exponent} of $q$ at $z=0$ is $r=\pm 2\pi i \sqrt{1-2a}= 2\pi i \theta$. For the reader's convenience, we remark that in \cite{AB18} a slightly different form of the Schwarzian derivative is used, leading to a different normalization for constants in the correspondence between differentials and monodromy of solutions.
    Observing that $\pm \theta = r_1-r_2$, and recalling the relation $f=\frac{u_1}{u_2}$, one has the following:
	\begin{enumerate}
		\item if $\theta \not \in \ZZ$ then (\ref{eq:schwarz}) has a solution of the form $f(z)=z^{\theta}M(z)$, where $M$ is holomorphic at $z=0$, $M(0)\neq 0$;
		\item if $\theta \in \ZZ$ then (\ref{eq:schwarz}) has a solution of the form $f(z)=z^{\theta}M(z)+C\log(z)$, where $C\in \CC$, and $M$ is holomorphic at $z=0$, $M(0)\neq 0$.
	\end{enumerate}

	For each $q$ one can regard a solution to (\ref{eq:schwarz}) as a developing map for a projective structure on $\pDisk$, equivariant with respect to the monodromy group of the equation. Notice that the holonomy of this structure (i.e. the monodromy of (\ref{eq:schwarz})) is a representation $\rho:\pi_1(\pDisk)\to \pslc$ which is just  the projectivization of the monodromy $\widetilde \rho:\pi_1(\pDisk)\to \slc$ of (\ref{eq:IIODE}). If $\gamma$ denotes a simple loop in $\pDisk$ around $z=0$, then the action of the monodromy is given by the linear fractional transformation $\hol(\gamma)\cdot z=e^{2\pi i \theta}z+2\pi i C$.
	
	A direct computation using the above description of solutions to (\ref{eq:schwarz}) leads to the following statement. Here continuous extensions to the origin should be thought in the sense of the end-extension topology introduced in \S\ref{subsec:tame_rel_ell}. 
	
	\begin{lemma}\label{lem:local_theory}
		In the above notation, 
		  the following hold:
		\begin{enumerate}
			\item if $\theta=0$ then $\hol(\gamma)$ is parabolic (necessarily $C\neq 0$);
			
			\item if $\theta \in \ZZ \setminus \{0\}$, then $\hol(\gamma)$ is trivial (if $C=0$) or parabolic (if $C\neq 0$);
			
			\item if $\theta \in \RR \setminus \ZZ$, then $\hol(\gamma)$ is elliptic;
			\item if $\theta \in \ZZ \oplus i\RR$, then $\hol(\gamma)$ is hyperbolic;
			\item if $\theta \in \CC\setminus (\ZZ \oplus i\RR)$, then $\hol(\gamma)$ is purely loxodromic.
		\end{enumerate}
		Moreover, if $\theta \in \RR\setminus \ZZ$, then 	a solution $f$ of (\ref{eq:schwarz}) extends continuously to $z=0$.
	\end{lemma}
	
	As the reader might expect,  projective structures in $\Modulit(\Sigma)$ relate to the elliptic case in the above statement.
	On the other hand, a solution $f$ of (\ref{eq:schwarz}) does not extend continuously to $z=0$ when $\theta=n \in \ZZ$ and $C\neq 0$ (i.e. when $f$ is of the form $f(z)=z^{n}M(z)+C\log(z)$), which can be seen by inspecting the behavior of $f$ along appropriately chosen sequences that spiral into the singularity. 
	A similar phenomenon occurs when $\theta\not\in\RR$.
	
\subsection{Meromorphic projective structures}\label{subsec:meromorphic_structures}
We now recall how to construct projective structures in terms of meromorphic quadratic differentials, and discuss its relationship with our space $\Modulit(\Sigma)$ of tame, relatively elliptic, and non-degenerate structures, introduced in \S \ref{subsec:tame_rel_ell}. This is analogous to the classical parametrization of complex projective structures on closed surfaces by holomorphic quadratic differentials (see \cite[\S 3]{DU} for an expository account). This section includes the proof of Theorem~\hyperref[thm:tamerelell=merodoublenonint]{E}.

Let us fix a complex structure $\overline X$ on the closed surface $\Sigmac$, and let $\overline \sigma_0$ be the $\cp$--structure on $\overline{ X}$ defined by the Poincar\'e uniformization, i.e. the unique conformal metric of constant curvature $-1,0$ or $1$, the exact value depending on the genus $g$ of $\overline{X}$. 
Let $X$ be the induced complex structure on $\Sigma=\Sigmac \setminus \{x_1,\dots,x_n\}$; notice $X$ is a punctured Riemann surface, i.e. each $x_j$ has a neighborhood biholomorphic to $\pDisk$. We consider the space $\Quadddp{X}$ of \textit{meromorphic quadratic differentials} with at worst double poles at the punctures of $X$; these are meromorphic sections of the line bundle $K_X^{ 2}$, where $K_X$ denotes the canonical bundle of $X$. More concretely, by slight abuse of notation, in suitable local complex coordinates around the puncture these differentials can be written as  
$$
q(z)=\left(\frac{a}{z^2}+O\left(\frac{1}{z}\right)\right) dz^2.
$$

The leading coefficient at a double pole is a well-defined invariant of a quadratic differential, i.e. does not depend on the chosen coordinates  (see \cite[\S 4.2]{S84}).
In particular, the local analysis developed in \S \ref{subsec:local_theory} applies, and provides a definition of exponents and reduced exponents of $q$ at a puncture.

Moreover the properties of the Schwarzian derivative ensure that the Schwarz equation $\sch f=q$ is well-defined on $X$, as soon as a background projective structure has been fixed, and we choose the Poincar\'e uniformization $\overline \sigma_0$. Local solutions are in general multi-valued, i.e. they should be considered as functions on the universal cover, equivariant with respect to some representation  $\hol_q:\pi(X)\to \pslc$, which is called the \textit{monodromy} of $q$.
We say a puncture is an \textit{apparent singularity} if $\hol_q(\gamma)$ is trivial for a peripheral loop $\gamma$ around the puncture.
It is a theorem of Luo (see \cite{L93}) that differentials without apparent singularities are locally determined by their monodromy. The analogous results for holomorphic quadratic differentials is due to Hejhal (see \cite{HE75}).

Following \cite[\S 3]{AB18} and \cite[\S 3.1]{GM19} we define a \textit{meromorphic projective structure}  to be the structure $\sigma_q$ induced by a meromorphic quadratic differentials $q\in \Quadddp X$ as follows: a developing map $\dev_q$ for $\sigma_q$ is given by taking a local solution to $\sch f=q$ and considering its analytic continuation as a function  on the universal cover; the monodromy of the differential provides the holonomy $\hol_q$ of the structure. The differential $q$ is recovered from $\sigma_q$ by computing the Schwarzian derivative of $\dev_q$ with respect to the background projective structure $\overline \sigma_0$.
	
For the sake of clarity, we emphasize that this correspondence between meromorphic quadratic differentials and meromorphic projective structures is not canonical, and does depend on the choice of a background projective structure. Changing this choice only translates the differentials by the vector space of holomorphic differentials, hence orders and leading coefficients of poles are well-defined invariant for the projective structure.
	
	We are now ready to provide a proof of the following correspondence. Here $\Modulit(\Sigma)$  is the space of tame, relatively elliptic and non-degenerate structures introduced in \S \ref{subsec:tame_rel_ell}.
	
	\begin{maintheoremc}{E}\label{thm:tamerelell=merodoublenonint}
    Let  $\sigma \in \Moduli(\Sigma)$ and let $X \in \Teich(\Sigma)$ be the underlying complex structure. Then $\sigma \in \Modulit(\Sigma)$ if and only if $X$ is a punctured Riemann surface and $\sigma$ is represented by a meromorphic quadratic differential on $X$ with double poles and reduced exponents in $\RR\setminus \ZZ$.
    \end{maintheoremc}
	\proof 
    We prove the backward direction first.
    Let $X$ be a punctured Riemann surface structure on $\Sigma$, and let $\sigma=\sigma_q$ for some meromorphic quadratic differential  $q\in \Quadddp X$  with reduced exponents  $\theta_i\in \RR\setminus \ZZ$.
    By Lemma~\ref{lem:local_theory}, since the $\theta_i$'s are real but not integers, the developing map for $\sigma$ extends continuously to the punctures (i.e. $\sigma$ is tame), and 
    the peripheral holonomy of $\sigma$ is elliptic at every puncture.
    In particular, the holonomy representation is known to be non-degenerate by \cite[Theorem 6.1]{AB18}, as there are no apparent singularities. Therefore $\sigma \in \Modulit(\Sigma)$.
    
    We now prove the forward direction.
	Let $\sigma \in \Modulit(\Sigma)$, and let $U$ be a neighborhood of a puncture $x$, which is some conformal annulus. We claim that its modulus is infinite.
	Let $\End\in \Ends(\Sigmaw)$ be an end covering $x$, and let $\widetilde U$ be the lift of $U$ around $\End$. 
	By Theorem~\hyperref[thm:maxneighborhood]{D} we can choose $U$ so that $\dev:\widetilde U\to D^*=\dev(\widetilde U)$ is a conformal covering map onto a punctured disk. 
	The family of curves $\Gamma$ in $D^*$ joining the boundary to the puncture has infinite extremal length, lifts to a family of curves in $\widetilde U$ joining  $\partial \widetilde{ U}$  to $\End$, and projects to a family of curves in $U$ joining the boundary to the puncture $x$.
	Since extremal length is conformally invariant, this   family has infinite extremal length in $U$, hence the modulus of $U$ is infinite. This shows that the complex structure $X$ underlying $\sigma$ is that of a punctured Riemann surface.
	
	Finally let us check the conditions on the differential are satisfied. 
	Let $\widetilde q=\sch(\dev)$; recall we have fixed the Poincar\'e uniformization $\overline{\sigma}_0$ as a reference projective structure on $\Sigmac$, and we are taking Schwarzian derivatives with respect to the induced structure on $\Sigma$. Since $\dev$ is a conformal immersion, possibly branching only at the ends, $\widetilde q$ is holomorphic on $\Sigmaw$, possibly with double poles at the ends. By the classical cocycle property  of the Schwarzian, $\widetilde  q$ descends to a meromorphic quadratic differential $q$ with at worst double poles on $\Sigma$. By Lemma~\ref{lem:local_theory}, since the peripheral holonomy is elliptic, the reduced exponents must be in $\RR \setminus \ZZ$.
	\endproof
	
For completeness, with respect to the list of cases in Lemma~\ref{lem:local_theory}, we observe the following.
Differentials with zero reduced exponents at all punctures correspond to parabolic projective structures (see \cite{K69,K7169,K71,DD17,HD19}).
Differentials with integer non-zero reduced exponents and trivial holonomy at the punctures  (apparent singularities) correspond to branched projective structures (see \cite{M72,CDF,CDHL,FR19}).
The next lemma implies that for structures in $\Modulit(\Sigma)$ the absolute value of the exponent at a puncture coincides with the value of the index, as defined in \S \ref{subsec:index_puncture}.

\begin{lemma}\label{lem:index_exponent}
If $q\in \Quadddp X$ has reduced exponent $\pm \theta\in \RR\setminus \ZZ$ at a puncture $x$, then the index of $\sigma_q$ at that puncture is $\Ind_\sigma(x)=2\pi|\theta|$.
\end{lemma}
\proof
Let $z$ be a coordinate around the puncture, let $\eta$ be a simple closed positively oriented peripheral loop around the puncture.
Up to normalizing by a M\"obius transformation, we can assume that a local determination of the developing map is given by $w=f(z)=\dev_{q_\sigma}(z)=z^\theta M(z)$, for $\theta>0$ and for some $M$ holomorphic and non-zero at $z=0$ (see \S \ref{subsec:local_theory}). 
Then the statement follows from the following computation in local coordinates:
$$
\int_{f(\eta)} \frac{dw}{w} = 
\int_\eta \dfrac{ \theta z^{\theta -1}M(z)+z^\theta M'(z)} {z^\theta M(z)}  dz = 
  \theta\int_\eta \dfrac{ dz} {z}  + \int_\eta \dfrac{ M'(z)} {M(z)}  dz   =  2\pi i \theta .
$$
where the second integral vanishes, because $M$ is holomorphic, and $\eta$ can be chosen to be small enough to enclose $z=0$ but no zero of $M$. % (recall $M(0)\neq 0$).
\endproof

\begin{remark}\label{rem:signing}
When the exponent (equivalently the reduced exponent) is not zero, a choice of a sign is called a \textit{signing} of the projective structure at that puncture, and can be used to define a framing from the holonomy representation (cf. \cite{AB18,G19}). This is in general an arbitrary choice. However, as observed in Corollary~\ref{cor:framed_hol}, continuously extending the developing map to the punctures always provides a canonical framing for structures in $\Modulit(\Sigma)$.	
\end{remark}

	%%%%%%%%%%%%%%%%%%%%%%%%%%%%%%%%%
	%%%%%%%%%%%%%%%%%%%%%%%%%%%%%%%%%
	%%%%%%%%%%%%%%%%%%%%%%%%%%%%%%%%%
	
	\section{Structures on the thrice--punctured sphere}\label{sec:thrice_punct_sphere}
	
	In this chapter we prove Theorems~\hyperref[thm:framed hol grafting]{A} and~\hyperref[thm:mainthm]{B} about grafting structures on the \emph{thrice--punctured sphere} $S := \SS^2 \setminus\{x_\alpha,x_\beta,x_\gamma\}$. 
	This is the oriented topological space obtained from the $2$--dimensional unit sphere $\SS^2$ by removing three distinct points $\{x_\alpha,x_\beta,x_\gamma\} \subset \SS^2$. The points $\{x_\alpha,x_\beta,x_\gamma\}$ are the \emph{punctures} of $S$. 
	(For the easier case of the twice-punctured sphere we refer the reader back to Remark~\ref{rem:twice}.)
	The \emph{fundamental group} $\pi_1(S)$ of $S$ is isomorphic to the free group on two generators $\FF_2$. Once and for all we fix the presentation
	$$
	\pi_1(S) = \langle \alpha,\beta,\gamma \ | \ \alpha\beta\gamma=1 \rangle \cong \FF_2,
	$$
	where each generator $\delta \in \{\alpha,\beta,\gamma\}$ can be represented by a peripheral loop (also denoted by $\delta$) around $x_\delta$, oriented to travel around the puncture in the anticlockwise direction. Furthermore, we denote by $\End_\delta \in \Ends(\Sw)$ the end in the end-extended universal cover $\Swe$ of $S$, that is fixed by $\delta$.
	
	In this setting, we observe that $\Modulic(S)$ is the space of complex projective structures whose underlying conformal structure is that of $\cp\setminus \{0,1,\infty\}$. 
	The $\pslc$--character variety can be explicitly described (see~\cite[Remark 4.4]{HP04} for details).
	A conjugacy class of representations is said to be \emph{non-degenerate relatively elliptic} if it is the class of a non-degenerate relatively elliptic representation. 
	It follows from Theorem~\hyperref[thm:tamerelell=merodoublenonint]{E} and \cite[Theorem 1.1]{G19} that any non-degenerate relatively elliptic conjugacy class arises from the holonomy of a structure in $\Modulit(S)$. We will see that the structure can be chosen to be of a special type (cf. Corollary~\ref{cor:gen_rel_ell}).
	
	\begin{remark}
	A relatively elliptic representation of $\pi_1(S)$ is degenerate if and only if its image is a subgroup of rotations around two fixed points, i.e. a group of coaxial rotations (cf. \cite[\S2.4]{G19}).
	\end{remark} 
	
	The main result of this chapter is a complete description of $\Modulit(S)$. We begin in \S\ref{subsec:tri_structures} by constructing some structures in $\Modulit(S)$, called \emph{triangular structures}, which will be our key examples. Then in \S\ref{subsec:grafting theorem} we show that $\Modulit(S)$ is precisely the space of complex projective structures obtained by grafting triangular structures.

	%%%%%%%%%%%%%%%%%%%%%%%%%%%%%%%%%
	%%%%%%%%%%%%%%%%%%%%%%%%%%%%%%%%%
	%%%%%%%%%%%%%%%%%%%%%%%%%%%%%%%%%
	
	\subsection{Triangular structures}\label{subsec:tri_structures}
	
	In this section we construct a family of structures in  $\Modulit(S)$ which will be the main reference example for the rest of the paper.
	
\begin{figure}[t]
		\centering
		\includegraphics[width=\textwidth]{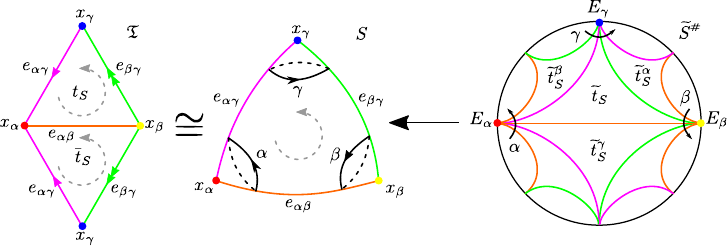}
		\caption{The ideal triangulation $\mathfrak{T}$ of the thrice--punctured sphere $S$, and its lift to the end-extended universal cover $\Swe$.}
		\label{fig:sphere_triangulation}
	\end{figure}

	First, we fix the following \emph{ideal triangulation} $\Tri$ of $S$ (cf. Figure~\ref{fig:sphere_triangulation}). For every distinct pair $\delta,\delta' \in \{\alpha,\beta,\gamma\}$, let $e_{\delta\delta'}$ be a simple arc on $S$ from $x_\delta$ to $x_{\delta'}$. The collection of arcs $\{e_{\alpha\beta},e_{\beta\gamma},e_{\alpha\gamma}\}$ are the \emph{ideal edges} of $\Tri$, and subdivide $S$ into two \emph{ideal triangles} $t_S$ and $\overline{t}_S$. The orientation of $S$ induces an orientation on $t_S$ (resp. $\overline{t}_S$) so that the punctures are ordered as $(x_\alpha,x_\beta,x_\gamma)$ (resp. $(x_\alpha,x_\gamma,x_\beta)$) on its boundary. The ideal triangulation $\Tri$ lifts to a triangulation $\wTri$ of $\Swe$. We notice that the restriction of $\wTri$ to $\Sw$ is an ideal triangulation of $\Sw$. We denote by $\widetilde{t}_S$ the unique triangle in $\wTri$ with vertices $\{\End_\alpha,\End_\beta,\End_\gamma\}$, and by $\widetilde{t}_S^\delta$ the unique triangle adjacent to $\widetilde{t}_S$ that does not have $\End_\delta$ as its vertex. It is easy to check that $\widetilde{t}_S$ projects onto $t_S$, while $\left\{\widetilde{t}_S^\alpha,\widetilde{t}_S^\beta,\widetilde{t}_S^\gamma \right\}$ all project onto $\overline{t}_S$.
	
	Recall that $\triangle \subset \RR^3$ is the standard $2$--simplex (cf.~\S\ref{subsec:tri_immersions}). Let $\tau : \triangle \rightarrow \cp$ be a non-degenerate triangular immersion, with vertices $(V_a,V_b,V_c)$ and angles $(a,b,c)$. Let $\Circles_\tau = (\C_{ab},\C_{bc},\C_{ac})$ be the configuration of circles determined by $\tau$, defined such that $V_x,V_y \in \C_{xy}$, for all distinct pairs $x,y \in \{a,b,c\}$. From Corollary \ref{cor:circles_reps} we have a relatively elliptic representation associated to $\Circles_\tau$ given by
	$$
	\rho_{\tau} := \rho_{\Circles_\tau}:\pi_1(S)\to \pslc,
	$$
	$$
	\rho_{\tau}(\alpha):=J_{ac}J_{ab}, \quad \rho_{\tau}(\beta):=J_{ab}J_{bc}, \quad \rho_{\tau}(\gamma):=J_{bc}J_{ac},
	$$
	where $J_{xy}$ denotes the reflection of $\cp$ in $\C_{xy}$. Notice that, if $\tau$ embeds onto a Euclidean, hyperbolic or spherical triangle with angles rational multiples of $\pi$, then the image of this representation is a discrete Euclidean, hyperbolic or spherical group; however a generic choice of $\tau$ results in a non-discrete subgroup of $\pslc$.

	The \emph{triangular structure} $\sigma_\tau \in \Moduli(S)$  associated to the triangular immersion $\tau:\triangle\to \cp$ is the structure defined by the developing pair $(\dev_{\tau},\rho_{\tau})$, where the developing map is constructed as follows. Recall that $(V_1,V_2,V_3)$ are the vertices of $\triangle$. Consider the following maps:
	\begin{enumerate}
	    \item $\varphi:\widetilde{t}_S\to \triangle$: the unique simplicial map mapping $(\End_\alpha,\End_\beta,\End_\gamma)$ to $(V_1,V_2,V_3)$;
	    
	    \item $\varphi^\gamma:\widetilde{t}_S^\gamma \to \triangle$: the unique simplicial map mapping $(\End_\beta,\End_\alpha,\End_{\beta\gamma\beta^{-1}})$ to $(V_1,V_2,V_3)$;
	    
	    \item $\iota : \triangle \to \triangle$: the unique (orientation reversing) simplicial map mapping $(V_1,V_2,V_3)$ to $(V_2,V_1,V_3)$;
	    
	    \item $\tau^\gamma := J_{ab} \circ \tau \circ \iota$ be the \emph{triangular immersion conjugate to $\tau$}, mapping $(V_1,V_2,V_3)$ to $(V_b,V_a,J_{ab}(V_c))$.
	\end{enumerate}	
	Then we define
	$$
	\restr{(\devt)_{\tau}}{\widetilde{t}_S} :=\tau\circ \varphi, \qquad \text{and} \qquad   \restr{(\devt)_{\tau}}{\widetilde{t}_S^\gamma}:=\tau^\gamma \circ \varphi^\gamma.
	$$
	Since this defines $(\devt)_{\tau}$ on a fundamental domain for the action of $\pi_1(S)$ on $\Swe$, we can then extend it by equivariance with respect to the representation $\rho_{\tau}$ to obtain a global $(\devt)_\tau:\Swe\to \cp$. The developing map $\dev_\tau$ is the restriction of $(\devt)_\tau$ to $\Sw$. Notice that, when $\tau$ is an embedding, this is the pillowcase structure obtained by doubling $\tau(\triangle)$.
	
	By construction, triangular structures are non-degenerate, tame and their holonomy representations are relatively elliptic. We record this in the following lemma.

	\begin{lemma}\label{lem:triangular_structures}
		Let $\tau$ be a non-degenerate triangular immersion and let $\sigma_\tau$ be the associated triangular structure. Then $\sigma_\tau \in \Modulit(S)$.
	\end{lemma}
	
\begin{figure}[t]
		\centering
		\includegraphics[width=\textwidth]{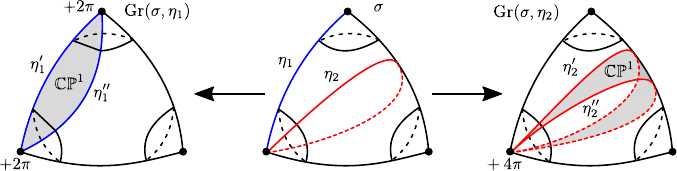}
		\caption{An edge-grafting and a core-grafting on a structure $\sigma$.}
		\label{fig:edge_core_grafting}
	\end{figure}

	Triangular immersions that are especially simple, e.g. embeddings, carry some obvious curves that one can graft along, namely the edges $e_{\delta\delta'}$ of the triangulation $\Tri$. Other graftable curves are those joining one puncture to itself by crossing the triangle. We introduce the following terminology, motivated by these observations (see \S\ref{subsec:tame_rel_ell} for the general definition of this surgery).
	Let $\sigma \in \Modulit(S)$ and let $\eta:I\to S$ be a graftable curve. The grafting along $\eta$ will be called an \textit{edge-grafting}  if $\eta$ joins two different punctures, and a  \textit{core-grafting} if it starts and ends at the same puncture and separates $S$ into two punctured disks. The inverse surgery will be called \textit{edge-degrafting} and \textit{core-degrafting} respectively (cf. Figure~\ref{fig:edge_core_grafting}).

\begin{figure}[t]
    \centering
    \includegraphics[width=\textwidth]{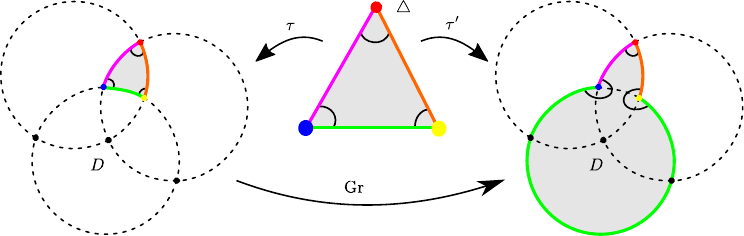}
    \caption{An edge-grafting on an embedded structure $\sigma$.}
    \label{fig:edge_grafting_example}
    \end{figure}

	\begin{example}\label{ex:edge_graft}
    Some embedded triangular structures allow for an easy description of edge-grafting. Let $\tau,\tau':\triangle\to \cp$ be two triangular embeddings such that $\sigma_{\tau'}$ differs from $\sigma_{\tau}$ by the insertion of a disk $D$ along one of the edges (cf. Figure~\ref{fig:edge_grafting_example}, first two pictures). Then $\sigma_{\tau'}$ is isomorphic to the structure obtained by edge-grafting $\sigma_\tau$ along that edge. 
    Indeed reflecting in the edges of $\tau'(\triangle)$ we obtain a copy of $\cp$ obtained by the union of $D$ and its complement. Since $D$ is included in $\tau'(\triangle)$, its complement is contained in a suitable reflection of it; the union of $D$ and its complement  gives precisely a grafting region on $\sigma_{\tau'}$. This grafting procedure can be iterated by thinking of immersions as membranes spread over $\cp$, obtained by including additional disks across the edges that are being grafted. This is a particularly concrete way of thinking about edge-grafting triangular structures.
	\end{example}
	
	A triangular structure is said to be \emph{Euclidean/hyperbolic/spherical atomic} if it comes from a Euclidean/hyperbolic/spherical atomic triangular immersion (cf. end of \S\ref{subsec:tri_immersions}). The terminology is motivated by the main theorem (cf. Theorem~\hyperref[thm:mainthm]{B}), which states that every tame and relatively elliptic $\cp$--structure is obtained by grafting an atomic structure. 
	
	\begin{lemma}\label{lem:grafting_atomic}
	   Let $\sigma$ be an atomic triangular structure with indices $\Ind_{\sigma} := (2a,2b,2c)$. Let $e_{\delta\delta'}$ be the edge of the triangle of $\Tri$ in $S$ connecting the two distinct punctures $x_\delta$ and $x_{\delta'}$. Let $e_{\delta}$ be a simple ideal arc in $S$ connecting the puncture $x_\delta$ to itself by crossing the edge opposite to $x_\delta$. For $G_{\alpha\beta},G_{\alpha\gamma},G_{\alpha},G_{\beta\gamma} \in \NN$, consider the formal sums
	   $$
	   \eta := G_{\alpha\beta} e_{\alpha\beta} + G_{\alpha\gamma}e_{\alpha\gamma} + G_{\alpha}e_{\alpha}, \quad \text{and} \quad \eta' :=  G_{\alpha\beta}e_{\alpha\beta} + G_{\alpha\gamma}e_{\alpha\gamma} + G_{\beta\gamma}e_{\beta\gamma}.
	   $$
	   If $\sigma$ is spherical or hyperbolic, then $\sigma$ is graftable along both $\eta$ and $\eta'$, up to small deformations. If $\sigma$ is Euclidean and we further assume that that $a \in (0,3\pi)$ while $b,c \in (0,\pi)$. Then
	\begin{enumerate}
	    \item \label{item:graftable_1} If $a \in (0,\pi)$ and $-a+b+c = \pi$, then $\sigma$ is graftable along $\eta'$, but not along any arc isotopic to $e_{\alpha}$.
	    \item \label{item:graftable_2} If $a \in (\pi,2\pi)$ and $a-b-c = \pi$, then $\sigma$ is graftable along $\eta$, but not along any arc isotopic to $e_{\beta\gamma}$.
	    \item \label{item:graftable_3} If $a \in (2\pi,3\pi)$, then $\sigma$ is graftable along $\eta$, but not along any arc isotopic to $e_{\beta\gamma}$.
	    \item \label{item:graftable_5} Otherwise $\sigma$ is graftable along both $\eta$ and $\eta'$.
	\end{enumerate}
	\end{lemma}
	
\begin{figure}[t]
    \centering
    \includegraphics[width=0.9\textwidth]{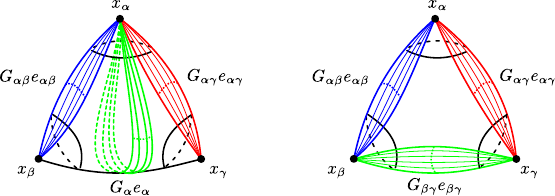}
    \caption{The multi-curves $\eta$ (on the left), and $\eta'$ (on the right).}
    \label{fig:multi_grafting}
    \end{figure}

	\proof
	We begin by noticing that $\eta$ (resp. $\eta'$) can be realized as a group of pairwise disjoint arcs in $S$ (cf. Figure~\ref{fig:multi_grafting}), therefore we only need to check that $\sigma$ is graftable once along each arc (cf. Remark~\ref{rem:multiple_grafting}).
	
	If $\sigma$ comes from a triangular immersion $\tau$ supported by a spherical configuration, then $\sigma$ is graftable along both $\eta$ and $\eta'$ because the triangular immersion $\tau$ is an embedding (cf. Figures~\ref{fig:tiny_angles_hs} (right) and ~\ref{fig:small_angles_hs} (right)), hence each simple ideal arc develops injectively into $\cp$.
	
\begin{figure}[t]
    \centering
    \includegraphics[width=\textwidth]{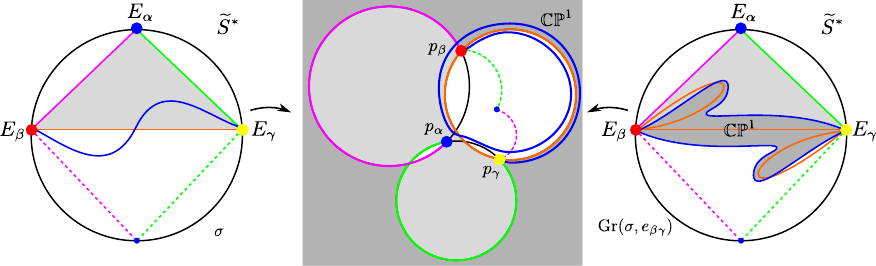}
    \caption{An edge-grafting on the hyperbolic atomic structure coming from a hyperbolic atomic triangular immersion as in~Figure~\ref{fig:small_angles_h2} (2)(i).}
    \label{fig:edge_graft_degenerate_hyp}
    \end{figure}

	Similarly, if $\tau$ is supported by a hyperbolic configuration, then $\tau$ is an embedding unless it is as in Figure~\ref{fig:small_angles_h2} (2)(i). These are immersions where one angle is in $(\pi,2\pi)$, say for example $a$, and $a-b-c>\pi$. In these situations, the edge $e_{\beta\gamma}$ (opposite to the large angle $a$) is not graftable on the nose, as the developing map develops it surjectively to a circle. However any arbitrarily small deformation of it is graftable (cf. Figure~\ref{fig:edge_graft_degenerate_hyp}).
	
	Finally, suppose that $\tau$ is supported by a Euclidean configuration. Here we further assume $a \in (0,3\pi)$ while $b,c \in (0,\pi)$, namely that if there is an angle larger than $\pi$, then it is $a$. Here we have an issue only when a puncture is mapped to the common intersection point $y$ of the Euclidean configuration. If $a \in (0,\pi)$ and $-a+b+c = \pi$ (case~\eqref{item:graftable_1}), then the puncture $x_\alpha$ develops to $y$ and it is not possible to core-graft along any arc isotopic to $e_\alpha$ (cf. Figure~\ref{fig:tiny_angles_e} (right)). On the other hand, every edge is injectively developed, and therefore $\sigma$ is graftable along $\eta'$. If $a \in (\pi,2\pi)$ and $a-b-c = \pi$ (case~\eqref{item:graftable_1} and Figure~\ref{fig:small_angles_e2} (2) (ii)), then both $x_\beta,x_\gamma$ are mapped to $y$, thus $\sigma$ is not graftable along any arc isotopic to $e_{\beta\gamma}$, and in particular along $\eta'$. However, $e_\alpha$ is injectively developed, hence $\sigma$ is graftable along $\eta$. Case~\eqref{item:graftable_3} is similar to the previous one (cf. Figure~\ref{fig:additional_euclidean}).
	The remaining Euclidean cases are embeddings where $x_\alpha$ never maps to $y$, hence all relevant arcs are injectively developed.
	\endproof

	We conclude this section with a simple consequence of Lemma~\ref{lem:triangular_structures}, namely that almost every non-degenerate framed relatively elliptic representation is the framed holonomy representation of an atomic triangular structure. 
	Recall that a framing of a representation $\rho$ is a $\rho$--equivariant map $\mathcal F:\Ends(\Sigmaw) \to \cp$ from the space of ends to $\cp$, and that for structures in $\Modulit(S)$ there is a canonical framing given by a continuous extensions of the developing map (cf. Corollary~\ref{cor:framed_hol}). 
	We remark that  \cite[Theorem 1.2]{G19} states that a non-degenerate framed representation is the holonomy of a signed meromorphic projective structure with respect to some framing, while here we realize these framed representations with respect to this canonical framing (compare the discussion in Remark~\ref{rem:signing}).
	To simplify the statement of the following result, we say that a framing $\mathcal F$ is \emph{pathological} if $\mathcal F$ maps the entire set of ends to a single point. In our context, the holonomy representation of a triangular structure is pathological if and only if the underlying configuration of circles is Euclidean and the framing consists only of the point at infinity. Therefore the holonomy representation of an atomic triangular structure is never pathological. Note that a pathological framing is not considered degenerate according to the definition in \S\ref{subsec:tame_rel_ell}.
	
	\begin{corollary}\label{cor:gen_rel_ell}
		Every non-degenerate framed relatively elliptic representation that is not pathological is the framed holonomy representation of an atomic triangular structure. In particular, $\Chart(S)=\Hol (\Modulit(S))$.
	\end{corollary}
	\proof
	Suppose $\rho$ is a non-degenerate relatively elliptic representation, with a non-pathological framing $\mathcal F$. Then $(\rho(\alpha),\rho(\beta),\rho(\gamma))$ is an ordered triple of elliptic transformations with trivial product. As $\rho$ is non-degenerate, $(\rho(\alpha),\rho(\beta),\rho(\gamma))$ share at most one common fixed point.  By Corollary~\ref{cor:circles_reps}, there is a unique non-degenerate configuration of circles $\Circles := (\C_{ab},\C_{bc},\C_{ac})$ associated to $(\rho(\alpha),\rho(\beta),\rho(\gamma))$. By construction
	$$
	p_\alpha := \mathcal F(\End_\alpha) \in \C_{ab} \cap \C_{ac}, \quad 
	p_\beta := \mathcal F(\End_\beta) \in \C_{ab} \cap \C_{bc}, \quad \text{and} \quad
	p_\gamma := \mathcal F(\End_\gamma) \in \C_{bc} \cap \C_{ac}.
	$$
	We are going to show that there is an atomic triangular immersion $\tau$ supported by $\Circles$, with vertices $(p_\alpha,p_\beta,p_\gamma)$. As a consequence, the framed holonomy representation of its associated triangular structure $\sigma_\tau$ is $(\rho, \mathcal F)$, proving the first part of the corollary. If $\Circles$ is a spherical configuration, the points $(p_\alpha,p_\beta,p_\gamma)$ are the vertices of a unique triangular region $R$ in $\cp \setminus \Circles$. Depending on the cyclic order of $(p_\alpha,p_\beta,p_\gamma)$ on the boundary of $R$, we either take $\tau$ to map onto $R$, or to map onto the complement of $R$ in a disk (cf. Figure~\ref{fig:small_angles_hs} on the right). If $\Circles$ is a hyperbolic configuration, we refer to Table~\ref{tab:hyperbolic} to check that any framing is realized by at least one triangular immersion $\tau$. Finally, Table~\ref{tab:Euclidean} shows that any framing that is not pathological, namely $(-,-,-)$ and $(-,-,-)^*$, can be realized by at least one triangular immersion $\tau$.
	
	The last statement of the corollary follows from the observation that every non-degenerate relatively elliptic conjugacy class $[\rho]$ has a class representative $\rho$ that can be framed with a non-degenerate and non-pathological framing.
	\endproof

\subsection{Grafting Theorems A and B}\label{subsec:grafting theorem}
We are now ready to prove the main results about the Grafting Conjecture. A key step will be being able to recognize structures based on their indices, which we are able to do thanks to the description of $\Modulit(S)$ in terms of meromorphic differentials (cf. Theorem~\hyperref[thm:tamerelell=merodoublenonint]{E}).

Up to isomorphism, there is a unique complex structure on the thrice--punctured sphere, namely that of $\cp\setminus \{0,1,\infty\}$. The space of meromorphic quadratic differentials with double poles at $0,1$ and $\infty$ can be described as follows:
$$
\left\lbrace q_\Theta=\left(\frac{1-\theta_1^2}{2z^2}+\frac{1-\theta_2^2}{2(z-1)^2}+\frac{\theta_1^2+\theta^2_2-\theta_3^2-1}{2z(1-z)}\right) dz^2 \ | \ \Theta=(\theta_1,\theta_2,\theta_3)\in \CC^3 \right\rbrace .
$$
A direct computation shows that $q_\Theta$ has double poles at $0,1,\infty$  with reduced exponents $\theta_1,\theta_2,\theta_3$, respectively.
% It follows from the discussion in \S\ref{subsec:meromorphic_structures} that the structures in $\Modulit(S)$ are precisely those meromorphic projective structures defined by the meromorphic quadratic differentials $q_\Theta$ with reduced exponents $\theta_i\in \RR\setminus \ZZ$ for $i=1,2,3$. 
In particular, the indices of the structure defined by the differential $q_\Theta$ are $(2\pi \left|\theta_1\right|,2\pi \left|\theta_2\right|,2\pi \left|\theta_3\right|)$ (cf. Lemma~\ref{lem:index_exponent}). Therefore we obtain the following statement.

\begin{proposition}\label{prop:indices_determine_structure}
If $\sigma,\sigma'\in \Modulit(S)$ have the same indices, then $\sigma=\sigma'$.
\end{proposition}
\proof
By Theorem~\hyperref[thm:tamerelell=merodoublenonint]{E} we know that $\sigma=\sigma_q,\sigma'=\sigma_{q'}$ for some meromorphic differentials  $q,q'\in \Quadddp S$, with real non-integer reduced exponents at each puncture. 
Since the index at each puncture is the same, by Lemma \ref{lem:index_exponent} the exponent at each puncture is also the same (up to sign). So $q,q'$ have the same leading coefficient at each puncture, but this determines them completely, so $q=q'$.
\endproof

Notice that the developing maps of structures obtained with $\theta_i\in (0,1)$  correspond to Schwarz triangle maps. The special cases in which $\theta_i=\frac{1}{p_i},p_i\in \ZZ$ correspond to the classic uniform  tilings of the sphere, Euclidean or hyperbolic plane.
In the general case $\theta_i\in \RR \setminus \ZZ$, the associated holonomy representations are not discrete, and the groups are not isomorphic to triangle groups. 

A direct application of Proposition~\ref{prop:indices_determine_structure} to Lemmas~\ref{lem:triangular_reg_for_small_angles} and~\ref{lem:additional_euclidean} allows us to easily characterize atomic structures through their indices.

\begin{lemma}\label{lem:angles_of_atomic_structures}
    Let $\sigma \in \Modulit(S)$ with indices $(2a,2b,2c)$.
    Then $\sigma$ is atomic if and only if (up to relabelling the punctures):
    \begin{enumerate}
        \item either $a \in (0,2\pi)$ and $b,c \in (0,\pi)$;
        \item or $a \in (2\pi,3\pi)$ and $b,c \in (0,\pi)$ and $a-b-c = \pi$.
    \end{enumerate}
\end{lemma}
	
\proof
    Atomic structures are defined in such a way that their indices satisfy the above conditions (cf. Lemmas~\ref{lem:triangular_reg_for_small_angles} and~\ref{lem:additional_euclidean}). But more importantly, every triple of numbers $(2a,2b,2c)$ satisfying those conditions is the triple of indices of an atomic structure, see for example Tables~\ref{tab:hyperbolic},~\ref{tab:spherical} and ~\ref{tab:Euclidean}. The fact that there are no other structures with those indices follows by Proposition~\ref{prop:indices_determine_structure}. 
\endproof
	
	As observed in Corollary \ref{cor:framed_hol}, the holonomy representation of a structure in $\Modulit(S)$ carries a natural framing, given by the extension of the developing map to the punctures. Edge-grafting and core-grafting do not change the holonomy representation, nor this framing (see Lemma \ref{lem:grafting_properties}). 
	
	\begin{maintheoremc}{B}\label{thm:mainthm}
	Every $\sigma \in \Modulit(S)$ is obtained by a sequence of edge- and core-graftings on an atomic triangular structure with the same framed holonomy.
	\end{maintheoremc}
	
	\proof
	
	Let $\sigma \in \Modulit(S)$, let   $2a :=\Ind_{\sigma}(x_\alpha)$, $2b :=\Ind_{\sigma}(x_\beta)$, $2c :=\Ind_{\sigma}(x_\gamma)$ be its indices. 
	Without loss of generality we can assume that $a\geq b\geq c$. Indeed we can rename the punctures so that $\Ind_{\sigma}(x_\alpha)$ is the largest index, and the case where $a \geq c \geq b$ follows by a similar argument.

	Let $k_a=\left \lfloor{\frac{a}{\pi}}\right \rfloor$,
	$k_b=\left \lfloor{\frac{b}{\pi}}\right \rfloor$,
	$k_c=\left \lfloor{\frac{c}{\pi}}\right \rfloor  \in \NN$. We are going to reduce the triple $(a,b,c)$ to a triple $(a',b',c')$ by subtracting as many integer multiple of $\pi$ as possible in a certain controlled way, until $(a',b',c')$ satisfies the conditions of Lemma~\ref{lem:triangular_reg_for_small_angles}, that is
	\begin{equation}\label{eq:small_indices}
	a'\in (0,\pi) \cup (\pi,2\pi), \quad \text{and} \quad b',c' \in (0,\pi).
	\end{equation}
 
    We distinguish two cases:
    \begin{itemize}
	    \item[(i)] If $k_a\geq k_b+k_c$, let  
	    $$
	    G_{\alpha\gamma} := k_c, \quad G_{\alpha\beta} := k_b, \quad G_{\alpha} := \left \lfloor{\frac{k_a-(k_b+k_c)}{2}}\right \rfloor, \quad G_{\beta\gamma} := 0.
	    $$
	    
	    \item[(ii)] If $k_a < k_b+k_c$, let $L := k_a - k_b$, $L' := k_c + k_b - k_a$ and
	    $$
	    G_{\alpha\gamma} := L + \left \lfloor \frac{L'}{2} \right \rfloor , \quad G_{\alpha\beta} := k_b - \left \lceil \frac{L'}{2} \right \rceil, \quad G_{\alpha} := 0, \quad G_{\beta\gamma} := \left \lceil \frac{L'}{2} \right \rceil.
	    $$
	\end{itemize}
	    Either way, let 
	$$
	a':=a-\pi (G_{\alpha\gamma}+G_{\alpha\beta}+2 G_{\alpha}), \ b':=b-\pi( G_{\beta\gamma}+G_{\alpha\beta}), \ c':=c-\pi (G_{\alpha\gamma}+ G_{\beta\gamma}).
	$$
    It is easy to check that $G_{\alpha\gamma}, G_{\alpha\beta}, G_{\alpha}, G_{\beta\gamma} \geq 0$, and
    $$
    G_{\beta\gamma}+G_{\alpha\beta} = k_b, \quad G_{\alpha\gamma}+ G_{\beta\gamma} = k_c, \quad G_{\alpha\gamma}+G_{\alpha\beta}+2 G_{\alpha} \in \{k_a, k_a - 1\},
    $$
    therefore~\eqref{eq:small_indices} is satisfied, and by Lemma~\ref{lem:triangular_reg_for_small_angles} there is a triangular immersion $\tau$ with angles $\left(a',b',c'\right)$. Let $\sigma_\tau$ be the associated triangular structure. By construction $\sigma_\tau$ is atomic with indices $(2a',2b',2c')$, thus it is left to check if $\sigma_\tau$ grafts to $\sigma$.
    
    Recall we have fixed an ideal triangulation $\Tri$ of $S$. 
    Let $e_{\delta\delta'}$ be the edges of $\Tri$  connecting the two distinct punctures $x_\delta$ and $x_{\delta'}$. Let $e_{\delta}$ be a simple ideal arc in $S$ connecting the puncture $x_\delta$ to itself by crossing the edge opposite to $x_\delta$. Consider the multi-curve
	   $$
	   \mu := G_{\alpha\beta} e_{\alpha\beta} + G_{\alpha\gamma} e_{\alpha\gamma} + G_{\alpha} e_{\alpha} + G_{\beta\gamma} e_{\beta\gamma}.
	   $$
	If $\mu$ is graftable then grafting $\sigma_\tau$ along $\mu$ would yield a structure with indices $(2a,2b,2c)$ and the same framed holonomy as $\sigma_\tau$ (cf. Lemma~\ref{lem:grafting_properties}). It follows from  Proposition~\ref{prop:indices_determine_structure} that $\sigma =  \Gr(\sigma_\tau,\mu)$, so it is left to check if $\sigma_\tau$ is graftable along $\mu$.
	   
    Depending on the above cases, we remark that at least one of $G_{\beta\gamma},G_{\alpha}$ is $0$, hence $\mu$ is either $\eta$ or $\eta'$ in the notation of Lemma~\ref{lem:grafting_atomic}.
    
    If $G_{\beta\gamma} = G_{\alpha} = 0$ then $\mu = \eta = \eta'$ and every atomic triangular structure $\sigma_\tau$ is graftable along $\mu$.
    
    If $G_{\alpha} > 0$ then $G_{\beta\gamma}= 0$ and $\mu = \eta$. Lemma~\ref{lem:grafting_atomic} covers every case except the Euclidean case where $a' \in (0,\pi)$ and $-a'+b'+c' = \pi$. In this case we must consider a different atomic structure $\sigma_\tau'$ and curve $\mu'$, as $\sigma_\tau$ is not graftable along $e_\alpha$. Let
    $$
    a'':=a'+2\pi, \qquad \ b'':=b', \qquad c'':=c',
    $$
    $$
    G_\alpha' := G_\alpha - 1, \quad \text{and} \quad \mu' := G_{\alpha\beta} e_{\alpha\beta} + G_{\alpha\gamma} e_{\alpha\gamma} + G_{\alpha}' e_{\alpha}.
    $$
    By construction $a'' \in (2\pi,3\pi)$, $b'',c'' \in (0,\pi)$ and $a''-b''-c'' = \pi$, therefore there is an atomic triangular structure $\sigma_\tau'$ with indices $(2a'',2b'',2c'')$ (cf. Lemma~\ref{lem:additional_euclidean}). Furthermore, the structure $\sigma_\tau'$ is graftable along $\mu'$ (cf. Lemma~\ref{lem:grafting_atomic}). Grafting $\sigma_\tau'$ yields a structure with indices $(2a,2b,2c)$, which must be $\sigma$ by Proposition~\ref{prop:indices_determine_structure}, concluding this case.
    
    Lastly, suppose that $G_{\beta\gamma} > 0$. This time $G_{\alpha} = 0$ and $\mu = \eta'$. Recall that $a' < 2\pi$, hence the only case that is not covered by Lemma~\ref{lem:grafting_atomic} is the Euclidean case where $a' \in (\pi,2\pi)$ and $a'-b'-c' = \pi$. We are once again forced to consider a different atomic structure as $\sigma_\tau$ is not graftable along $e_{\beta\gamma}$. Let
    $$
    a'':=a'-\pi, \qquad \ b'':=b'+\pi, \qquad c'':=c',
    $$
    $$
    G_{\alpha\gamma}' := G_{\alpha\gamma} + 1, \quad G_{\beta\gamma}' := G_{\beta\gamma} - 1, \quad \text{and} \quad \mu' := G_{\alpha\beta} e_{\alpha\beta} + G_{\alpha\gamma}' e_{\alpha\gamma} + G_{\beta\gamma}' e_{\beta\gamma}.
    $$
    By construction $b'' \in (\pi,2\pi)$, $a'',c'' \in (0,\pi)$ and $-a''+b''+c'' = \pi$, therefore there is an atomic triangular structure $\sigma_\tau'$ with indices $(2a'',2b'',2c'')$ (cf. Lemma~\ref{lem:triangular_reg_for_small_angles}). The structure $\sigma_\tau'$ is graftable along $\mu'$ according to Lemma~\ref{lem:grafting_atomic} part $(5)$ applied to the triple $(b'',c'',a'')$. Once again, grafting $\sigma_\tau'$ along $\mu'$ yields a structure with indices $(2a,2b,2c)$, which must be $\sigma$ by Proposition~\ref{prop:indices_determine_structure}, concluding the proof.

\endproof
	
	Theorem~\hyperref[thm:mainthm]{B} has two interesting consequences. The first is the promised characterization of atomic structures in terms of grafting.

	\begin{corollary}\label{cor:atomic=ungrafted}
	   A structure $\sigma \in \Modulit (S)$ is atomic if and only if it is not degraftable.
	\end{corollary}
	
	\proof
	For one implication, let $\sigma$ be a structure which cannot be degrafted. Then by Theorem \hyperref[thm:mainthm]{B} it must be atomic.
	
	For the reverse implication, let $\sigma$ be atomic. Suppose by contradiction that $\sigma$ was degraftable to some structure $\sigma'$. Recall that core-grafting increases one index by $4\pi$ and edge-grafting increases two indices by $2\pi$. Then $\sigma$ cannot be one of the atomic structures coming from the atomic triangular immersions of Lemma~\ref{lem:triangular_reg_for_small_angles}, as its indices would be too small. It follows that $\sigma$ is the atomic triangular structure associated to an atomic triangular immersion $\tau$ from Lemma~\ref{lem:additional_euclidean}. Without loss of generality we may assume that the largest index of $\sigma$ is at $x_\alpha$, while the other two are less than $2\pi$, so that
	$$
	\Ind_{\sigma}(x_\alpha) \in (4\pi,6\pi),\quad \Ind_{\sigma}(x_\beta),\Ind_{\sigma}(x_\gamma) \in (0,2\pi), \quad \text{and} \quad \Ind_{\sigma}(x_\alpha) - \Ind_{\sigma}(x_\beta)  - \Ind_{\sigma}(x_\gamma) = 2\pi.
	$$
	Then $\sigma$ cannot be obtained by edge-grafting $\sigma'$, and the only option is that $\sigma'$ is a core-degrafting at $x_\alpha$ on $\sigma$. In particular $\Ind_{\sigma'}(x_\alpha) = \Ind_{\sigma}(x_\alpha) - 4\pi \in (0,2\pi)$ and
	$$
	\Ind_{\sigma'}(x_\alpha),\Ind_{\sigma'}(x_\beta),\Ind_{\sigma'}(x_\gamma) \in (0,2\pi), \quad \text{and} \quad -\Ind_{\sigma'}(x_\alpha) + \Ind_{\sigma'}(x_\beta)  + \Ind_{\sigma'}(x_\gamma) = 2\pi.
	$$
	It follows that $\sigma'$ is an atomic triangular structure (cf. Lemma~\ref{lem:angles_of_atomic_structures}), coming from a triangular immersion $\tau'$ enclosed in a Euclidean configuration (cf. Lemma~\ref{lem:tiny_angles}). But this is impossible because $\sigma'$ is not core-graftable at $x_\alpha$ (cf. Lemma~\ref{lem:grafting_atomic} part~\eqref{item:graftable_1}), giving the desired contradiction.
	\endproof

    Next, we obtain that edge-grafting and core-grafting (together with the inverse operations) account for all the possible deformations that preserve the holonomy as a framed representation.

\begin{maintheoremc}{A}\label{thm:framed hol grafting}
  Two structures in $\Modulit(S)$ have the same framed holonomy if and only if it is possible to obtain one from the other by some combination of graftings and degraftings along ideal arcs.
\end{maintheoremc}

\proof
One direction is clear by Lemma \ref{lem:grafting_properties}. For the reverse implication, suppose $\sigma,\sigma' \in \Modulit (S)$ have the same framed holonomy. 
By Theorem \hyperref[thm:mainthm]{B}, the structure $\sigma$ (resp. $\sigma'$) can be degrafted to an atomic structure $\sigma_0$ (resp. $\sigma_0'$) having the same framed holonomy.

Let $\tau_0,\tau_0'$ be the atomic triangular immersions defining $\sigma_0,\sigma_0'$, with angles $(a_0,b_0,c_0)$ and  $(a_0',b_0',c_0')$, respectively. Since these structures have the same framed holonomy, up to conjugation we can assume that $\tau_0,\tau_0'$ are supported by the same configuration of circles $\Circles$ (cf. Corollary~\ref{cor:circles_reps}), and that $\tau_0(V_j)=\tau_0'(V_j)$, for $j=1,2,3$. By Corollary~\ref{cor:atomic_angles_relations} we are in one of the following two cases:
\begin{enumerate}
 \item \label{item:same_angles} $(a_0,b_0,c_0)=(a_0',b_0',c_0')$;
 \item \label{item:alg_graft_angles} $(a_0-a_0',b_0-b_0',c_0-c_0')=(\pi,-\pi,0)$ up to permutation.
\end{enumerate}
In the first case $\sigma_0$ and $\sigma_0'$ have the same indices, hence $\sigma_0 = \sigma_0'$ by Proposition~\ref{prop:indices_determine_structure}, and we are done. For the second case, let us fix the permutation
$$
(a_0-a_0',b_0-b_0',c_0-c_0')=(\pi,-\pi,0),
$$
as the other cases are similar. Then in particular $a_0,b_0' \in (\pi,2\pi)$ while $a_0',b_0,c_0,c_0' \in (0,\pi)$. Let $\sigma_1$ (resp. $\sigma_1'$) be the triangular structure obtained by grafting $\sigma_0$ along $e_{\beta\gamma}$ (resp. $\sigma_0'$ along $e_{\alpha\gamma}$). These structures exist by Lemma~\ref{lem:grafting_atomic} (with respect to $\eta'$), and they both have indices
$$
(2a_0,2b_0+2\pi,2c_0+2\pi)=(2a_0' + 2\pi,2b_0',2c_0'+2\pi).
$$
We explicitly observe that Lemma~\ref{lem:grafting_atomic} has only two cases in which $\eta'$ is not graftable, and a direct inspection of Table~\ref{tab:Euclidean} shows that those two structures are covered by the case $(a_0,b_0,c_0)=(a_0',b_0',c_0')$ above (see Remark~\ref{rem:rigid_signs}).
It follows that $\sigma_1 = \sigma_1'$ by Proposition~\ref{prop:indices_determine_structure}, completing the proof.
\endproof
	 
\begin{figure}
		\centering
		\includegraphics[width=0.4\textwidth]{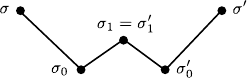}
		\caption{To prove Theorem~\hyperref[thm:framed hol grafting]{A} we find a path of graftings and degraftings from $\sigma$ to $\sigma'$, passing through atomic structures.}
    \label{fig:grafting_path}
\end{figure}

	\clearpage
\appendix

\section{Tables of atomic triangular immersions}

\begin{table}[h]
	\begin{center}
 \hspace*{-1.75cm}
		\begin{tabular}{|c|c|c|c|c|c|c|c|}
			\hline
			\multicolumn{3}{|c|}{Angles Range} &  &  & Target angles & & \\
			\hline
			$a$ & $b$ & $c$ & Conditions & Type & $(\hat a, \hat b, \hat c)$ & Signs & Figure \\
			\hline \hline
			$(0,\pi)$ & $(0,\pi)$ & $(0,\pi)$ & $a+b+c < \pi$ & H & $(a,b,c)$ & $(+,+,+)$ & Figure~\ref{fig:tiny_angles_hs} (left) \\
			\hline \hline
			$(0,\pi)$ & $(0,\pi)$ & $(0,\pi)$ & $a+\pi < b+c$ & H & $(a,\pi-b,\pi-c)$ & $(+,-,-)$ & Figure~\ref{fig:small_angles_hs} (left) \\
			\hline
			$(0,\pi)$ & $(0,\pi)$ & $(0,\pi)$ & $b+\pi < a+c$ & H & $(\pi-a,b,\pi-c)$ & $(-,+,-)$ & Figure~\ref{fig:small_angles_hs} (left) \\
			\hline
			$(0,\pi)$ & $(0,\pi)$ & $(0,\pi)$ & $c+\pi < a+b$ & H & $(\pi-a,\pi-b,c)$ & $(-,-,+)$ & Figure~\ref{fig:small_angles_hs} (left) \\
			\hline \hline
			$(\pi,2\pi)$ & $(0,\pi)$ & $(0,\pi)$ & $a+b+c > 3\pi$ & H & $(2\pi-a,\pi-b,\pi-c)$ & $(-,-,-)$ & Figure~\ref{fig:small_angles_h2} (1) \\
			\hline
			$(0,\pi)$ & $(\pi,2\pi)$ & $(0,\pi)$ & $a+b+c > 3\pi$ & H & $(\pi-a,2\pi-b,\pi-c)$ & $(-,-,-)$ & Figure~\ref{fig:small_angles_h2} (1) \\
			\hline
			$(0,\pi)$ & $(0,\pi)$ & $(\pi,2\pi)$ & $a+b+c > 3\pi$ & H & $(\pi-a,\pi-b,2\pi-c)$ & $(-,-,-)$ & Figure~\ref{fig:small_angles_h2} (1) \\
			\hline \hline
			$(\pi,2\pi)$ & $(0,\pi)$ & $(0,\pi)$ & $a-b-c > \pi$ & H & $(2\pi-a,b,c)$ & $(-,+,+)$ & Figure~\ref{fig:small_angles_h2} (2) \\
			\hline
			$(0,\pi)$ & $(\pi,2\pi)$ & $(0,\pi)$ & $-a+b-c > \pi$ & H & $(a,2\pi-b,c)$ & $(+,-,+)$ & Figure~\ref{fig:small_angles_h2} (2) \\
			\hline
			$(0,\pi)$ & $(0,\pi)$ & $(\pi,2\pi)$ & $-a-b+c > \pi$ & H & $(a,b,2\pi-c)$ & $(+,+,-)$ & Figure~\ref{fig:small_angles_h2} (2) \\
			\hline \hline
			$(\pi,2\pi)$ & $(0,\pi)$ & $(0,\pi)$ & $a-b+c < \pi$ & H & $(a-\pi,\pi-b,c)$ & $(+,-,+)$ & Figure~\ref{fig:small_angles_h2} (3i) \\
			\hline
			$(0,\pi)$ & $(\pi,2\pi)$ & $(0,\pi)$ & $a+b-c < \pi$ & H & $(a,b-\pi,\pi-c)$ & $(+,+,-)$ & Figure~\ref{fig:small_angles_h2} (3i) \\
			\hline
			$(0,\pi)$ & $(0,\pi)$ & $(\pi,2\pi)$ & $-a+b+c < \pi$ & H & $(\pi-a,b,c-\pi)$ & $(-,+,+)$ & Figure~\ref{fig:small_angles_h2} (3i) \\
			\hline \hline
			$(\pi,2\pi)$ & $(0,\pi)$ & $(0,\pi)$ & $a+b-c < \pi$ & H & $(a-\pi,b,\pi-c)$ & $(+,+,-)$ & Figure~\ref{fig:small_angles_h2} (3ii) \\
			\hline
			$(0,\pi)$ & $(\pi,2\pi)$ & $(0,\pi)$ & $-a+b+c < \pi$ & H & $(\pi-a,b-\pi,c)$ & $(-,+,+)$ & Figure~\ref{fig:small_angles_h2} (3ii) \\
			\hline
			$(0,\pi)$ & $(0,\pi)$ & $(\pi,2\pi)$ & $a-b+c < \pi$ & H & $(a,\pi-b,c-\pi)$ & $(+,-,+)$ & Figure~\ref{fig:small_angles_h2} (3ii) \\
			\hline
		\end{tabular}
	\end{center}
	\caption{Table of atomic triangular immersions of hyperbolic type.}
	\label{tab:hyperbolic}
\end{table}

\begin{table}[h]
	\begin{center}
  \hspace*{-1.85cm}
		\begin{tabular}{|c|c|c|c|c|c|c|c|}
			\hline
			\multicolumn{3}{|c|}{Angles Range} &  &  & Target angles & & \\
			\hline
			$a$ & $b$ & $c$ & Conditions & Type & $(\hat a, \hat b, \hat c)$ & Signs & Figure \\
			\hline \hline
			$(0,\pi)$ & $(0,\pi)$ & $(0,\pi)$ & $a+b+c > \pi$ & S & $(a,b,c)$ & $(+,+,+)$ & Figure~\ref{fig:tiny_angles_hs} (right) \\
			& & & $a+\pi > b+c$ & & & & \\
			& & & $b+\pi > a+c$ & & & & \\
			& & & $c+\pi > a+b$ & & & & \\
			\hline \hline
			$(\pi,2\pi)$ & $(0,\pi)$ & $(0,\pi)$ & $3\pi > a+b+c $ & S & $(2\pi-a,\pi-b,\pi-c)$ & $(-,-,-)$ & Figure~\ref{fig:small_angles_hs} (right) \\
			& & & $a+b > \pi+c$ & & & & \\
			& & & $a+c > \pi+b$ & & & & \\
			& & & $\pi > a-b-c$ & & & & \\
			\hline
			$(0,\pi)$ & $(\pi,2\pi)$ & $(0,\pi)$ & $3\pi > a+b+c $ & S & $(\pi-a,2\pi-b,\pi-c)$ & $(-,-,-)$ & Figure~\ref{fig:small_angles_hs} (right) \\
			& & & $a+b > \pi+c$ & & & & \\
			& & & $b+c > \pi+a$ & & & & \\
			& & & $\pi > -a+b-c$ & & & & \\
			\hline
			$(0,\pi)$ & $(0,\pi)$ & $(\pi,2\pi)$ & $3\pi > a+b+c $ & S & $(\pi-a,\pi-b,2\pi-c)$ & $(-,-,-)$ & Figure~\ref{fig:small_angles_hs} (right) \\
			& & & $b+c > \pi+a$ & & & & \\
			& & & $a+c > \pi+b$ & & & & \\
			& & & $\pi > -a-b+c$ & & & & \\
			\hline
		\end{tabular}
	\end{center}
	\caption{Table of atomic triangular immersions of spherical type.}
	\label{tab:spherical}
\end{table}

\begin{table}[h]
	\begin{center}
  \hspace*{-2.15cm}
		\begin{tabular}{|c|c|c|c|c|c|c|c|}
			\hline
			\multicolumn{3}{|c|}{Angles Range} &  &  & Target angles & & \\
			\hline
			$a$ & $b$ & $c$ & Conditions & Type & $(\hat a, \hat b, \hat c)$ & Signs & Figure \\
			\hline \hline
			$(0,\pi)$ & $(0,\pi)$ & $(0,\pi)$ & $a+b+c = \pi$ & E & $(a,b,c)$ & $(+,+,+)$ & Figure~\ref{fig:tiny_angles_e} (left) \\
			\hline \hline
			$(0,\pi)$ & $(0,\pi)$ & $(0,\pi)$ & $-a+b+c = \pi$ & E & $(a,\pi-c,\pi-b)$ & $(-,+,+)^*$ & Figure~\ref{fig:tiny_angles_e} (right) \\
			\hline
			$(0,\pi)$ & $(0,\pi)$ & $(0,\pi)$ & $a-b+c = \pi$ & E & $(\pi-a,\pi-c,b)$ & $(+,-,+)^*$ & Figure~\ref{fig:tiny_angles_e} (right) \\
			\hline
			$(0,\pi)$ & $(0,\pi)$ & $(0,\pi)$ & $a+b-c = \pi$ & E & $(\pi-a,c,\pi-b)$ & $(+,+,-)^*$ & Figure~\ref{fig:tiny_angles_e} (right) \\
			\hline \hline
			$(\pi,2\pi)$ & $(0,\pi)$ & $(0,\pi)$ & $a+b+c = 3\pi$ & E & $(2\pi-a,\pi-c,\pi-b)$ & $(+,+,+)^*$ & Figure~\ref{fig:small_angles_e2} (1) \\
			\hline
			$(0,\pi)$ & $(\pi,2\pi)$ & $(0,\pi)$ & $a+b+c = 3\pi$ & E & $(\pi-a,\pi-c,2\pi-b)$ & $(+,+,+)^*$ & Figure~\ref{fig:small_angles_e2} (1) \\
			\hline
			$(0,\pi)$ & $(0,\pi)$ & $(\pi,2\pi)$ & $a+b+c = 3\pi$ & E & $(\pi-a,2\pi-c,\pi-b)$ & $(+,+,+)^*$ & Figure~\ref{fig:small_angles_e2} (1) \\
			\hline \hline
			$(\pi,2\pi)$ & $(0,\pi)$ & $(0,\pi)$ & $a-b-c = \pi$ & E & $(2\pi-a,c,b)$ & $(+,-,-)^*$ & Figure~\ref{fig:small_angles_e2} (2) \\
			\hline
			$(0,\pi)$ & $(\pi,2\pi)$ & $(0,\pi)$ & $-a+b-c = \pi$ & E & $(a,c,2\pi-b)$ & $(-,+,-)^*$ & Figure~\ref{fig:small_angles_e2} (2) \\
			\hline
			$(0,\pi)$ & $(0,\pi)$ & $(\pi,2\pi)$ & $-a-b+c = \pi$ & E & $(a,2\pi-c,b)$ & $(-,-,+)^*$ & Figure~\ref{fig:small_angles_e2} (2) \\
			\hline \hline
			$(\pi,2\pi)$ & $(0,\pi)$ & $(0,\pi)$ & $a-b+c = \pi$ & E & $(a-\pi,\pi-b,c)$ & $(+,-,+)$ & Figure~\ref{fig:small_angles_e2} (3i) \\
			\hline
			$(0,\pi)$ & $(\pi,2\pi)$ & $(0,\pi)$ & $a+b-c = \pi$ & E & $(a,b-\pi,\pi-c)$ & $(+,+,-)$ & Figure~\ref{fig:small_angles_e2} (3i) \\
			\hline
			$(0,\pi)$ & $(0,\pi)$ & $(\pi,2\pi)$ & $-a+b+c = \pi$ & E & $(\pi-a,b,c-\pi)$ & $(-,+,+)$ & Figure~\ref{fig:small_angles_e2} (3i) \\
			\hline \hline
			$(\pi,2\pi)$ & $(0,\pi)$ & $(0,\pi)$ & $a+b-c = \pi$ & E & $(a-\pi,b,\pi-c)$ & $(+,+,-)$ & Figure~\ref{fig:small_angles_e2} (3ii) \\
			\hline
			$(0,\pi)$ & $(\pi,2\pi)$ & $(0,\pi)$ & $-a+b+c = \pi$ & E & $(\pi -a,b-\pi,c)$ & $(-,+,+)$ & Figure~\ref{fig:small_angles_e2} (3ii) \\
			\hline
			$(0,\pi)$ & $(0,\pi)$ & $(\pi,2\pi)$ & $a-b+c = \pi$ & E & $(a,\pi-b,c-\pi)$ & $(+,-,+)$ & Figure~\ref{fig:small_angles_e2} (3ii) \\
			\hline \hline
			$(2\pi,3\pi)$ & $(0,\pi)$ & $(0,\pi)$ & $a-b-c = \pi$ & E & $(a-2\pi,\pi-b,\pi-c)$ & $(+,-,-)$ & Figure~\ref{fig:additional_euclidean} \\
			\hline
			$(0,\pi)$ & $(2\pi,3\pi)$ & $(0,\pi)$ & $-a+b-c = \pi$ & E & $(\pi-a,b-2\pi,\pi-c)$ & $(-,+,-)$ & Figure~\ref{fig:additional_euclidean} \\
			\hline
			$(0,\pi)$ & $(0,\pi)$ & $(2\pi,3\pi)$ & $-a-b+c = \pi$ & E & $(\pi-a,\pi-b,c-2\pi)$ & $(-,-,+)$ & Figure~\ref{fig:additional_euclidean} \\
			\hline
		\end{tabular}
	\end{center}
	\caption{Table of atomic triangular immersions of Euclidean type.}
	\label{tab:Euclidean}
\end{table}

%%%%%%%%%%%%%%%%%%%%%%%%%%%%%%%%%
%%%%%%%%%%%%%%%%%%%%%%%%%%%%%%%%%
%%%%%%%%%%%%%%%%%%%%%%%%%%%%%%%%%

\clearpage 

\printbibliography

\end{document}